\documentclass[a4paper]{amsart} 
\usepackage{amssymb,times, amscd,amsmath,amsthm, xypic}
\author{J\"org Sch\"urmann and Shoji Yokura$^{(*)}$}
\address{J. Sch\"urmann: Westf. Wilhems-Universit\"at, SFB 478 ``Geometrische Strukturen in der Mathematik", Hittorfstr. 27, 
48149 M\"unster, Germany}
\email {jschuerm@math.uni-muenster.de}
\address{S. Yokura: Department of Mathematics and Computer Science, Faculty of Science, University of Kagoshima, 21-35 Korimoto 1-chome, Kagoshima 890-0065, Japan}
\email {yokura@sci.kagoshima-u.ac.jp}
\title{A survey of characteristic classes of singular spaces}
\thanks {(*) Partially supported by Grant-in-Aid for Scientific Research(No.17540088), the Japanese Ministry of Education, Science, Sports and Culture} 

\begin{document} 
\numberwithin{equation}{section}
\newtheorem{thm}[equation]{Theorem}
\newtheorem{pro}[equation]{Proposition}
\newtheorem{prob}[equation]{Problem}
\newtheorem{cor}[equation]{Corollary}
\newtheorem{lem}[equation]{Lemma}
\theoremstyle{definition}
\newtheorem{ex}[equation]{Example}
\newtheorem{defn}[equation]{Definition}
\newtheorem{rem}[equation]{Remark}
\renewcommand{\rmdefault}{ptm}
\def\alp{\alpha}
\def\be{\beta}
\def\jeden{1\hskip-3.5pt1}
\def\om{\omega}
\def\bigstar{\mathbf{\star}}
\def\ep{\epsilon}
\def\vep{\varepsilon}
\def\Om{\Omega}
\def\la{\lambda}
\def\La{\Lambda}
\def\si{\sigma}
\def\Si{\Sigma}
\def\Cal{\mathcal}
\def\ga{\gamma}
\def\Ga{\Gamma}
\def\de{\delta}
\def\De{\Delta}
\def\bF{\mathbb{F}}
\def\bH{\mathbb H}
\def\bPH{\mathbb {PH}}
\def \bB{\mathbb B}
\def \bK{\mathbb K}
\def \bG{\mathbf G}
\def \bL{\mathbf L}
\def\bM{\mathbf M}
\def\bN{\mathbb N}
\def\bR{\mathbb R}
\def\bP{\mathbb P}
\def\bZ{\mathbb Z}
\def\bC{\mathbb  C}
\def \bQ{\mathbb Q}
\def\op{\operatorname}

\begin{abstract} A theory of characteristic classes of vector bundles and smooth manifolds plays an important role in the theory of smooth manifolds. An investigation of reasonable notions of characteristic classes of singular spaces started since a systematic study of singular spaces such as singular algebraic varieties. We make a quick survey of characteristic classes of singular varieties, mainly focusing on the functorial aspects of some important ones such as the singular versions of the Chern class, the Todd class and the 
Thom--Hirzebruch's L-class. 
Then we explain our recent ``motivic" characteristic classes, which in a sense unify these three different theories of characteristic classes.
We also discuss bivariant versions of them and characteristic classes of proalgebraic varieties, which are related to the motivic measures/integrations.
Finally we explain some recent work on ``stringy" versions of these theories,
together with some references for ``equivariant" counterparts.
\end{abstract}

\maketitle

\centerline {\em Dedicated to Jean-Paul Brasselet on  the occasion of his 60th birthday}

\tableofcontents

\section{Introduction}\label{introduction}
Characteristic classes are usually certain kinds of cohomology classes for vector bundles over spaces and characteristic classes of smooth manifolds are defined via their tangent bundles. The most basic ones are {\em Stiefel--Whitney\/}, {\em Euler\/} and {\em Pontrjagin\/} classes  in the real case, and  {\em Chern\/} classes in the complex case. They were introduced in 1930's and 1940's  and constructed in a topological manner, i.e., via the obstuction theory, and in a differential-geometrical manner, i.e., via the Chern--Weil theory. Various important characteristic classes of vector bundles and invariants of manifolds are expressed as polynomials of them. The theory of cohomological characteristic classes were used for classifying manifolds and the study of structures of manifolds.\\ 

In 1960's a systematic study of singular spaces was started by {\em R. Thom\/}, {\em H. Whitney\/}, {\em H. Hironaka\/},
{\em S. \L ojasiewicz\/}, et al.; they studied triangulations, stratifications, resolution of singularities (in characteristic zero) and so on.
Already in 1958 {\em R. Thom\/} introduced in \cite{Thom}  {\em rational Pontrjagin and L-classes\/} for rational PL-homology manifolds.
In 1965 {\em M.-H. Schwartz\/} defined in \cite{Schwartz1} certain characteristic classes using obstruction theory of the so-called radial vector fields; the  {\em Schwartz class\/} is defined for a singular complex variety embedded in a complex manifold as a cohomology class of the manifold supported on the singular variety. In 1969, {\em D. Sullivan\/} \cite{Sullivan} proved that a real analytic space is mod 2 Euler space, i.e., the Euler--Poincar\'e characteristic of the link of any point is even, which implies that the sum of simplices in the first barycentric subdivision of any triangulation is mod 2 cycle. This enabled Sullivan to define the ``singular" {\em Stiefel--Whitney class\/} as a mod 2 homology class, which is equal to the Poincar\'e dual of the above cohomological Stiefel--Whitney class for a smooth variety. \\

{\em P. Deligne\/} and  {\em A.Grothendieck\/} (cf. \cite{Sullivan}) conjectured the unique existence of the  {\em Chern class\/} version of the Sullivan's Stiefel--Whitney class, and in 1974  {\em R. MacPherson\/} \cite{MacPherson1} proved their conjecture affirmatively. Motivated by MacPherson's proof of the conjecture,  {\em P. Baum\/},  {\em W. Fulton\/} and  
{\em R. MacPherson\/} \cite {Baum-Fulton-MacPherson} proved the so-called ``singular Riemann--Roch theorem", which is nothing but the  {\em Todd class\/} transformation in the case of singular varieties.\\

{\em M. Goresky\/} and  {\em R. MacPherson\/} (\cite {Goresky-MacPherson1}, \cite {Goresky-MacPherson2})
have introduced {\em Intersection Homology Theory\/}, by using the notion of ``perversity".  In \cite {Goresky-MacPherson1} they extended the work of \cite{Thom} to stratified spaces
with even (co)dimensional strata and introduced a {\em homology $L$-class\/} $L_*^{\op {GM}}(X)$ such that if $X$ is nonsingular 
it becomes the Poincar\'e dual of the original Thom--Hirzebruch $L$-class:
$L_*^{\op {GM}}(X) = L^{*}(TX) \cap [X].$ 
In \cite{Si} this was further extended to so-called stratified ``Witt-spaces",
whose intersection (co)homology complex (for the middle perversity) becomes self-dual
(compare also with \cite{Ban} for a more recent extension).
Later, {\em S. Cappell\/} and {\em J. Shaneson\/} \cite {Cappell-Shaneson1}(see also \cite {Cappell-Shaneson2} and \cite {Shaneson}) introduced 
a {\em homology $L$-class\/} transformation $L_*$, which turns out to be a natural 
transformation from the abelian group $\Omega(X)$ (see \S 7) of cobordism classes of selfdual constructible complexes to the rational homology group \cite{BSY2} (cf.\cite{Yokura-TAMS}).\\
 
In the case of singular varieties, the characteristic cohomology classes have been individually extended to the corresponding characteristic homology classes without any unifying theory of characteristic classes of singular varieties, unlike the case of smooth manifolds and vector bundles.
Only very recently such a unifying theory of ``motivic characteristic classes" for singular
spaces appeared in our work \cite{BSY2}.
The purpose of the present paper is to make a quick survey on the development of characteristic classes and the up date situation of characteristic classes of singular spaces. 
This includes our motivic characteristic classes, bivariant versions, characteristic classes of proalgebraic varieties and finally  ``stringy" versions of these theories,
together with some references for ``equivariant" counterparts.\\

The present survey is a kind of extended and up-dated version of MacPherson's survey article  \cite{MacPherson2} of more than 30 years ago. There are other surveys, e.g., \cite{Aluffi-lecture note}, \cite{Brasselet-lecture note}, 
\cite{Parusinski-lecture note}, \cite{Schuermann-lecture note}, \cite{Suwa-lecture note}  on characteristic classes of singular varieties written from different viewpoints.\\
$ $\\
{\it Acknowledgements.} 
It is a pleasure to thank P. Aluffi, J.-P. Brasselet, A. Libgober, P. Pragacz, J. Seade, T. Suwa, W. Veys and A. Weber
for valuable conversations about different aspects of this subject. 

This survey is a combined, modified and extended version of the author's two talks at ``Singularities in Geometry and Topology" (the 5th week of Ecole de la Formation Permanente du CNRS - Session r\'esidentielle de la FRUMAM)
held at Luminy, Marseille, during the period of 21 February -- 25 February 2005. The authors would like to thank the organizers of the conference for inviting us to give these talks. The second named author (S.Y.) also would like to thank the staff of ESI (Erwin Schr\"odinger International Institute for Mathematical Physics, Vienna, Austria), where a part of the paper was written in August 2005, for providing a nice atmosphere in which to work.  

\section{Euler--Poincar\'e characteristic}\label{EP}
The simplest, but most fundamental and most important topological invariant of a compact topological space is the {\em Euler number\/} or {\em Euler--Poincar\'e characteristic\/}. Its definition is quite simple; for a compact triangulable space or more generally for a cellular decomposable space $X$, it is defined to be the alternating sum of the numbers of cells and denoted by $\chi(X)$:
\begin{equation}
\chi(X) = \sum_n \; (-1)^i \sharp (n-\op{cells}). \label{chi-1}
\end{equation}
By the homology theory, the Euler--Poincar\'e characteristic turns out to be equal to the alternating sum of Betti numbers, i.e.,
\begin{equation}
\chi(X) = \sum_n \;(-1)^n \op {dim} H_n(X; \bR). \label{chi-2}
\end{equation}
With this fact, the Euler--Poincar\'e characteristic is defined for any topological space as long as the right-hand-side of 
(\ref{chi-2}) is defined, e.g. for locally compact semialgebraic sets. Note that taking the alternating sum is essential in the definition (\ref{chi-1}), but it is not the case in the definition (\ref{chi-2}). The following general form is called the 
{\em Poincar\'e polynomial\/}:
$$P_t(X):= \sum_n \op {dim} H_n(X; \bR) t^n ,$$
which is also a topological invariant. The Euler--Poincar\'e characteristic has the following properties:
\begin{enumerate}
\item $\chi(X) = \chi(X')$ if $X \cong X'$,
\item $\chi (X) = \chi(X,Y) + \chi (Y)$ for any closed  subspace $Y \subset X$, where the relative Euler--Poincar\'e characteristic 
$\chi(X,Y)$ is defined by the relative homology groups $H_*(X,Y)$,
\item $\chi(X \times Y) = \chi(X) \cdot \chi(Y)$.
\end{enumerate}
For a fiber bundle $f: X \to Y$ we have $\chi(X) = \chi(F) \cdot \chi(Y)$, if the Euler characteristic $\chi(F)$ of all fibers $F$ is constant, e.g. $Y$ is connected.
This generalizes the above property (3).  
The same properties also hold for the {\em Euler characteristic with compact support\/}
\begin{equation}
\chi_{c}(X) := \sum_n \;(-1)^n \op {dim} H^{n}_{c}(X; \bR), \label{chi-c1}
\end{equation}
together with the following {\em additivity} property
\begin{equation}
\chi_{c}(X)- \chi_{c}(Y) = \chi_{c}(X,Y) = \chi_{c}(X \setminus Y)  \label{chi-c2}
\end{equation} 
for any closed  subspace $Y \subset X$, where the relative Euler characteristic with compact support $\chi_{c}(X,Y)$ is defined by the relative cohomology groups $H^{*}_{c}(X,Y)
=H^{*}_{c}(X\setminus Y)$. Of course $\chi(X)=\chi_{c}(X)$ for $X$ compact.

\begin{rem} For two topological spaces $X, Y$, let $X + Y$ denote the topological sum, which is the disjoint sum, we clearly have
$$\chi(X + Y) = \chi(X) + \chi(Y).$$
However, we should note that for a closed subspace $Y \subset X$ the following 
additivity property does not hold in general:
\begin{equation}
\chi(X) = \chi(X \setminus Y) + \chi(Y), \label{chi-3}
\end{equation}
although $X = (X \setminus Y) + Y$ as a set, since the topological sum $Y + (X \setminus Y)$ is not equal to the original topological space $X$. In other words, $\chi(X, Y) \not = \chi(X \setminus Y)$ in general.

However, in the category of complex algebraic varieties, the above formula (\ref{chi-3}) holds, i.e., for any closed subvariety $Y \subset X$ we have that $\chi(X) = \chi(X \setminus Y) + \chi(Y)$. The key geometric reason for the equality 
$\chi(X) = \chi(X \setminus Y) + \chi(Y)$ is that a closed subvariety $Y$ always has a neighborhood deformation retract $N$ such that the Euler--Poincar\'e characteristic of the ``link" $\chi(N \setminus Y)$ vanishes due to a result of Sullivan (see 
\cite[Exercise, p.95, comments on p.141-142]{Fulton-toric}). In other words 
$\chi(X \setminus Y)=\chi_{c}(X \setminus Y)$ in the complex algebraic context, which also can be
extended and proved in the language of complex algebraically constructible functions
(see \cite[\S 6.0.6]{Schuermann-book}).
\end{rem}

\section{Characteristic classes of vector bundles}\label{characteristic class}
Very nice references for this section are the books
\cite{Milnor-Stasheff, Hirzebruch2, Hus, Stong}.
A characteristic class of vector bundles over a topological space $X$ is defined to be a map from the set of isomorphism classes of vector bundles over $X$ to the cohomology group (ring) $H^*(X; \Lambda)$ with a coefficient ring $\Lambda$, which is supposed to be compatible with the pullback of vector bundle and cohomology group for a continuous map. Namely, it is an assignment 
$c \ell: \op {Vect} (X) \to H^*(X; \Lambda)$ such that the following diagram commutes for a continuous map $f: X \to Y$:
$$\CD\op {Vect}(Y)@> {cl} >> H^*(Y; \Lambda)\\@V {f^*}VV @VV {f^*}V\\\op {Vect}(X)@>> {cl} > H^*(X; \Lambda). \endCD$$
Here $\op {Vect}(W)$ is the set of isomorphism classes of vector bundles over $W$. \\

The theory of characteristic classes started in Stiefel's paper \cite{Stiefel}, in which he considered the problem of the existence of tangential frames, i.e., linearly independent vector fields on a differentiable manifold. And at the same year H. Whitney defined such characteristic classes for sphere bundles over a simplicial complex \cite{Whitney}, and some time later he
invented cohomology and proved his important ``sum formula" \cite{Whitney2}. Then Pontrjagin \cite{Pontr} introduced other characteristic classes of real vector bundles, based on the study of the homology of real Grassmann manifolds.
Finally Chern \cite{Chern1, Chern2} defined similar characteristic classes of complex vector bundles.\\

The most fundamental characteristic classes of a real vector bundle $E$ over $X$ are the
{\em Stiefel-Whitney classes} $w^{i}(E)\in H^i(X; \bZ_2)$, 
{\em Pontrjagin classes} $p^{i}(E)\in H^{4i}(X; \bZ[1/2])$, and for a complex vector bundle
$E$ the {\em Chern classes} $c^{i}(E)\in H^{2i}(X; \bZ)$.
These characteristic classes $c \ell^{i}(E)\in H^*(X; \Lambda)$ are described axiomatically in a unified way:
\begin{defn}
The Stiefel Whitney resp. Pontrjagin classes of real vector bundles,
resp. Chern classes of complex vector bundles, is the operator assigning to each real
(resp. complex) vector bundle  $E \to X$ cohomology classes  
$$c \ell^{i}(E) := \begin{cases}
w^{i}(E) &\in H^i(X; \bZ_2)\\
p^{i}(E) &\in H^{4i}(X; \bZ[1/2])\\
c^i(E) &\in H^{2i}(X; \bZ)  
\end{cases}$$
of the base space  $X$ such that the following four axioms are satisfied:\\
\noindent {\bf Axiom-1}: (finiteness) For each vector bundle  $E$ one has $c \ell^{0}(E):= 1$ and 
$c \ell^{i}(E) = 0$  for $i > \op{rank} E$
(in fact $p^{i}(E)=0$ for $i>[\op {rank} E/2]$). $c \ell^{*}(E) := \sum _i\; c \ell^{i}(E)$ is called the corresponding {\em total characteristic class}.
In particular $c \ell^{*}(0_{X})=1$ for the zero vector bundle $0_{X}$ of rank zero.\\
\noindent {\bf Axiom-2}: (naturality)  
One has $ c \ell^{*}(F) = c \ell^{*}(f^{*}E) = f^*c \ell^{*}(E)$
for any cartesian diagram
$$\CD F\simeq f^{*}E @> {}>> E\\@V{}VV  @VV{}V\\Y@>>f> X \:.\endCD$$
\noindent {\bf Axiom-3}: (Whitney sum formula) 
$$c \ell^{*}(E \oplus F) = c \ell^{*}(E)c \ell^{*}(F)\:,$$
or more generally
$$c \ell^{*}(E) = c \ell^{*}(E')c \ell^{*}(E'')$$
for any short exact sequence $0\to E'\to E\to E''\to 0$ of vector bundles.\\
\noindent {\bf Axiom-4}: (normalization or the ``projective space" condition) For the canonical (i.e., the dual of the tautological) line bundle $\ga_n^1(\bK):= \Cal O_{\bold P^n(\bK)}(1)$  over the projective space $\bold P^n(\bK)$ (with $\bK=\bR, \bC$) one has:
\begin{enumerate}
\item[($w^1$):] $w^1(\ga_n^1(\bR))$ is non-zero.
\item[($p^1$):] $p^1(\ga_n^1(\bC))= c^1(\ga_n^1(\bC))^{2}$.
\item[($c^1$):] $c^1(\ga_n^1(\bC)) = [\bold P^{n-1}(\bC)]  \in H^2(\bold P^n(\bC); \bZ)$ 
is the cohomology class represented by the hyperplane $\bold P^{n-1}(\bC)$.
\end{enumerate}
\end{defn}

\begin{rem}
We use the superscript notation $c \ell^{*}$  for contravariant functorial characteristic classes of vector bundles in cohomology, to distinguish them from the 
subscript notation $c \ell_{*}$  for covariant functorial characteristic classes of singular spaces in homology, which we consider later on. Also note that in topology any short exact sequence
of vector bundles over a reasonable (i.e. paracompact) space splits (by using a metric on $E$).
But this is not the case in the algebraic or complex analytic context, where one should ask
the ``Whitney sum formula" for short exact sequences. 
\end{rem}

The existence of such a class for vector bundles of rank $n$ can be shown, for example, with the help of a classifying space, i.e., the infinite dimensional Grassmanian manifolds
$\bold G_{n}(\bK^{\infty})$ (with $\bK=\bR, \bC$),
and the fact that the cohomology ring of this Grassmanian manifold is a polynomial ring 
$$H^{*}(\bold G_{n}(\bK^{\infty});\Lambda)=
\begin{cases}
\bZ_{2}[w^1, w^2, \cdots, w^n] & \text{for $\bK=\bR$ and $\Lambda= \bZ_{2}$,}\\
\bZ[1/2][p^1, p^2, \cdots, p^{[n/2]}] & \text{for $\bK=\bR$ and $\Lambda= \bZ[1/2]$,}\\
\bZ[c^1, c^2, \cdots, c^n] & \text{for $\bK=\bC$ and $\Lambda= \bZ$.}
\end{cases}$$

 The most important axiom is Axiom-2 and the uniqueness of such a class follows then form 
 Axiom-3 and Axiom-4. By the so-called ``splitting principle" one can assume 
 (after pulling back to a suitable bundle, whose pullback on the cohomology level is injective)
 that a given non-zero vector bundle $E$ splits into a sum of line (or $2$-plane) bundles.
 These line (or $2$-plane) bundles are then called the ``Chern roots" of $E$.
 Then Axiom-3 reduces the calculation of characteristic classes to the case of line bundles (for $c \ell=w,c$) or real 2-plane bundles
(for $c \ell=p$). By naturality these are then fixed by Axiom-4, since
$$\bold G_{1}(\bK^{\infty})= \lim_{k} \:\bold P^k(\bK) \quad \text{(for $\bK=\bR, \bC$),}$$
for the case $c \ell=w,c$, or from the fact that the canonical projection 
$$\lim_{k} \:\bold P^k(\bC) \to \bold G_{2}(\bR^{\infty})$$
is the orientation double cover for the case $c \ell=p$.\\ 

From the axioms one gets in all cases $w^1,p^1$ and $c^1$ are {\em nilpotent} on finite
dimensional spaces, and $c \ell^{*}(E)=1$ for a trivial vector bundle
$E$. Note that a real oriented line bundle is always trivial so that a real line
bundle $L\to X$ has no interesting characteristic class $c \ell^{j}(L)= 0\in  H^{j}(X; \bZ[1/2])$
for $j>0$. Just pullback to an orientation double cover $\pi: \tilde{X}\to X$ so that
$\pi^{*}L$ is orientable with $\pi^{*}: H^{j}(X; \bZ[1/2])\to H^{j}(\tilde{X}; \bZ[1/2])$
injective (since $2\in \bZ[1/2]$ is invertible). In particular a real vector bundle $E$ of rank
$r$ is orientable if and only if $w^1(E)=w^1(\Lambda^r\;E)=0$.\\

If a characteristic class $c \ell^{*}: \op {Vect} (X) \to H^*(X; \Lambda)$ satisfies the Whitney sum condition
$$c \ell^{*}(E \oplus F) = c \ell^{*}(E) c \ell^{*}(F)
\quad \text{with} \quad  c \ell^{*}(0_{X})=1 \:,$$
then $c \ell^{*}$ is called a {\em multiplicative\/} characteristic class.
Another important multiplicative characteristic class of an {\em oriented\/} real vector bundle
$E\to X$ of rank $r$ is the {\em Euler class} $e(E)\in H^{r}(X;\bZ)$, with
$e(E)\; mod\;2= w^{r}(E)$, $e(E)^2= p^{r/2}(E)$ for $r$ even 
 and $e(E)= c^{r}(E)$ in case $E$ is given by a complex vector bundle $E$
of rank $r$. But the {\em Euler class} is not a {\em normalized} characteristic class with
 $c \ell^{0}(L)= 1$.\\

The {\em Stiefel-Whitney, Pontrjagin and Chern classes} are essential in the sense that any {\em multiplicative} characteristic class $c \ell^{*}$ over finite dimensional base spaces is uniquely expressed as a polynomial (or power series) in these classes, i.e. the ``splitting principle" implies: 

\begin{thm}\label{charpower}
Let $\Lambda$ be a $\bZ_{2}$-algebra (resp. a $\bZ[1/2]$-algebra) for the case of real
vector bundles, or a $\bZ$-algebra for the case of complex vector bundles.
Then there is a one-to-one correspondence between 
\begin{enumerate}
\item {\em multiplicative} characteristic classes $c \ell^{*}$ over finite dimensional base spaces,  and 
\item {\em formal power series} $f\in \Lambda[[z]]$
\end{enumerate} 
such that
$c \ell^{*}(L)=f(w^1(L))$ or $c \ell^{*}(L)=f(c^1(L))$ for any real or complex line bundle $L$
(resp. $c \ell^{*}(L)=f(p^1(L))$ for any real $2$-plane bundle $L$). In this case $f$ is called
the {\em characteristic power series\/} of the corresponding multiplicative characteristic class $c \ell^{*}_{f}$.
\end{thm}

\begin{rem} 
For the result above it is important that characteristic classes of vector bundles
live in cohomology so that one can build new classes by multiplication (i.e. by the cup-product)
of the basic ones. This is not possible in the case of characteristic classes of singular spaces, which live in homology (except in the case of homology manifolds where Poincar\'{e}
duality is available).
\end{rem}

Moreover $c \ell^{*}_{f}$ is invertible with inverse $c \ell^{*}_{1/f}$, if $f\in \Lambda[[z]]$
is invertible, i.e. if $f(0)\in \Lambda$ is a unit (e.g. $f$ is a normalized power series
with $f(0)=1$).
Then the corresponding multiplicative characteristic class $c \ell^{*}$  extends 
over finite dimensional base spaces $X$  to a natural transformation of groups
$$c \ell^{*}: (\bold K(X),\oplus)\to (H^{*}(X;\Lambda),\cup)$$ 
on the Grothendieck group $\bold K(X)$
of real or complex vector bundles over $X$.\\

\section{Characteristic classes of smooth manifolds}\label{manifolds}
Let us now switch to smooth manifolds, which will be an important intermediate step
on the way to characteristic classes of singular spaces.
For a smooth (or almost complex) manifold $M$ its real (or complex) tangent bundle $TM$ is available and a characteristic
class $cl^{*}(TM)$ of the tangent bundle $TM$ is called a {\em characteristic cohomology class} $cl^{*}(M)$ of the manifold $M$. We also use the notation 
$$cl_{*}(M):=cl^{*}(TM)\cap [M]\in H^{BM}_{*}(M;\Lambda)$$
for the corresponding {\em characteristic homology class} of the manifold $M$,
with $[M]\in H^{BM}_{*}(M;\Lambda)$ the fundamental class in Borel-Moore homology
of the (oriented) manifold $M$. Note that $H^{BM}_{*}(X;\Lambda)=H_{*}(X;\Lambda)$
for $X$ compact.

\begin{rem} Using a relation to suitable cohomology operations, i.e. Steenrod squares,
Thom \cite{Thom1} has shown that the Stiefel-Whitney classes $w^{*}(M)$ of a smooth manifold $M$
are {\em toplogical\/} invariants. 
Later he introduced in \cite{Thom}  {\em rational Pontrjagin and L-classes\/} for compact rational PL-homology manifolds so that the {\em rational\/} Pontrjagin classes $p^{*}(M)\in H^{*}(M;\bQ)$ of a closed smooth manifold $M$ are {\em combinatorial or piecewise linear\/} invariants. 
A deep result of Novikov \cite{Nov} implies the {\em topological\/} invariance of these
{\em rational\/} Pontrjagin classes $p^{*}(M)\in H^{*}(M;\bQ)$ of a smooth manifold $M$.
\end{rem}

For a {\em closed oriented\/} manifold $M$ one has the interesting formula
\begin{equation} \label{Eulerclass}
deg(e(M))=\int_M \:e(TM) \cap [M] = \chi(M) \:,
\end{equation}
which justifies the name ``Euler class".
For a closed complex manifold $M$ this formula becomes
$$deg(c_{*}(M))=\int_M \;c^{*}(TM) \cap [M] = \chi(M) \:,$$
which is called the {\em Gauss--Bonnet--Chern Theorem}
(compare \cite{Chern3}). In this sense, the Chern class is a higher cohomology class version of the Euler--Poincar\'e characteristic. Similarly
$$deg(w_{*}(M))=\int_M \;w^{*}(TM) \cap [M] = \chi(M)\;mod\;2 $$
for any closed manifold $M$.\\

More generally let $Iso(n-dim.\; G-mfd.)$ be the set of isomorphism classes of smooth closed
(and oriented) pure $n$-dimensional manifolds $M$ for $G=O$ (or $G=SO$), 
or of pure $n$-dimensional weakly (``$=$ stably") almost complex manifolds $M$ for $G=U$,
i.e. $TM\oplus \bR_{M}^N$ is a complex vector bundle (for suitable $N$, with $\bR_{M}$ the trivial real line bundle over $M$). Then
$$Iso(G-mfd.)_{*}:= \bigoplus_{n}\: Iso(n-dim.\; G-mfd.)$$
becomes a commutative graded semiring with addition and multiplication given by disjoint union
and exterior product, with $0$ and $1$ given by the classes of the empty set and one point
space. Moreover any multiplicative characteristic class $c \ell_{f}$
coming from the power series $f$ in the variable $z=w^1, p^1$ or $c^1$ induces by
$$M\mapsto deg(c \ell_{f*}(M)):= \int_M \;c \ell^{*}_{f}(TM) \cap [M]$$
a semiring homomorphism 
$$\Phi_{f}: Iso(G-mfd.)_{*}\to \Lambda =
\begin{cases}
\text{a $\bZ_{2}$-algebra for $G=O$ and $z=w^1$,}\\
\text{a $\bZ[1/2]$-algebra for $G=SO$ and $z=p^1$,}\\
\text{a $\bZ$-algebra for $G=U$ and $z=c^1$.}
\end{cases}$$

Let $\Omega^{G}_{*}:= Iso(G-mfd.)_{*}/\sim$ be the corresponding {\em cobordism ring}
of closed ($G=O$) and oriented ($G=SO$) or weakly (``$=$ stably") almost complex manifolds ($G=U$),
with $M\sim 0$ for a closed pure $n$-dimensional $G$-manifold $M$ if and only if 
there is a compact pure $n+1$-dimensional $G$-manifold $B$ with boundary $\partial B\simeq M$.
Note that this is indeed a ring with $-[M]=[M]$ for $G=O$ or $-[M]=[-M]$ for $G=SO,U$,
where $-M$ has the opposite orientation of $M$. Moreover, for $B$ as above with
$\partial B\simeq M$ one has 
$$TB|\partial B \simeq TM\oplus \bR_{M}$$
so that 
$c \ell^{*}_{f}(TM)=i^{*}c \ell^{*}_{f}(TB)$ for $i: M\simeq \partial B \to B$ the closed inclusion of
the boundary. This also explains the use of the stable tangent bundle for the definition of a
stably or weakly almost complex manifold.
By a simple argument due to Pontrjagin one then gets 
$$M\sim 0 \quad \Rightarrow  \quad deg(c \ell_{f*}(TM))=\int_M \;c \ell^{*}_{f}(TM) \cap [M]=0$$
so that any multiplicative characteristic class $c \ell^{*}_{f}$
coming from the power series $f$ in the variable $z=w^1, p^1$ or $c^1$ induces
a ring homomorphism called {\em genus\/}
\begin{equation} \label{genus}
\Phi_{f}: \Omega^{G}_{*}\to \Lambda =
\begin{cases}
\text{a $\bZ_{2}$-algebra for $G=O$ and $z=w^1$,}\\
\text{a $\bZ[1/2]$-algebra for $G=SO$ and $z=p^1$,}\\
\text{a $\bZ$-algebra for $G=U$ and $z=c^1$.}
\end{cases}
\end{equation}
In fact for $\Lambda$ a $\bQ$-algebra this induces a one-to-one correspondence between
\begin{enumerate}
\item normalized power series $f$ in the variable $z=p^1$ (or $c^1$),
\item {\em normalized and multiplicative} characteristic classes $c \ell_{f}^{*}$ over finite dimensional base spaces,  and
\item genera $\Phi: \Omega^{G}_{*}\to \Lambda$ for $G=SO$ (or $G=U$). 
\end{enumerate}
Here one uses the following structure theorem.

\begin{thm}\label{cob-structure} 
\begin{enumerate}
\item {\em (Thom)} $\Omega^{SO}_{*}\otimes \bQ = \bQ[[\bold P^{2n}(\bC)]|n\in \bN]$
is a polynomial algebra in the classes of the complex even dimensional projective spaces.
\item {\em (Milnor)} $\Omega^{U}_{*}\otimes \bQ = \bQ[[\bold P^{n}(\bC)]|n\in \bN]$
is a polynomial algebra in the classes of the complex  projective spaces.
\end{enumerate}
\end{thm}

In particular, the corresponding genus $\Phi_{f}$ with values in a $\bQ$-algebra $\Lambda$,
or the corresponding normalized and multiplicative characteristic class $c \ell_{f}^{*}$,
is uniquely fixed by the values $\Phi_{f}(M)=\int_M \;c \ell^{*}_{f}(TM) \cap [M]$
for all (complex even dimensional) complex projective spaces $M=P^{n}(\bC)$.
These are best codified by the {\em logarithm\/} $g\in \Lambda[[t]]$ of $\Phi_{f}$:
\begin{equation} \label{log}
g(t):= \sum_{i=0}^{\infty} \; \Phi_{f}(P^{i}(\bC))\cdot \frac{t^{i+1}}{i+1} \:.
\end{equation}
Moreover, a genus $\Phi_{f}: \Omega^{U}_{*}\otimes \bQ \to \Lambda$ factorizes over the
canonical map 
$$\Omega^{U}_{*}\otimes \bQ \to\Omega^{SO}_{*}\otimes \bQ$$ 
if and only if
$f(z)$ is an even power series in $z=c^1$, $f(z)=h(z^2)$ with $z^2= (c^1)^2=p^1$.
Consider for example the {\em signature\/} $\sigma (M)$ of the cup-product pairing on the middle dimensional cohomology of the closed oriented manifold $M$ of real dimension $4n$,
with $\sigma (M):=0$ in all other dimensions. This defines a genus
$\sigma: \Omega^{SO}_{*}\otimes \bQ\to \bQ$, as observed by Thom, 
with $\sigma (P^{2n}(\bC))=1$ for all $n$.
The signature genus comes from the
normalized power series $h(z)=\sqrt z /\tanh(\sqrt z)$ in the variable $z=p^1$
(or $f(z)=z/\tanh(z)$ in the variable $z=c^1$), whose corresponding
characteristic class $c \ell^{*}=L^{*}$ is by definition the Hirzebruch-Thom $L$-class.
This is the content of the famous
{\em Hirzebruch's Signature Theorem\/} (compare also with \cite{Hirzebruch3}):
$$\sigma (M)  = \int _M L^{*}(TM) \cap [M].$$

\begin{rem} The first structure theorem about cobordism rings due to Thom is the description
of $\Omega^{O}_{*}$ as a polynomial algebra 
$\bZ_{2}[[M^n]|n\in \bN, n+1\neq 2^k]$
in the classes of suitable closed manifolds $M^n$ of dimension $n$,
with one generator in each dimension $n$ with $n+1$ not a power of $2$.
Then each genus $\Omega^{O}_{*}\to \Lambda$ to a $\bZ_{2}$-algebra $\Lambda$ is coming
form a normalized and multiplicative characteristic class $c \ell_{f}^{*}$,
but this correspondence is not injective.
\end{rem}

The value $\Phi(M)$ of a genus $\Phi$ on the closed manifold $M$ is also called a
characteristic number of $M$. All these numbers can be used to classify closed manifolds
up to cobordism.

\begin{thm} \begin{enumerate}
\item {\em (Pontrjagin--Thom)} Two closed $C^{\infty}$-manifolds are cobordant (i.e., represent the same element in $\Omega^{O}_{*}$) if and only if all their Stiefel--Whitney numbers are the same.
\item {\em (Thom--Wall)} Two closed oriented $C^{\infty}$-manifold are corbordant up to two-torsion (i.e., represent the same element in
$\Omega^{SO}_{*}\otimes \bZ[1/2]$) if and only if all their Pontrjagin numbers are the same.
\item {\em (Milnor--Novikov)} Two closed stably or weakly almost complex manifold are cobordant (i.e., represent the same element in
$\Omega^{U}_{*}$) if and only if all their Chern numbers are the same.
\end{enumerate}
\end{thm}

\section{Hirzebruch--Riemann--Roch and Grothendieck--Riemann--Roch}\label{HRR-GRR}
Let $X$ be a non-singular complex projective variety and $E$ a holomorphic vector bundle over $X$. 
Note that in this context we do not need to distinguish between holomorphic and algebraic vector bundles,
and similarly for coherent sheaves, by the so-called ``GAGA-principle" \cite{Serre-GAGA}.
Then the Euler--Poincar\'e characteristic of $E$ is defined by
$$ \chi (X,E) = \sum_{i \geq 0}(-1)^i \op {dim}_{\bC}H^i(X;\Om (E)),$$ 
where $\Om(E)$  is the coherent sheaf of germs of sections of  $E$. {\em J.-P. Serre\/} conjectured in his letter to 
{\em Kodaira and Spencer\/} (dated September 29, 1953) that there exists a polynomial $P(X, E)$ of Chern classes of the base variety $X$ and the vector bundle $E$ such that$$ \chi(X, E) = \int_X P(X,E) \cap [X].$$
Within three months (December 9, 1953) {\em F. Hirzebruch\/} solved this conjecture affirmatively: the above looked for polynomial 
$P(X,E)$ can be expressed as
$$P(X,E) = ch^{*}(E) td^{*}(X)$$
where $ch^{*}(E)$ is the total {\em Chern character\/} of $E$ and $td^{*}(TX)$ is the total {\em Todd class} of the tangent bundle $TX$ of  $X$ .  Let us recall that the cohomology classes $ch^{*}(V)$ and $td^{*}(V)$ are defined as follows:
$$ch^{*}(V) = \sum_{i=1}^{\op {rank} V} e^{\alp _i} \in H^{2*}(X;\bQ)$$
and
$$td^{*}(V) = \prod _{i=1}^{\op {rank} V} \frac{\alp _i}{1-e^{-\alp _i}} \in H^{2*}(X;\bQ)$$
where $\alp _i$'s are the Chern roots of  $V$ .
So $td^{*}$ is just the normalized and
multiplicative characteristic class corrsponding to the normalized power series $f(z)=z/(1-e^{-z})$ in $z=c^{1}$.
Similarly the  Chern character defines a contravariant natural transformation of rings
$$ch^{*}:  (\bold K(X),\oplus, \otimes)\to (H^{2*}(X;\bQ),+,\cup)$$ 
on the Grothendieck group $\bold K(X)$
of complex vector bundles over $X$.
Then we have the following celebrated theorem of Hirzebruch:

\begin{thm}(Hirzebruch--Riemann--Roch)
\begin{equation}
\chi (X,E) = T(X,E) := \int_X \left (ch^{*}(E) td^{*}(X) \right ) \cap [X]\:. 
\tag{{\bf HRR}}
\end{equation}
\end{thm}
$T(X,E)$ is called the $T$-characteristic (\cite{Hirzebruch2}). For a more detailed historical aspect of {\bf HRR}, 
see \cite{Hirzebruch3}.

\begin{rem} The $T$-characteristic $T(X,E)$ is {\em a priori} a rational number  by the definitions of the Todd 
class and Chern character, but it has to be {\em an integer} as a consequence of the HRR. The $T$-characteristic 
$T(X,E)$ of a complex vector bundle $E$ can be defined for any almost complex manifold and Hirzebruch 
\cite {Hirzebruch1} asked if the T-genus $T(X):= T(X, \jeden)$ with $\jeden$ being a trivial line bundle is always an integer. Of course this follows from HRR and the later result of Quillen,
that $\Omega^{U}_{*}\otimes \bQ$ is generated by complex projective algebraic manifolds.
The identity
$$z/(1-e^{-z}) = e^{z/2}\cdot z/(2\sinh{(z/2)})$$
allows one to introduce the Todd class
$$ Td^{*}(X) := e^{c^1(TX)/2} \cdot  \hat A^{*}(TX)\:,$$
and therefore also the $T$-characteristic $T(X,E)$,
more generally for a so-called {\em Spin$^c$-manifold} $X$.
Here  $\displaystyle \hat A$ is the so-called {\em A hat genus or characteristic class\/} corresponding to the even normalized power series
$f(z)=z/(2\sinh{(z/2)})$ in the variable $z=c^{1}$ or to
$f(z)= \sqrt{z}/(2\sinh{(\sqrt{z}/2)})$ in the variable $z=p^{1}$.
The $T$-characteristic $T(X,E)$ of a complex vector bundle $E$ is then an {\em integer\/}
by an application of the {\em Atiyah-Singer Index theorem\/} \cite{AS}
for a suitable {\em Dirac operator\/} (compare \cite[p.197, Theorem 26.1.1]{Hirzebruch1}).
\end{rem}

{\em A. Grothendieck\/} (cf. \cite{Borel-Serre}) generalized {\bf HRR} for non-singular quasi-projective algebraic varieties over any field and proper morphisms with Chow cohomology ring theory instead of ordinary cohomology theory
(compare also with \cite[chapter 15]{Fulton-intersection}). For the complex case we can still take the ordinary cohomology theory (or the homology theory by the Poincar\'e duality). Here we stick ourselves to {\em complex projective algebraic varieties\/} for the sake of simplicity. For a variety $X$, let $\bold G_0(X)$ denote the Grothendieck group of algebraic coherent sheaves on $X$ and for a morphism $f:X \to Y$ the pushforward $f_! :\bold G_0(X) \to \bold G_0(Y)$ is defined by
$$f_!(\Cal F) := \sum_{i \geq 0} (-1)^i \bold R ^if_* \Cal F ,$$
where $\bold R ^if_* \Cal F$ is (the class of ) the higher direct image sheaf of $\Cal F$. Then $\bold G_0$  is a covariant functor with the above pushforward (see \cite{Grothendieck1} and \cite{Manin}). 
Let similarly $\bold K^{0}(X)$ be the Grothendieck group of complex algebraic vector bundles over $X$
so that one has a canonical contravariant transformation of rings $\bold K^{0}(\quad)\to \bold K(\quad)$
to the Grothendieck group of complex vector bundles.
Note that on a smooth algebraic manifold the canonical map $\bold K^{0}(\quad) \to \bold G_0(\quad)$ taking the sheaf of sections is an isomorphism. With this isomorphism one can define characteristic classes of any algebraic coherent sheaf.
Then Grothendieck showed the existence of a natural transformation from the covariant functor $\bold G_0$ to the $\bQ$-homology covariant functor $H_{2*}(\quad;\bQ)$ (see \cite{Borel-Serre}):
\begin{thm}(Grothendieck--Riemann--Roch)  
Let the transformation $\tau_{*}:\bold G _0 (\quad) \to H_{2*}(\quad;\bQ)$ be defined by  
$\tau_{*}(\Cal F) = td^{*}(X)ch^{*}(\Cal F)\cap [X]$ for any smooth variety $X$. Then $\tau_{*}$ is actually natural, i.e., for any morphism $f:X \to Y$  the following diagram commutes:
$$\CD\bold G_0(X) @>\tau_{*} >> H_{2*}(X;\bQ)\\@V f_!VV  @VV f_* V\\\bold G_0(Y) @>>\tau > H_{2*}(Y;\bQ)\endCD$$
i.e.,                                      
\begin{equation}
 td^*(T_Y)ch^*(f_! \Cal F)\cap [Y] = f_*(td^*(TX)ch^*(\Cal F)\cap[X]).
\tag{{\bf GRR}}
\end{equation}
\end{thm}
Clearly  {\bf HRR}  is induced from  {\bf GRR} by considering a map from $X$to a point. Note that the target of the transformation of the original {\bf GRR} is the cohomology $H^{2*}(\quad;\bQ)$ with the Gysin homomorphism instead of the homology $H_{2*}(\quad;\bQ)$, but, by the definition of the Gysin homomorphism the original {\bf GRR} can be put in as above.

\section{The Generalized Hirzebruch--Riemann--Roch}\label{gHRR}
In Hirzebruch's book \cite[\S 12.1 and \S 15.5]{Hirzebruch2} he has generalized the characteristics $\chi (X,E)$ and $T(X,E)$ to the so-called {\em $\chi _y$-characteristic\/} $\chi _y(X,E)$ and {\em $T_y$-characteristic\/} $T_y(X,E)$ as follows, using a parameter $y$ (see also \cite [Chapter 5]{HBJ}).
\begin{defn}     
$$\aligned\chi_y(X,E) : 
& = \sum_{p \geq 0}\left (\sum_{q\geq 0}(-1)^q \op {dim}_{\bold C}H^q(X, \Om (E)\otimes \La ^p T^*X) \right )y^p\\         
& = \sum_{p \geq 0}\chi (X, E\otimes \La ^p T^*X))y^p \endaligned$$
where $T^*X$ is the dual of the tangent bundle $TX$, i.e., the cotangent bundle of $X$.
$$\aligned T_y(X,E) &:= \int _X \widetilde td_{(y)}(TX) ch_{(1+y)}(E) \cap [X],\\
\widetilde {td_{(y)}}(TX) &:= \prod _{i=1}^{\op {dim}X}\left (\frac {\alp _i(1+y)}{1-e^{-\alp _i(1+y)}} - \alp _iy \right ), \\
ch_{(1+y)}(E) &:= \sum_{j=1}^{\op {rank} ~E} e^{\be _j(1+y)}, \endaligned$$
where ${\alp _i}'s$ are the Chern roots of $TX$ and ${\be _j}'s$ are the Chern roots of $E$ .\end{defn}
F. Hirzebruch \cite[\S 21.3]{Hirzebruch2} showed the following theorem:
\begin{thm} (The generalized Hirzebruch--Riemann--Roch)
\begin{equation}
\chi _y(X,E) = T_y(X,E).
\tag{{\bf g-HRR}}
\end{equation}
\end{thm}
The {\bf g-HRR} can be shown as follows, using {\bf HRR}:
$$\aligned \chi_y(X,E) & = \int_X \sum_{p \geq 0}\chi (X, E\otimes \La ^p T^*X))y^p \quad \text {(by the definition)} \\
& = \int_X \sum_{p \geq 0}\left (ch^*(E\otimes \La ^p T^*X)td^*(X) \cap [X] \right) y^p \quad \text {(by {\bf HRR})}\\
& = \int_X \left (\sum_{p \geq 0}ch^*(E\otimes \La ^p T^*X)td^*(X) y^p \right ) \cap [X] \\
& = \int_X \left (ch^*(E)td^*(X) \sum_{p \geq 0}ch^*(\La ^p T^*X)y^p \right ) \cap [X] \\
& = \int_X \left (ch^*(E)td^*(X)\prod_{i =1}^{\op {dim}X}\left (1 + y e ^{- \alp_i} \right) \right ) \cap [X] \\
& = \int_X \left (\sum_{j=1}^{\op {rank} ~E} e^{\be _j}\prod_{i =1}^{\op {dim}X}\left (1 + y e ^{- \alp_i}\right) \frac {\alp_i}{1 - e^{-\alp_i}} \right ) \cap [X].  \endaligned$$
However, the power series $\displaystyle \left (1 + y e ^{- \alp_i}\right) \frac {\alp_i}{1 - e^{-\alp_i}}$ is not a normalized power series because the $0$-degree part is $1 + y$, not $1$. So, by dividing this non-normalized power series by $1 +y$ and furthermore by changing $\be_j$ to $\be_j(1+y)$ and $\alp_i$ to $\alp_i(1+y)$, which does not change the value of $\chi_y(X,E)$ at all, and by noticing that 
$$ \frac {1 + y e ^{- \alp_i(1+y)}}{1+y} \frac {\alp_i}{1 - e^{-\alp_i(1+y)}} = 
\frac {\alp _i(1+y)}{1-e^{-\alp _i(1+y)}} - \alp _iy,$$
we can see that the right hand side of the last equation is $T_y(X,E)$ (compare \cite[p.61-62]{HBJ}). In fact the same argument shows that a non-normalized power series
$f(z)$ with $a:=f(0)\in \Lambda$ a unit induces the same genus as the normalized power series
$f(az)/a$.

\begin{rem} \label{rem:AS}
The generalized Hirzebruch Riemann-Roch theorem is also true for a holomorphic vector bundle $E$ over a compact complex manifold $X$,
by an application of the {\em Atiyah-Singer Index theorem\/} \cite{AS}. 
\end{rem}

The above {\em modified Todd class\/} $\widetilde td_{(y)}$ is the normalized and multiplicative characteristic class
corresponding to the normalized power series (in $z=c^1$):
$$f(z)=f_{y}(z)= \frac {z(1+y)}{1-e^{-z(1+y)}} - zy \in \bQ[y][[z]]\:.$$
The associated genus $\chi_{y}: \Omega^{U}_{*}\to \bQ[y]$ is called the Hirzebruch {\em $\chi_{y}$-genus\/}.
A simple residue calculation in \cite[Lemma 1.8.1]{Hirzebruch2} implies for all $n\in \bN$:
\begin{equation}\label{chiy}
\chi_{y}(P^{n}(\bC))= \sum_{i=0}^{n}\; (-y)^i \in \bZ[y]\subset \bQ[y]\:.
\end{equation}
So these values on $P^{n}(\bC)$ fix the $\chi_{y}$-genus and the modified Todd class $\widetilde td_{(y)}$.
Moreover, the normalized power series $f_{y}(z)$ specializes to
$$f_{y}(z)= \begin{cases}
1+z \quad &\text{for $y=-1$},\\
z/(1-e^{-z}) \quad &\text{for $y=0$},\\
z/\tanh{(z)} \quad &\text{for $y=1$}.
\end{cases} $$
So the modified Todd class $\widetilde td_{(y)}$
defined above unifies the following three important characteristic cohomology classes:\\
(y = -1) \quad  the total Chern class    
$$\widetilde td_{(-1)}(TX) = c^{*}(TX),$$
(y = 0) \quad the total Todd class
$$\widetilde td_{(0)}(TX) = td^{*}(TX),$$
(y = 1) \quad the total Thom--Hirzebruch  L-class
$$\widetilde td_{(1)}(TX)  = L^{*}(TX).$$
Therefore, when $E$ = the trivial line bundle, for these special values $y = -1, 0, 1$ the {\bf g-HRR} reads as follows:\\
(y = -1) \quad  {\it Gauss--Bonnet--Chern Theorem}:
$$ \chi (X) = \int _X c^{*}(TX) \cap [X],$$
(y = 0) \quad {\it Riemann--Roch Theorem}: denoting $\chi_a(X):= \chi(X, \Cal O_X)$, called the arithmetic genus of $X$, to avoid the confusion with the above topological Euler--Poincar\'e characteristic $\chi (X)$,
$$\chi_a (X) = \int _X td^{*}(TX) \cap [X],$$
(y = 1) \quad {\it Hirzebruch's Signature Theorem}:
$$\sigma (X)  = \int _X L^{*}(TX) \cap [X].$$

\begin{rem} (Poincar\'e--Hopf Theorem) The above Gauss--Bonnet--Chern Theorem due to Chern \cite {Chern3} is a generalization of the original Gauss--Bonnet theorem saying that the integration of the Guassian curvature is equal to $2 \pi$ times the topological Euler--Poincar\'e characteristic. There is another well-known differential-topological formula concerning the topological Euler--Poincar\'e characteristic. That is the so-called {\it Poincar\'e --Hopf theorem}, saying that the index of a smooth vector field $V$ with only isolated singularites on a smooth compact manifold $M$ is equal to the topological Euler--Poincar\'e characteristic of the manifold 
$M$;
$$\op {Index} (V) = \chi (M),$$
where the index $\op {Index} (V)$ is defined to be the sum of the indices of the vector field at the isolated singularities. 
Compare with \cite{Milnor} for a beautiful introduction to the Poincar\'e --Hopf theorem.
Note that the Gauss--Bonnet--Chern Theorem follows from the Poincar\'e--Hopf theorem (cf. \cite{Witten} and \cite{Zhang}).
\end{rem}

\section{Characteristic classes of singular varieties}\label{singular characteristic classes}
In the following we consider for simplicity only compact spaces.
For a singular algebraic or analytic variety $X$ its tangent bundle is not available any longer because of the existence of singularities, thus one cannot define its characteristic class $c \ell_{*}(X)$ as in the previous case of manifolds, although a ``tangent-like"  bundle such as Zariski tangents is available. A main theme for defining reasonable characteristic classes for singular varieties is that reasonable ones should be interesting enough; for example, they should be geometrically or topologically interesting and quite well related to other well-known invariants of varieties and singularities (e.g., see \cite{MacPherson2}).\\

The theory of characteristic classes of vector bundles is {\em a natural transformation\/} from the contravariant functor $\op {Vect}$ to the contravariant cohomology functor $H^*(\quad; \Lambda)$. This {\em naturality} is an important guide for developing various theories of characteristic classes for singular varieties. The known {\em functorial characteristic classes\/} for singular spaces
are {\em covariant\/} functorial maps
$$c \ell_{*}: A(X)\to H_{*}(X;\Lambda) $$
from a suiable covariant theory $A$ depending on the choice of $c \ell_{*}$. Moreover, there is always a {\em distinguished element\/}
$\jeden_{X}\in A(X)$ such that the corresponding {\em characteristic class of the singular space\/} $X$ is defined as
$c \ell_{*}(X):=c \ell_{*}(\jeden_{X})$. Finally one has the {\em normalization\/}
$$c \ell_{*}(\jeden_{M})= c \ell^{*}(TM)\cap [M] \in H_{*}(X;\Lambda) $$
for $M$ a smooth manifold, with $c \ell^{*}(TM)$ the corresponding characteristic cohomology class of $M$.
This justifies the notation $c \ell_{*}$ for this homology class transformation, which should be seen as a relative homology class
version of the following {\em characteristic number\/} of the singular space $X$:
$$\sharp (X):=c \ell_{*}(const_{*}\jeden_{X})= const_{*}(c \ell_{*}(\jeden_{X})) \in H_{*}(\{pt\};\Lambda)=\Lambda \:,$$
with $const: X\to \{pt\}$ a constant map. Note that the {\em normalization\/} implies for $M$ smooth:
$$\sharp (M)=deg(c \ell_{*}(M))=\int_M \;c \ell^{*}(TM) \cap [M]$$
so that this is consistent with the notion of characteristic number of the smooth manifold $M$ as used before.\\

\subsection{Stiefel-Whitney classes $w_{*}$}
The first example of functorial characteristic classes is the theory of singular Stiefel--Whitney homology classes due to {\em Dennis Sullivan\/} \cite {Sullivan} (also see \cite{Fulton-MacPherson}). A crucial fact about the original Stiefel--Whitney class is the following fact: if $T$ is any triangulation of a manifold $X$, then the sum of all the simplices of the first barycentric subdivision is a {\em mod 2 cycle} and its homology class is equal to the Poincar\'e dual of the Stiefel--Whitney class.  In \cite{Sullivan} D. Sullivan observed that also a {\em singular real algebraic variety\/} $X$ is a {\em mod 2 Euler space\/}, i.e. the link of any point
of $X$ has even Euler characteristic. And this condition implies that
the sum of all the simplices of the first barycentric subdivision of any triangulation of $X$ is always a {\em mod 2 cycle\/} and he defined its homology class to be the singular Stiefel--Whitney class of the variety $X$. Then, with an insight of {\em Deligne\/}, Sullivan's Stiefel--Whitney homology classes where enhanced as a natural tansformation from a certian covariant functor to the mod 2 homology theory.\\

Let $X$ be a complex (or real) algebraic set  and let $F(X)$ (or $F^{mod 2}(X)$) be the abelian group of $\bZ$- (or $\bZ_{2}$-)valued complex (or real) algebraically constructible functions on a variety $X$. Then the assignment $F$ (or $F^{mod 2}$)
$: \Cal V \to \Cal A$ is a {\em contravariant\/} functor (from the category of algebraic  varieties to the category of abelian groups) by the usual functional pullback for a morphism $f: X \to Y$:
$ f^*(\alp) := \alp \circ f $.
For a constructible set $Z \subset X$, we define
$$\chi(Z; \alp) := \sum_{n \in \bZ} n\cdot  \chi_{c}(Z \cap \alp^{-1}(n)) \quad (mod\:2).$$
Then it turns out that the assignment $F$ (or $F^{mod 2}$)$: \Cal V \to \Cal A$ also becomes a {\em covariant\/} functor by the following pushforward defined by
$$ f_*(\alp)(y) := \chi (f^{-1}(y); \alp) \quad \text{for $y\in Y$.}$$
To show that this is welldefined (i.e. $f_*(\alp)$ is again constructible) and functorial requires, for example,
stratification theory (see \cite {MacPherson1}) or a suitable theory of constructible sheaves (see \cite{Schuermann-book}).
For later use we also point out, that here in the (semi-)algebraic context we do {\em not\/} need the assumption that
our spaces are compact or the morphism $f$ is proper for the defintion of $f_*$.
This properness of $f$ for the definition of $f_*$ is only needed in the corresponding (sub-)analytic context.\\

The above Sullivan's Stiefel--Whitney class is now the special case of the following {\em Stiefel--Whitney class transformation\/}: \begin{thm} \label{w*}
On the category of compact real algebraic varieties 
there exists a unique natural transformation
$$w_*: F^{mod 2}(\quad) \to H_*(\quad; \bZ_2)$$
satisfying the normalization condition that for a nonsingular variety $X$
$$w_*(\jeden_X) = w^{*}(TX) \cap [X]\:.$$
Here $\jeden_X:=1_{X}$ is the characteristic function of $X$.
\end{thm}
Note that $\sharp(X)=deg(w_{*}(\jeden_X))= \chi(X) \:mod\;2$ is just the Euler characteristic $mod\;2$ of
the singular space $X$.

\subsection{Chern classes $c_{*}$}
Based on {\em Grothendieck's\/} ideas or modifying Grothendieck's conjecture on a {\em Riemann--Roch type formula\/} concerning the constructible \'etale sheaves and Chow rings (see \cite[Part II, note ($87_{1}$), p.361 ff.]{Grothendieck2}), {\em Deligne\/} made the following conjecture --- this is usually simply phrased  ``Deligne and Grothendieck made the following conjecture" --- and 
{\em R. MacPherson\/} \cite{MacPherson1} proved it affirmatively: 
\begin{thm} There exists a unique natural transformation 
$$c_*: F(\quad) \to H_{2*}(\quad; \bZ) $$ 
from the constructible function covariant functor $F$ to the integral homology covariant functor (in even degrees)
$H_{2*}$, satisfying the``normalization" that the value of the characteristic function $\jeden_X:=1_{X}$ of a smooth complex algebraic
variety $X$ is the Poincar\'e dual of the total Chern cohomology class:
$$c_*(\jeden_X)= c^{*}(TX)\cap[X] \:.$$
\end{thm}
The main ingredients are {\em Chern--Mather classes, local Euler obstruction and ``graph construction"\/}. The uniqueness follows from the above normalization condition and resolution of singularities. For an algebraic version of the Chern class transformation
$c_*$ over a base field of characteristic zero (taking values in Chow groups), compare with \cite{Ken}.
MacPherson's approach \cite{MacPherson1} also works in the complex analytic context, since
the analyticity of the ``graph construction" was solved by Kwieci\'{n}ski in his
thesis \cite{Kwiecinski2}.

\begin{rem}(see \cite{KMY}) 
The individual component $c_i: F(\quad) \to H_{2i}(\quad)$ of the transformation $c_*: F(\quad) \to H_{2*}(\quad)$ is also a natural transformation and also any linear combination of these components is a natrual transfomation. Let us consider {\em projective\/} varieties. Then, {\em modulo torsion\/}, these linear combinations are the {\em only\/} natural tansformations from the covariant functor $F$ to the homology functor. In particular, the {\em rationalized\/} Chern--Schwartz--MacPherson class transformation 
$c_*\otimes \bQ$ is the only such natural tansformation satisfying the {\em weaker normalization consition\/} that for each complex projective space $\bold P$ the top dimensional component of $c_*(\bold P)$ is the fundamental class $[\bold P]$. A noteworthy feature of the proof of these statements is that one does not need to appeal to the resolution of singularities.\end{rem}  

{\em J.-P. Brasselet and M.-H. Schwartz\/} [BrSc] showed that the distinguished value $c_*(\jeden_X)$ of the
characteristic function of a complex  variety embedded into a complex manifold is isomorphic to the {\em Schwartz class\/} 
\cite{Schwartz1, Schwartz2} via the Alexander duality. Thus the above transformation $c_*$ is usually called the 
{\em Chern--Schwartz--MacPherson class transformation\/}. For a complex algebraic variety $X$, singular or nonsingular, 
we have the distinguished element $\jeden_X:=1_X$ and $c_*(X):=c_*(\jeden_X)$ is called the total Chern--Schwartz--MacPherson class of 
$X$. By considering mapping $X$ to a point, one gets 
$$ \chi(X) = deg(c_{*}(\jeden_X))=  \sharp(X) \:,$$
which is a singular version of the Gauss--Bonnet--Chern theorem.

\begin{rem} For a singular version of the Poincar\'e--Hopf theorem in terms of stratified vector fields see \cite{BLSS}, and for a version in terms of $1$-forms and characteristic cycles of constructible functions, compare for example with 
\cite [\S 5.0.3]{Schuermann-book} and \cite{Schuermann-index}. There are also other notions of Chern classes of a singular complex algebraic variety $X$: Chern-Mather classes $c_*^{Ma}(X)$ (\cite{MacPherson1}), Fulton- and Fulton-Johnson Chern classes
$c_*^F(X), c_*^{FJ}(X)$ (\cite{Fulton-Johnson} and \cite[Ex. 4.2.6]{Fulton-intersection}), and for ``stringy and arc Chern classes"
$c_*^{str}(X), c_*^{arc}(X)$ see subsection \ref{stringy-arc classes}. In many interesting cases these can be described as 
$c_*(\alpha_X)$ for a suitable constructible function $\alpha_X$ related to some geometric properties of the singular space
$X$ (compare \cite{Aluffi-lecture note, Brasselet-lecture note, Parusinski-lecture note, Schuermann-VRR, Schuermann-lecture note,
Schuermann-index, Suwa-lecture note}).  Of course $\alpha_X=1_X$ for $X$ smooth, but in general $\alpha_X \neq 1_X$ so that
the MacPherson Chern class transformation $c_*$ is the basic one, but in general $\jeden_X=1_X$
is not the only possible choice of a distinguished element $\jeden_X$!
\end{rem}

\subsection{Todd classes $td_{*}$}
Motivated by the formulation of the Chern--Schwartz--MacPherson class transformation, {\em P. Baum, W. Fulton and R. MacPherson\/} 
\cite{Baum-Fulton-MacPherson} have extended {\bf GRR} to singular varieties, by introducing the so-called 
{\em localized Chern character}
$ch_X^M(\Cal F)$ of a coherent sheaf $\Cal F$ with $X$ embedded into a non-singular quasi-projective variety $M$, as a substitute of 
$ch^{*}(F)\cap[X]$  in the above {\bf GRR}.  Note that if $X$ is smooth  $ch_X^X(\Cal F) = ch^{*}(F)\cap [X]$. In [BFM] they showed the following theorem:
\begin{thm} (Baum--Fulton--MacPherson's Riemann--Roch)\\ 
\noindent(i) $td_* (\Cal F) := td^{*}(i_M^*T_M) \cap ch_X^M(\Cal F)$  is independent of the embedding $i_M : X \to M$.\\
\noindent(ii) Let the transformation $td_*: \bold G_0(\quad) \to H_{2*}(\quad;\bQ)$ be defined by  
$$td_*(\Cal F) = td^{*}(i_M^*T_M)\cap ch_X^M(\Cal F)$$ 
for any variety $X$. Then $td_* $ is actually natural, i.e., for any morphism $f:X \to Y$  the following diagram commutes:
$$\CD\bold G_0(X) @>td_*>> H_{2*}(X;\bQ)\\@V f_!VV  @VV f_* V\\\bold G_0(Y) @>>td_* > H_{2*}(Y;\bQ)\endCD$$       
i.e., for any embeddings $i_M :X \to M$ and $i_N : Y \to N$ 
\begin{equation}
td^{*}(i_N^*T_N)\cap ch_Y^N(f_! \Cal F)= f_*(td^{*}(i_M^*T_M)\cap ch_X^M(\Cal F))\:.
\tag{BFM-RR}
\end{equation}
\end{thm}
For a complex algebraic variety $X$, singular or nonsingular, $td_*(X) := td_*(\Cal O_X)$ is called the Baum--Fulton--MacPherson's Todd homology class of $X$, i.e. the class of the structure sheaf is the distingiuished element $\jeden_X:=[\Cal O_X]$. And we get 
$$ \chi_a(X) = \int_X td_*(X)= \sharp(X) \:,$$
which is a singular version of the Riemann--Roch theorem. And in \cite{Baum-Fulton-MacPherson2} this Todd class transformation
is moreover factorized through complex K-homology, which maybe is the most natural formulation of this transformation.
For the algebraic version of the Todd class transformation $td_*$ over any base field compare with 
\cite[chapter 18]{Fulton-intersection}.

\begin{rem}[Euler homology class $e_{0}$]
Even though the formulation of the BFM--RR was motivated by that of the Chern--Schwartz--MacPherson class, it was proved in a completely different way. And now there is available a similar proof of MacPherson's theorem for the embedded context based on the theory of characteristic cycles $CC$ of constructible functions, with the Segre class $s_*CC$ of these conic characteristic cycles playing the role of the localized Chern character in the proof of Baum--Fulton--MacPherson. 
Here these characteristic cycles are conic Lagrangian cycles in $T^{*}M|X$, and the pullback 
$$e_{0}:=k^{*}CC: F(X)\to  H_{0}(X;\bZ)$$
by the zero section $k: X\to T^{*}M|X$ can be seen as a functorial {\em Euler homology class transformation\/}
even in the context of real geometry. In particular 
$$\chi(X)= deg(e(\jeden_X))= \sharp(X)$$
also in this context. For more details of this, see \cite {Schuermann-lecture note, Schuermann-index}.
Finally, this approach by characteristic cycles also gives a new approach to the
Stiefel-Whitney class transformation $w_{*}$ of Sullivan as observed and explained in \cite{FuMC}.
\end{rem}

\subsection{L-classes $L_{*}$}
Using the notion of ``perversity", {\em M. Goresky and R. MacPherson\/} (\cite{Goresky-MacPherson1}, \cite{Goresky-MacPherson2}) 
have introduced {\em Intersection Homology Theory\/}. In \cite{Goresky-MacPherson1} they 
introduced a homology $L$-class $L_*^{\op {GM}}(X)$ 
for stratified spaces $X$ with even (co)dimensional strata
such that if $X$ is nonsingular 
it becomes the Poincar\'e dual of the original Thom--Hirzebruch $L$-class:
$L_*^{\op {GM}}(X) = L^{*}(TX) \cap [X]$. And for {\em rational PL-homology manifolds\/},
their $L$-classes agree with the classes introduced by Thom long ago in \cite{Thom} as one of the first characteristic classes
of suitable singular spaces.\\

Later, {\em S. Cappell and J. Shaneson\/} \cite{Cappell-Shaneson1} (see also \cite{Cappell-Shaneson2} and \cite{Shaneson}) introduced 
a homology $L$-class transformation $L_*$, which turns out to be a natural 
transformation from the abelian group $\Omega(X)$ of cobordism classes of 
selfdual constructible complexes, whose definition we now explain, to the 
rational homology group \cite{BSY2} (cf. \cite{Yokura-TAMS}).\\

Let $X$ be a compact complex analytic (algebraic) space with $D^{b}_{c}(X)$ 
the bounded derived
category of complex analytically (algebraically) constructible complexes of 
sheaves of
$\bQ$-vector spaces (compare \cite{Kashiwara-Schapira} and \cite{Schuermann-book}).
So we consider bounded sheaf complexes $\Cal{F}$, which have locally 
constant cohomology
sheaves with finite dimensional stalks along the strata of a complex 
analytic (algebraic)
Whitney stratification of $X$.
This is a triangulated category with translation
functor $T=[1]$ given by shifting a complex one step to the left.
It also has a duality in the sense of Youssin \cite{You} induced by
the {\em Verdier duality functor\/} (compare \cite[Chap.4]{Schuermann-book} 
and \cite[Chap.VIII]{Kashiwara-Schapira}):
$$ D_{X}:= Rhom(\cdot, k^{!}\bQ_{pt}): \: D^{b}_{c}(X)\to D^{b}_{c}(X)\:,$$
with $k:X\to \{pt\}$ a constant map, together with its
{\em biduality isomorphism\/} $can: id\stackrel{\sim}{\to} D_{X}\circ 
D_{X}$.
A constructible complex $\Cal{F}\in ob(D^{b}_{c}(X))$ is called {\em 
selfdual\/},
if there is an isomorphism
$$ d: \Cal{F} \stackrel{\sim}{\to} D_{X}(\Cal{F}) \:.$$
The pair $(\Cal{F},d)$ is called {\em symmetric\/} or {\em 
skew-symmetric\/}, if
$$ D_{X}(d)\circ can = d \quad \text{or} \quad D_{X}(d)\circ can = -d \:.$$
Finally an isomorphism or {\em isometry\/} of selfdual objects 
$(\Cal{F},d)$ and
$(\Cal{F}',d')$
is an isomorphism $u$ such that the following diagram commutes:
$$ \begin{CD}
\Cal{F} @> u > \sim > \Cal{F}' \\
@V d VV @VV d' V \\
D_{X}(\Cal{F}) @< \sim < D_{X}(u) < D_{X}(\Cal{F}') \:.
\end{CD} $$

The isomorphism classes of such (skew-)symmetric selfdual complexes form a 
set,
which becomes a {\em monoid\/} with addition induced by the direct sum.
Using a definition of Youssin \cite{You}, the {\em cobordism groups\/} 
$\Omega_{\pm}(X)$
of (skew-)symmetric selfdual constructible complexes on $X$
are defined by introducing a suitable {\em cobordism relation\/} in terms 
of an {\em octahedral diagram\/}, i.e.  a diagram $(Oct)$ of the following 
form:

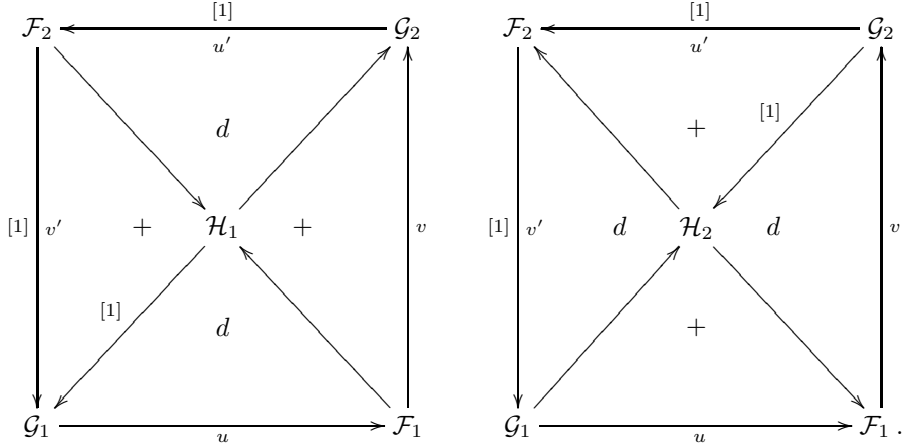
\xymatrix{\Cal{F}_{2} \ar[dddd]_{[1]}^{v'} \ar[ddr] & & 
\Cal{G}_{2} \ar[ll]^{u'}_{[1]} &  \Cal{F}_{2} \ar[dddd]_{[1]}^{v'}
 & &  \Cal{G}_{2} \ar[ll]^{u'}_{[1]} \ar[ddl]_{[1]} \\
& d & & & + \\
& + \quad\quad \Cal{H}_{1} \quad\quad + \ar[ddl]_{[1]} \ar[uur] & & &
 d\quad\quad \Cal{H}_{2}\quad\quad d \ar[uul] \ar[ddr] \\
& d & & & + \\
\Cal{G}_{1} \ar[rr]_{u} & & \Cal{F}_{1} \ar[luu] \ar[uuuu]_{v} 
& \Cal{G}_{1} \ar[rr]_{u} \ar[ruu]
 & & \Cal{F}_{1}  \ar[uuuu]_{v} \:.}

Here the morphism marked by $[1]$ are of degree one, the triangles marked 
$+$ are commutative,
and the ones marked $d$ are distinguished. Finally the two composite 
morphisms from $\Cal{H}_{1}$
to $\Cal{H}_{2}$ (via $\Cal{G}_{1}$ and $\Cal{G}_{2}$) have to be the same, 
and similarly for
the two composite morphisms from $\Cal{H}_{2}$
to $\Cal{H}_{1}$ (via $\Cal{F}_{1}$ and $\Cal{F}_{2}$).\\

Application of the duality functor $D:=D_{X}$ and a rotation by $180^{o}$ 
about the axis connecting upper-left and lower-right corner induces
another octahedral diagram $(RD\cdot Oct)$ such that $RD$ applied to 
$(RD\cdot Oct)$
gives the octahedral diagram $(D^{2}\cdot Oct)$ which one gets from $(Oct)$ 
by application of
$D^{2}$ (compare with \cite[p.387/388]{You}).
Then the octahedral diagram $(Oct)$ is called {\em symmetric\/}
or {\em skew-symmetric\/}, if there is an isomorphism $d: (Oct)\to (RD\cdot 
Oct)$ of
octahedral diagrams such that
$$ RD(d)\circ can = d \quad \text{or} \quad RD(d)\circ can = -d $$
as maps of octahedral diagrams $(Oct)\to (RD\cdot Oct)$.
Note that this induces in particular 
(skew-)symmetric dualities $d_{1}$ and $d_{2}$ of the corners $\Cal{F}_{1}$ 
and $\Cal{F}_{2}$,
and $(Oct,d)$ is called an {\em elementary cobordism\/} between 
$(\Cal{F}_{1},d_{1})$ and
$(\Cal{F}_{2},d_{2})$. This notion is a symmetric and reflexive relation.
$(\Cal{F},d)$ and $(\Cal{F}',d')$ are called {\em cobordant\/},
if there is a sequence 
$$ (\Cal{F},d)=(\Cal{F}_{0},d_{0}),\: (\Cal{F}_{1},d_{1}),\:\dots \:, 
(\Cal{F}_{m},d_{m})=(\Cal{F}',d') $$
with $(\Cal{F}_{i},d_{i})$ elementary cobordant to 
$(\Cal{F}_{i+1},d_{i+1})$ for
$i=0,\dots,m-1$. This {\em cobordism relation\/} is then an equivalence 
relation.\\

The {\em cobordism group\/} $\Omega_{\pm}(X)$ 
of selfdual constructible complexes on $X$ is the quotient of the monoid of 
isomorphism classes
of (skew-)symmetric selfdual complexes by this
cobordism relation. These are indeed abelian groups and not just monoids.\\

Consider now an algebraic (or holomorphic) map $f: X\to Y$,
with $X,Y$ compact so that $f$ is proper.
Then $Rf_{*}\simeq Rf_{!}$ maps $D^{b}_{c}(X)$ to $D^{b}_{c}(X)$.
Moreover, the {\em adjunction isomorphism\/}
$$ Rf_{*}Rhom(\Cal{F},f^{!}k^{!}\bQ_{pt}) \simeq
Rhom(Rf_{!}\Cal{F},k^{!}\bQ_{pt}) $$
induces the isomorphism
\begin{equation} \label{eq:dual1}
Rf_{*}D_{X} \stackrel{\sim}{\to} D_{Y}Rf_{!} \simeq  D_{Y}Rf_{*} 
\end{equation}
so that $Rf_{*}$ {\em commutes with Verdier-duality\/}.
In particular $Rf_{*}$ maps 
selfdual constructible complexes on $X$ to selfdual constructible complexes
on $Y$ inducing group homomorphisms
$$ f_{*}: \Omega_{\pm}(X)\to \Omega_{\pm}(Y);\: [(\Cal{F},d)]\mapsto 
[(Rf_{*}\Cal{F},Rf_{*}(d))] \:.$$

A simple example of a  selfdual constructible complex is the shifted 
constant sheaf
$\bQ_{Z}[n]$ on a complex manifold $Z$ of pure dimension $n$, with the 
duality
isomorphism induced from the {\em complex orientation\/} of $Z$
by Poincar\'e-Verdier duality:
$$ k^{!}\bQ_{pt} \simeq \bQ_{Z}[2n] \quad,\text{with $k:X\to \{pt\}$ a 
constant map.}$$
This is (skew-)symmetric for $n$ even (or odd).\\

Then the results of Cappell-Shaneson [CS1, \S 5]
can be reformulated as in \cite{BSY2}(cf. \cite[Corollary 2.3]{Yokura-TAMS}):
\begin{thm}[Cappell-Shaneson] \label{thm:Lclass}
For a compact complex analytic (or algebraic) space $X$ there is a {\em 
homology $L$-class transformation\/}
$$ L_{*}: \Omega(X) :=  \Omega_{+}(X)\oplus \Omega_{-}(X) \to H_{*}(X,\bQ) 
\:,$$
which is a group homomorphism functorial for the pushdown $f_{*}$
induced by a holomorphic (or algebraic) map.
The degree of $L_{0}((\Cal{F},d))$ is the {\em signature\/} of the induced 
pairing 
$$ H^{0}(X,\Cal{F})\otimes_{\bQ}\bR \times H^{0}(X,\Cal{F})\otimes_{\bQ}\bR 
\to \bR $$
(by definition this is $0$ for a skew-symmetric pairing).
Moreover, for $X$ smooth of pure dimension $n$ one has the normalization 
$$ L_{*}((\bQ_{X}[n],d)) = L^{*}(TX)\cap [X] \:. $$
\end{thm}

There is also a {\em uniqueness} statement in \cite[\S 5]{Cappell-Shaneson1}
for such an $L$-class transformation, but for this one has to go outside 
the
complex algebraic or analytic context.\\

For $X$ pure dimensional (otherwise one should only look at the top 
dimensional
irreducible components of $X$) one has the distinguished self-dual 
constructible
intersection cohomology complex $\jeden_{X}:=\Cal{IC}_{X}$, whose global cohomology
calculates the intersection (co)homology of Goresky-MacPherson.
By definition one gets $L_{*}(X):= L_{*}(\Cal{IC}_{X})=L_*^{\op {GM}}(X)$
so that 
$$ \int_X \: L_{*}(X)=\sharp (X)$$
is the signature of the global intersection (co)homology.

\begin{rem} Thom used in \cite{Thom} his combinatorial $L$-classes for the definition of {\em combinatorial 
Pontrjagin classes\/} of rational PL-homology manifolds. Note that in the context of rational homology manifolds,
{\em rational\/} $L$- and Pontrjagin classes carry the same information (i.e. can be deduced from each other).
But this is not the case for more singular spaces, and only a corresponding $L$-class transformation exists
for suitable singular spaces, but not a Pontrjagin class transformation.
\end{rem}

So all these theories of characteristic homology class transformations for singular spaces have the same formalism,
but their existence and construction is due to completely different underlying ideas: {\em mod 2 Euler spaces\/}
for $w_{*}$, {\em local Euler obstruction\/} for  $c_{*}$, {\em localized Chern character\/} for $td_{*}$ and {\em duality\/} for 
$L_{*}$.
Nevertheless it is natural to ask for another theory of characteristic homology classes of singular spaces,
which unifies these theories for complex algebraic varieties:
 
\begin{prob}\label{unification problem}(cf. \cite{MacPherson2} and \cite{Yokura-Banach Center}) 
Is there a ``unifying and singular version" $\displaystyle \boxed {\bold ?}_y$ of the generalized Hirzebruch--Riemann--Roch {\bf g-HRR} such that\\
\noindent(y = -1) $\boxed {\bold  ?}_{-1}$ gives rise to the rationalized Chern--Schwartz--MacPherson's class $c_* \otimes \bQ$,\\ \noindent(y = 0)  $\boxed {\bold ?}_{0}$ gives rise to the Baum--Fulton--MacPherson's Todd class $td_*$, and\\
\noindent(y = 1) $\boxed {\bold ?}_{1}$ gives rise to the Cappell--Shaneson's homology L-class $L_*$.
\end{prob}
An obvious obstacle for this problem is that the source covariant functors of these three natural transformations are all different.
And even if such a theory is not known, its {\em normalization condition\/} for a smooth complex algebraic manifold $M$ has to be
$$c \ell_{*}(\jeden_{M})= \widetilde td_{(y)}(TM)\cap [M]$$ 
by {\bf g-HRR} so that this transformation has to be called a
{\em Hirzebruch $\widetilde td_{(y*)}$- or $T_{y*}$-class transformation\/}.

\section{Relative Grothendieck rings of varieties and motivic characteristic classes}\label{relative GR}
A ``reasonable" answer for the above Problem \ref{unification problem} has been obtained in \cite{BSY2} via the so-called {\it relative Grothendieck ring of complex algebraic varieties over $X$}, denoted by $K_0(\Cal V / X)$. This ring was introduced by E. Looijenga in \cite{Looijenga} and further studied by F. Bittner in \cite{Bittner}.
The relative Grothendieck group $K_0(\Cal V / X)$ ( of morphisms over a variety $X$) is the quotient of the free abelian group of isomorphism classes of morphisms to $X$ (denoted by $[Y \to X]$ or $[Y  \xrightarrow {h} X]$), modulo the following {\it additivity} relation:
$$[Y   \xrightarrow {h}  X] = [Z \hookrightarrow Y  \xrightarrow {h}  X] +  [Y \setminus Z \hookrightarrow Y   \xrightarrow {h} X]$$
for $Z \subset Y$ a closed subvariety of $Y$. The ring structure is given by the fiber square: for 
$[Y   \xrightarrow {f} X], [W   \xrightarrow {g}  X] \in K_0(\Cal V /X)$
$$[Y  \xrightarrow {f} X] \cdot [W  \xrightarrow {g}  X]:= [Y \times_X W    \xrightarrow {f \times_X  g} X].$$ 
Here $Y \times_X W     \xrightarrow {f \times_X  g} X$ is $g \circ f' = f \circ g'$ where $f'$ and $g'$ are as in the following diagram
$$\CD Y \times_X W  @> {f'} >>  W'\\@V {g' }VV @VV g V\\ Y @> f >>X.  \endCD$$
The relative Grothendieck ring $K_0(\Cal V /X)$ has the unit $1_X := [X \xrightarrow {id_X} X]$,
which later becomes the distinguished element $\jeden_{X}:=[id_{X}]$. 
Similarly one gets an exterior product 
$$\times : K_0(\Cal V / X) \times K_0(\Cal V / Y) \to K_0(\Cal V / X\times Y) \:.$$
Note that when $X = \{pt\}$ is a point, then the relative Grothendieck ring $K_0(\Cal V/ \{pt\})$ is nothing but the usual Grothendieck ring 
$K_0(\Cal V)$ of $\Cal V$, which is the free abelian group generated by the isomorphism classes of varietiesmodulo the subgroup generated by elements of the form 
$[V] - [V'] - [V \setminus V']$ for a subvariety $V' \subset V$, and the ring structure is given by the Cartesian product of varieties.

\begin{rem} In some sense the Grothendieck ring $K_0(\Cal V)$ can be seen as an algebraic substitute for cobordism rings
$\Omega_{*}$ of smooth manifolds, based on the {\em additivity\/} instead of a cobordism relation.
\end{rem} 

For a morphism $f: X' \to X$, the pushforward
$$f_*: K_0(\Cal V /X') \to K_0(\Cal V /X)$$
is defined by
$$f_*[Y \xrightarrow {h} X'] := [Y \xrightarrow {f \circ h} X].$$
With this pushforward, the assignment $X \longmapsto K_0(\Cal V /X)$ is a covariant functor. The pullback 
$$f^*: K_0( \Cal V /X) \to K_0(\Cal V /X')$$
is defined as follows: for a fiber square
$$\CD Y' @> {g'} >>  X'\\@V {f' }VV @VV f V\\ Y @> g >>X  \endCD$$
the pullback $f^*[Y \xrightarrow {g} X]:= [Y' \xrightarrow {g'} X'].$
With this pullback, the assignment $X \longmapsto K_0(\Cal V /X)$ is a contravariant functor. Let $\op {Iso}^{\op {pr}}(\Cal {SV}/X)$ be the free abelain groups on isomorphism classes of proper morphisms from smooth varieties to a given variety $X$. Then we get the canonical quotient homomorphism
$$quo: \op {Iso}^{\op {pr}}(\Cal {SV}/X) \to  K_0( \Cal V /X)$$
which is surjective by the above additivity relation and Hironaka's resolution of singularities \cite{Hironaka}. And it turns out that the kernel of this surjective map is generated by the ``blow-up relation", more precisely we have the following theorem, which is due to F. Bittner [Bi, Theorem 5.1], based on the very deep ``weak factorization theorem" (\cite{AKMW} and \cite{Wlodarczyk}):
\begin{thm}\label {Bittner}
The relative Grothendieck group $K_0( \Cal V /X)$ is isomorphic to the quotient of the free abelian group 
$\op {Iso}^{\op {pr}}(\Cal {SV}/X)$ modulo the following ``blow-up relation"
$$[\emptyset \to X] := 0 \quad \text {and} \quad [Bl_YX' \to X] - [E \to X] = [X' \to X] - [Y \to X]$$
for any Cartesian ``blow-up" diagram
$$\CD E @> {i'} >>  Bl_YX'\\@V {\pi' }VV @VV {\pi} V\\ Y @> i >>X' @> f >>X, \endCD$$
with $i$ a closed embedding of smooth (pure dimensional) varieties and $f: X' \to X$ proper. Here $\pi: Bl_YX' \to X'$ is the blow-up of $X'$ along $Y$ with $E$ denoting the exceptional divisor.
\end{thm} 
From this theorem we can get the following corollary:
\begin{thm}\label{blowup} 
Let $B_*: \Cal V/k \to \Cal A$ be a functor from the category of reduced separated schemes of finite type over $\bC$ to the category of abelian groups such that\\ 
\noindent (i) $B_*(\emptyset):= {0}$,\\
\noindent (ii) it is covariantly functorial for proper morphisms, and\\ 
\noindent (iii) for any smooth variety $X$  there exists a distinguished element $d_X \in B_*(X)$ such that\\
(iii-1) for any isomorphism $h: X' \to X$, $h_*(d_{X'}) = d_X$ and \\
(iii-2) for any Cartesian ``blow-up" diagram as in the above Theorem \ref{Bittner} with $f = \op {id}_X$,
$$\pi_*(d_{Bl_YX}) -i_*\pi'_*(d_E) = d_X - i_*(d_Y) \in B_*(X).$$
Then we have by (iii-1) that
 there exists a unique natural transformation of covariant functors
$$\Phi: \op {Iso}^{\op {pr}}(\Cal {SV}/ \quad)  \to B_*(\quad)$$
satisfying the normalization condition that for smooth $X$
$$\Phi([X  \xrightarrow {\op {id}}  X]) = d_X,$$
and furthermore
by (iii-2) there exists a unique natural transformation of covariant functors
$$\widetilde \Phi: K_0( \Cal V / \quad) \to B_*(\quad)$$
satisfying the normalization condition that for smooth $X$
$$\widetilde \Phi ([X  \xrightarrow {\op {id}}  X]) = d_X.$$
\end{thm}
Then, using results of \cite[IV.1.2.1]{Gros} or \cite[Proposition 3.3]{GNA}, we can get the following corollary about a
{\em motivic Chern class transformation $mC_{*}$\/}.
\begin{cor} \label{mC*}
There exisits a unique natural transformation (with respect to proper maps)
$$mC_*:  K_0( \Cal V / \quad) \to \bG_0(\quad)\otimes \bZ[y]$$
 satisfying the normalization condition that for $X$ smooth
$$mC_*([X  \xrightarrow {\op {id}}  X]) = \sum_{i=0}^{\op {dim}X} [\Lambda^i T^*X]y^i =:
\Lambda_{y}([T^*X]) \cap [\Cal{O}_{X}].$$
Here $\Lambda_{y}(\quad)$ is the so-called {\em total $\Lambda$-class\/}.
\end{cor}

If we compose $mC_*|_{y=-1,0,1}$ with the natural transformation $\bG_0(\quad)\to \bold K_0^{top}(\quad)$
to topological K-homology constructed in \cite{Baum-Fulton-MacPherson2}, then $mC_*(X)$ unifies for $X$ smooth the following
K-theoretical homology classes:\\
(y=-1) \quad the top-dimensional Chern class $c^{top}_K(TX)\cap [X]_K$ in K-theory
$$mC_*|_{y=-1}([id_X])= \Lambda_{-1}([T^*X]) \cap [X]_K \:,$$
(y=0) \quad the fundamental class in K-homology of the complex manifold $X$ 
$$mC_*|_{y=0}([id_X])= [X]_K \:,$$
(y=1) \quad the class of the signature operator of the underlying spin$^{c}$ manifold of $X$
$$mC_*|_{y=1}([id_X])= \Lambda_{1}([T^*X]) \cap [X]_K \:.$$

Consider the twisted BFM--RR transformation 
$$td_{(1+y)}: \bG_0(X)\otimes \bZ[y] \to H_{2*}(X)\otimes Q[y, (1+y)^{-1}]$$ 
defined by 
$$td_{(1+y)}([\Cal F]) := \sum_{i \geq 0}td_i([\Cal F])(1+y)^{-i}$$ 
and extending it linearly with respect to $\bZ[y]$ 
(\cite {Yokura-Banach Center}). Using this twisted BFM--RR transformation $td_{(1+y)}$ and the above transformation $mC_*$, we define the {\em Hirzebruch class transformation $T_{y*}$\/} as the composite $T_{y*}:= td_{(1+y)} \circ mC_*$. Then we get the following theorem:
\begin{thm}\label{T_y_*} 
Let $K_0(\Cal V/X)$ be the Grothendieck group of complex algebraic varieties over $X$. Then there exists  
a unique natural transformation (with respect to proper maps)
$${T_y}_*: K_0(\Cal V/\quad) \to H^{BM}_{2*}(\quad)\otimes \bQ[y] \subset
H^{BM}_{2*}(\quad)\otimes \bQ[y, (1+y)^{-1}] $$
such that for $X$ nonsingular
$${T_y}_*([X \xrightarrow  {\op {id}} X]) = \widetilde td_{(y)}(TX) \cap [X].$$
\end{thm}

\begin{rem} The transformations $mC_*$ and $T_{y*}$ can also be defined in the same way in the {\em algebraic context\/} over a base field of characteristic zero, using the algebraic version of the Todd tranformation $td_*$ as in \cite[chapter 18]{Fulton-intersection},
and in the {\em compactifiable complex analytic context\/}, using the analytic version of the Todd tranformation $td_*$
constructed in \cite{Levy} (compare with \cite{BSY2} for more details).
\end{rem}

For a later use, we observe that ${T_y}_*$ commutes with the exterior product
(and similarly for $mC_*$), i.e., the following diagram commutes:
$$\CD K_0(\Cal V/X) \times K_0(\Cal V/Y) @>  \times >> K_0(\Cal V/X \times Y) \\
@V {T_y}_* \times {T_y}_*VV @VV {T_y}_* V\\
H_{2*}(X) \otimes \bQ[y] \times H_{2*}(Y)\otimes \bQ[y] @> \times  >> 
H_{2*}(X \times Y)\otimes \bQ[y]. \endCD$$
And we have the following theorem for a compact complex algebraic variety $X$:
\begin{thm}\label{unification theorem}
\noindent (y = -1)  There exists a unique natural transformation $\epsilon: K_0(\Cal V/\quad) \to F(\quad)$ such that for $X$ nonsingular $\epsilon([X \xrightarrow {\op {id}} X]) = 1_X.$ And the following diagram commutes
$$\xymatrix{K_0(\Cal V/X)  \ar[dr]_ {T_{-1}}\ar[rr]^ {\epsilon} && F(X) \ar[dl]^{c_*}\\& H_{2*}(X)\otimes \bQ \:.}$$
\noindent (y = 0)  There exists a unique natural transformation $\ga: K_0(\Cal V/\quad) \to \bold G_0(\quad)$ such that for $X$ nonsingular $\ga([X \xrightarrow {\op {id}} X]) = [\Cal O_X].$ And the following diagram commutes
$$\xymatrix{K_0(\Cal V/X)  \ar[dr]_ {T_{0}}\ar[rr]^ {\ga} && \bold G_0(X) \ar[dl]^{td_*}\\& H_{2*}(X)\otimes \bQ \:.}$$
\noindent  (y = 1)  There exists a unique natural transformation $\omega: K_0(\Cal V/\quad) \to \Omega (\quad)$
such that for $X$ nonsingular $\omega([X \xrightarrow {\op {id}} X]) = \left [\bQ_X[\op {dim}X] \right ].$ And the following diagram commutes
$$\xymatrix{K_0(\Cal V/X)  \ar[dr]_ {T_{1}}\ar[rr]^ {\omega} && \Omega(X) \ar[dl]^{L_*}\\& H_*(X)\otimes \bQ \:.}$$
\end{thm}
An original proof of the above Theorem \ref{T_y_*} uses Saito's theory of mixed Hodge modules \cite{Saito} instead of the above Theorem \ref{Bittner}. And an even more elementary proof can be given based on some classical results of \cite{DuBois} about the so-called DuBois complex of a singular complex algebraic variety. Only the proof of the case $(y = 1)$ of the above Theorem \ref{unification theorem} depends, up to now, on the Bittner's theorem, i.e., the above Theorem \ref{Bittner}, in other words, on the ``weak factorization theorem" (\cite{AKMW} and \cite{Wlodarczyk}). 
Also note that the transformation $\epsilon$ is defined for {\em any\/} algebraic map of not
necessarily compact algebraic varieties, and it also commutes with pullback and (exterior) products. 
For more details, see \cite{BSY2}.

\begin{rem}
The reader should be warned that the transformations $\ga$ and $\omega$ above do {\em not\/} preserve the distinguished
elements in general. For any compact singular complex algebraic variety $X$ one has $\epsilon([id_{X}])=1_{X}$ so that the 
{\em Hirzebruch class\/}
$T_{y*}(X):=T_{y*}([id_{X}])$ specializes to $T_{-1*}(X)=c_{*}(X)\in H_{2*}(X;\bQ)$. But in general
$$ \ga([id_{X}])\neq [\Cal O_X]\in G_{0}(X)
\quad \text{and} \quad T_{0*}(X) \neq td_{*}(X) \:.$$
But $T_{0*}(X) = td_{*}(X)$ if $X$ has at most ``Du Bois singularities", e.g. ``rational singularities"like, for example,
toric varieties. Similarly 
$$ \omega([id_{X}])\neq [\Cal IC_X]\in \Omega(X)
\quad \text{and} \quad T_{1*}(X) \neq L_{*}(X) $$
in general, but we {\em conjecture\/} that $T_{1*}(X) = L_{*}(X)$ for $X$ a {\em rational homology manifold\/}.
\end{rem}

Moreover, the Hirzebruch characteristic class $\widetilde td_{(y)}=T^{*}_{y}$ is the {\em most general\/} 
normalized and multiplicative characteristic class of complex vector bundles
$$c \ell_{f}^{*}: \op {Vect} (X) \to H^{2*}(X; \Lambda) \:,$$
with $\Lambda$ a $\bQ$-algebra, which satisfies the condition of Theorem \ref{blowup}
with 
$$d_{X}:=c \ell_{f}^{*}(TX)\cap [X] \in H^{BM}_{2*}(X;\Lambda)$$ 
for $X$ smooth. In fact, the correspondig {\em genus\/} $\Phi_{f}$ factorizes as
\begin{equation} \label{eq:genus} \begin{CD} 
\op {Iso}^{\op {pr}}(\Cal {SV}/\{pt\}) @>>> \Omega^{U}_{*}\otimes\bQ  \\ 
@VVV  @VV \Phi_{f} V\\
K_{0}(\Cal V) 
@> \Phi_{f} >> \Lambda=H_{2*}(\{pt\};\Lambda)\:.
\end{CD} \end{equation}
Moreover, the characteristic class $c \ell_{f}^{*}$ or its genus $\Phi_{f}$
is uniquely determined by 
$$\Phi_{f}([P^{n}(\bC)])
= \int_{P^{n}(\bC)} \: (c \ell_{f}^{*}(TP^{n}(\bC)))\cap [P^{n}(\bC)])$$ 
for all $n$.
But if $\Phi_{f}$ also factorizes over $K_{0}(\Cal V)$ then we get from the decomposition 
$$P^{n}(\bC) = \{pt\} \cup \bC \cup \cdots \cup \bC^n$$
by ``additivity'' and ``multiplicativity'' (and compare with equation (\ref{chiy})):

\begin{equation} \label{eq:genus2}
\Phi_{f}([P^{n}(\bC)]) = 1+ (-y) + \cdots + (-y)^{n} 
\quad \text{with} \quad y:= 1- \Phi_{f}([P^{1}(\bC)]) \:.
\end{equation}
So $\Phi_{f}$ is a specialization of the {\em Hirzebruch $\chi_{y}$-genus\/}
corresponding to the {\em Hirzebruch characteristic class\/} $T^{*}_{y}$.
Of course here we use a decomposition into the {\em non-compact\/} manifolds $\bC^n$, which 
``is classically forbidden for a genus", with $y=- \Phi_{f}([\bC])$.

\begin{rem} So {\em additivity\/} is the underlying principle which ``singles out" those 
normalized and multiplicative characteristic classes $c \ell_{f}^{*}$, which have (so far)
a functorial extension to singular spaces. Also
note that the specialization $y=1$ corresponding to the {\em signature
genus\/} $sign=\chi_{1}$ and the {\em characteristic $L$-class transformation\/} $L^{*}=T^{*}_{1}$
is the only one that factorizes by the canonical map
$\Omega^{U}_{*}\otimes\bQ \to \Omega^{SO}_{*}\otimes\bQ$ over the {\em cobordism ring
$\Omega^{SO}_{*}$ of oriented manifolds\/},
since $[P^{1}(\bC)]= 0 \in \Omega^{SO}_{*}$.
In particular this ``explains" why there is no functorial {\em Pontrjagin class transformation\/}
for singular spaces.
\end{rem}

For $X$ a compact complex algebraic variety one can also deduce from Theorem \ref{blowup} the Chern class transformation
$$c_{*}: K_{0}(\Cal V /X) \to H_{2*}(X;\bZ) \:,$$
on the relative Grothen\-dieck group $K_{0}(\Cal V/X)$ without appealing to MacPherson's theorem,
since the distinguished element 
$$d_{X}:= c^{*}(TX)\cap [X] \in H_{2*}(X;\bZ)$$
of a smooth space $X$ satisfies the corresponding conditions.
Condition (iii-1) follows from the projection formula, and condition (iii-2) is an
easy application (by pushing down to $X$) of the classical ``blowing up formula
for Chern classes'' \cite[Theorem 15.4]{Fulton-intersection} . And recent work of Aluffi \cite{Aluffi-limits}
can be interpreted as showing that this
transformation $c_{*}$ factorizes over $\epsilon: K_0(\Cal V/\quad) \to F(\quad)$.

\section{Bivariant Characteristic classes}\label{bivariant characteristic classes}
In \cite{Fulton-MacPherson} (also, see \cite{Fulton-intersection}) {\em W. Fulton and R. MacPherson\/} introduced the notion of {\em Bivariant Theory\/}, which is a simultaneous generalization of a pair of covariant and contravariant functors. Most pairs of covariant and contravariant theories, e.g., such as homology theory, K-theory, etc., extend to bivariant theories. A bivariant theory $\bB$ on a suitable category $\Cal C$ (with a distinguished class of so-called ``proper" or ``confined"  maps) with values in the category of abelian groups is an assignment to each morphism $ X \xrightarrow {f}Y$ in the category $\Cal C$ an abelian group
$\bB(X \xrightarrow {f} Y)$, which is equipped with the following three basic operations:\\
\noindent (Product operations): For morphisms $f: X \to Y$ and $g: Y\to Z$, the product operation
$$\bullet: \bB( X \xrightarrow {f} Y) \otimes \bB( Y \xrightarrow {g} Z) \to\bB( X \xrightarrow {gf} Z)$$
is  defined.\\
\noindent (Pushforward operations): For morphisms $f: X \to Y$and $g: Y \to Z$ with $f$ proper, the pushforward operation
$$f_{\bigstar}: \bB( X \xrightarrow {gf}Z) \to \bB( Y \xrightarrow {g} Z) $$
is  defined.\\
\noindent (Pullback operations): For a fiber (or more generally a so-called independent) square
$$\CD X' @> g' >> X \\@V f' VV @VV f V\\Y' @>> g > Y, \endCD$$
the pullback operation
$$g^{\bigstar} : \bB( X \xrightarrow {f} Y) \to \bB( X' \xrightarrow {f'}Y') $$
is  defined. And these three operations are required to satisfy {\em seven compatibility axioms\/} (see \cite[Part I, \S 2.2]{Fulton-MacPherson} for details). 
In particular, the class of ``proper" maps has to be stable under composition and base change, and should contain all
identity maps.
Let $\bB, \bB'$ be two bivariant theories on such a category $\Cal C$. Then a {\em Grothendieck transformation\/} from $\bB$ to $\bB'$
$$\ga : \bB \to \bB'$$
is a collection of homomorphisms
$$\bB(X \to Y) \to \bB'(X \to Y) $$
for a morphism $X \to Y$ in the category $\Cal C$, which preserves the above three basic operations: \\
\hskip1cm (i) \quad $\ga (\alp \bullet_{\bB} \be) = \ga (\alp) \bullet _{\bB'} \ga (\be)$,\\ 
\hskip1cm  (ii) \quad $\ga(f_{\bigstar}\alp) = f_{\bigstar} \ga (\alp)$, and \\
\hskip1cm  (iii) \quad $\ga (g^{\bigstar} \alp) = g^{\bigstar} \ga (\alp)$.\\

$\bB_*(X):= \bB(X \to pt)$ 
and $\bB^*(X) := \bB(X \xrightarrow {\op {id}} X)$ become a covariant functor for proper maps and a contravariant functor, respectively. And a Grothendieck transformation $\ga: \bB \to \bB'$ induces natural transformations $\ga_*: \bB_* \to \bB_*'$ and $\ga^*: \bB^* \to {\bB'}^*$
such that $\ga_*$ commutes with the (bivariant) exterior product, i.e. the following diagram commutes:
$$\begin{CD} \bB_*(X)\times \bB_{*}(Y) @> \times >> \bB_*(X\times Y)\\
@V \ga_*\times \ga_* VV @VV \ga_*  V \\
\bB_*'(X)\times \bB_{*}'(Y) @> \times >> \bB_*'(X\times Y) \:.
\end{CD} $$
If we have a Grothendieck transformation $\ga: \bB \to \bB'$, then via a bivariant class $b \in \bB( X  \xrightarrow {f} Y)$ we get the commutative diagram
\begin{equation} \label{VRR1}
\CD \bB_*(Y) @> {\ga_*}>> \bB'_*(Y) \\@V {b \bullet} VV @VV {\ga(b) \bullet }V\\ \bB_*(X) @>> {\ga_*}> \bB'_*(X). \endCD
\end{equation}
This is called {\em the Verdier-type Riemann--Roch formula associated to the bivariant class $b$\/}.\\
$ $\\
{\bf Bivariant Todd class transformation $\tau$.\/}
The most important (and motivating) example of such a Grothendieck transformation of bivariant theories is the {\em bivariant
Riemann-Roch transformation\/} $\tau$ from the
{\em bivariant algebraic K-theory $\bK_{\op {alg}}$ of perfect complexes\/} to {\em rational bivariant homology\/} $\bH_{\bQ}$
$$\tau : \bK_{\op {alg}} \to \bH_{\bQ}$$
constructed in \cite[Part II]{Fulton-MacPherson} in the complex quasi-projective context. Here $\bH_{\bQ}$ is the bivariant homology theory corresponding to usual cohomology
with rational coefficients constructed in \cite[\S 3.1]{Fulton-MacPherson} for more general cohomology theories.
Then the associated contravariant theory $\bH_{\bQ}^*(X)=H^{*}(X;\bQ)$ is the usual cohomology, and the associated covariant theory
$\bH_{\bQ*}(X)=H^{BM}_{*}(X;\bQ)$ is the Borel-Moore homology. Similarly $\bold \bK^*_{\op {alg}}\simeq \bold K^{0}$ is the
Grothendieck group of algebraic vector bundles, and
$\bold \bK_{\op {alg}*}\simeq \bold G_{0}$ is the Grothendieck group of algebraic coherent sheaves.
Then the associated contravariant transformation $\tau^*$ is the {\em Chern character\/}
$$ch^*: \bold \bK_{\op {alg}}^*(\quad)\simeq \bold K^{0}(\quad)  \to H^*(\quad; \bQ)\simeq \bH_{\bQ}^*(\quad) \:,$$ 
and the associated covariant transformation 
$$\tau_*: \bold \bK_{\op {alg}*} (\quad)\simeq \bold G_{0}(\quad)  \to H^{BM}_*(\quad; \bQ)\simeq \bH_{\bQ*}(\quad)$$ 
is just Baum--Fulton--MacPherson's Todd class transformation $td_{*}$ constructed in \cite{Baum-Fulton-MacPherson}.
And the bivariant transformation $\tau$ unifies many different known Riemann-Roch type theorems.
In particular for a {\em smooth\/} morphism $f: X\to Y$ of possible singular varieties one has
$$\jeden_{f}:=[\Cal{O}_{X}] \in  \bK_{\op {alg}}(X \xrightarrow {f} Y) \:,$$
with $\tau(\jeden_{f}) =td^{*}(T_{f})\bullet [f]$.
Here $T_{f}$ is the vector bundle of tangent spaces of fibers of $f$, and $[f]\in \bH_{\bQ}(X \xrightarrow {f}  Y)$
is the {\em canonical orientation\/} of the smooth morphism $f$. Then the  Verdier-type Riemann--Roch formula 
(\ref{VRR1}) associated to 
$\jeden_{f}$ becomes the usual {\em Verdier Riemann-Roch theorem\/} for the Todd class transformation $td_*$:
\begin{equation} \label{VRRtodd}
td_*(f^* \be) = td^{*}(T_{f})\cap f^{!}td_{*}(\be) \quad \text{for $\be\in \bold G_{0}(Y)$.}
\end{equation}
Here $f^{!}=[f] \bullet : H^{BM}_{*}(Y;\bQ)\simeq \bH_{\bQ*}(Y)\to \bH_{\bQ*}(X)\simeq H^{BM}_{*}(X;\bQ)$
is the {\em smooth pullback\/} in Borel-Moore homology. And for an {\em algebraic version\/} of this bivariant 
Riemann-Roch transformation $\tau$ compare with \cite[Ex. 18.3.16]{Fulton-intersection}.\\
$ $\\
{\bf Bivariant Stiefel-Whitney class transformation $\omega$.\/}
In the context of real geometry (e.g. the piecewise linear, (semi-)algebraic or subanalytic context) one has the following
interesting example of a bivariant theory (with ``proper" the usual meaning). Here
{\em Fulton--Mac\-Pherson's bivariant group $\bF^{mod 2}(X \xrightarrow {f} Y)$ of $\bZ_2$-valued constructible functions\/} consists of all the constructible functions on $X$ which satisfy the local Euler condition with respect to $f$.  Here a $\bZ_2$-valued constructible function $\alp \in F^{mod 2}(X)$ is said to satisfy the {\em local Euler condition with respect to $f$\/}, if for any point $x \in X$ and for any local embedding $(X, x) \to (\bold R^N,0)$ the equality 
$$\alp(x) = \chi \left (B_{\epsilon} \cap f^{-1}(z);\alp \right) \, mod \, 2 $$ 
holds, where $B_{\epsilon}$ is a sufficiently small {\em open\/} ball of the origin $0$ with radius $\epsilon$ and $z$ is any point close to 
$f(x)$ (cf. \cite{Brasselet1}, \cite{Sabbah}). In particular, if $\jeden_f := 1_X$ belongs to the bivariant group $\bF^{mod 2}(X \xrightarrow {f} Y)$, then the morphism $f: X \to Y$ is called an {\em Euler morphism\/}. 
For $f: X\to \{pt\}$ a constant map this just means (by the ``local conic structure" of $X$),
that $X$ is a {\em mod 2 Euler space\/}, i.e. the link 
$\partial B_{\epsilon}\cap X$ of any point $x \in X$ has vanishing Euler characteristic modulo 2:
$$
\aligned
\chi(\partial B_{\epsilon}\cap X) & =
\chi_{c}(\partial B_{\epsilon}\cap X) \\
& = 1-\chi_{c}(B_{\epsilon}\cap X) \\
& = 1-\chi(B_{\epsilon}\cap X;1_{X}) = 0 \quad mod\;2 
\endaligned
$$
Also a {\em smooth\/} morphism, or a locally trivial fibration with fiber a mod 2 Euler space, is always an Euler morphism.\\

The three operations on $\bF^{mod 2}(X \xrightarrow {f} Y)$ are defined as follows:\\
(i) the product operation
$$\bullet: \bF^{mod 2}( X \xrightarrow {f} Y)\otimes \bF^{mod 2}( Y \xrightarrow {g} Z) \to\bF^{mod 2}( X \xrightarrow {gf} Z)$$
is  defined by $\alp \bullet \be:= \alp \cdot f^*\be.$\\
(ii) the pushforward operation $f_{\bigstar}: \bF^{mod 2}( X \xrightarrow {gf} Z) \to \bF^{mod 2}( Y \xrightarrow {g} Z) $
is the usual pushforward $f_*$, i.e., 
$$f_{\bigstar}(\alp)(y):= \chi (f^{-1}(\{y\});\alp) \, mod \, 2.$$\\
(iii) for a fiber square
$$\CD X' @> g' >> X \\@V f' VV @VV f V\\Y' @>> g > Y, \endCD$$
the pullback operation
$g^{\bigstar} : \bF^{mod 2}( X \xrightarrow {f} Y) \to \bF^{mod 2}( X' \xrightarrow {f'}Y') $
is the functional pullback ${g'}^*$, i.e., 
$$g^{\bigstar}(\alp)(x'):= \alp (g'(x')).$$
Note that for $f$ proper and any {\em bivariant\/} constructible function $\alp \in \bF^{mod 2}( X \xrightarrow {f} Y)$, the Euler--Poincar\'e characteristic 
$\displaystyle \chi \bigl (f^{-1}(y); \alp \bigr )$ of $\alp$ restricted to each fiber $f^{-1}(y)$ is {\em locally constant\/} on $Y$ mod 2 (by the local Euler condition for $f_{*}(\alp)$). \\

The correspondence $s\bF^{mod 2}(X \to Y):= F^{mod 2}(X)$ assigning to a morphism $f: X \to Y$ the abelian group $F^{mod 2}(X)$ of the source variety $X$, whatever the morphism $f$ is, becomes a bivariant theory with the same operations above. This bivariant theory is called the {\em simple\/} bivariant theory of constructible functions (see \cite{Schuermann-partial GT} and \cite{Yokura-TopAppl}). In passing, what we then need to do for showing that the Fulton--MacPherson's group of $\bZ_2$-valued constructible functions satisfying the local Euler condition with respect to a morphism is a bivariant theory, is to show that the local Euler condition with respect to a morphism is preserved by each of these three operations. \\

For later use let us point out the abstract properties needed for the definition of a {\em simple bivariant theory\/}
\cite[Definition, p.25-26]{Schuermann-partial GT}:\\
{\bf (SB1)} We have a contravariant functor $G: \Cal C \to Rings$
with values in the category of rings with unit.\\
{\bf (SB2)} $G$ is also covariantly functorial with respect to proper maps
(as a functor to the category of Abelian groups).\\
{\bf (SB3)} $G$  satisfies the {\em two-sided projection-formula\/}, i.e. for $f:X\to Y$ proper and $\alpha \in G(Y)$ and $\beta\in G(X)$, 
$$f_{*}((f^*\alpha)\cup \beta)=\alpha \cup (f_{*}\beta)\:,$$
i.e., $f_*$ is a left $G(Y)$-module and
$$f_{*}(\beta \cup (f^*\alpha))= (f_{*}\beta)\cup \alpha \:,$$
i.e., $f_*$ is a right $G(Y)$-module. (Note that we do not assume $(G,\cup)$ is (graded) commutative so that both versions of the usual projection formula are needed.)\\
{\bf (SB4)} $F$ has the {\em base-change property\/} 
$g^{*}f_{*}=f'_{*}g'^{*} : \: G(X)\to G(Y')$
for any fiber (or independent) square
\begin{displaymath} \begin{CD}
X' @>g' >> X \\
@V f' VV  @VV f V \\
Y' @> g >> Y \:,
\end{CD} \end{displaymath}
with $f,f'$ proper.\\

Then one gets a (simple) bivariant theory
$s\bG$ by $s\bG( X \xrightarrow {f} Y):=G(X)$, with the obvious push-down and pull-back transformations as above.
Finally the bivariant product
$$\bullet: s\bG( X \xrightarrow {f} Y) \times s\bG(Y \xrightarrow {g} Z) \to s\bG(X \xrightarrow {gf} Z)$$
is just given by $\alpha\bullet \beta := \alpha \cup f^{*}(\beta)$, with $\cup$ the given product
of the ring-structure. Note that this construction does not only apply to constructible functions $G(\quad)=F^{mod 2}(\quad)$,
but also to the relative Grothendieck group of complex algebraic varieties $G(\quad)=K_0(\Cal V/\quad)$, even if we allow all algebraic morphisms as ``proper" morphisms.\\

Let $\bH^{mod 2}(X \xrightarrow {f} Y)$ be Fulton--MacPherson's {\em bivariant homology theory\/} with $\bZ_2$ coefficients, constructed from the corresponding cohomology theory in \cite[\S 3.1]{Fulton-MacPherson} so that $\bH^{mod 2,*}(X)=H^{*}(X;\bZ_{2})$ and
 $\bH^{mod 2}_*(X)=H^{BM}_{*}(X;\bZ_{2})$. Then Fulton--MacPherson \cite[Theorem 6A]{Fulton-MacPherson} showed in the {\em piecewise linear context\/} the following theorem, which is a bivariant version of the singular Stiefel--Whitney class transformation $w_*: F^{mod 2}(\quad) \to H^{BM}_*(\quad: \bZ_2)$:
\begin{thm} 
There existis a unique Grothendieck transformation
$$\omega: \bF^{mod 2} \to \bH^{mod 2}$$
satisfying the normalization condition that for a morphism from a smooth variety $X$ to a point
$$\omega(1_X) = w^{*}(TX) \cap [X] \in \bH^{mod 2}_*(X)=H^{BM}_{*}(X;\bZ_{2}) \:.$$
\end{thm}

\begin{rem} As to the bivariant mod 2 constructible functions, in the context of real geometry, the definition and the theory of them 
can be given in any of the following categories: the $PL$-category, the (semi-)algebraic category and the subanalytic category. Note that the above bivariant Stiefel--Whitney class transformation is only proved and known in the $PL$-category.
\end{rem}
$ $\\
{\bf Bivariant Chern class transformation $\gamma$.\/}
Instead of mod 2 constructible functions, in the complex analytic or algebraic context we certainly have similarly the bivariant group 
$\bF(X \to Y)$ of $\bZ$-valued constructible functions satisfying the local Euler condition with
values in $\bZ$ (and not only in $\bZ_{2}$) and the bivariant homology theory $\bH(X \to Y)$ with integer coefficients, and  W. Fulton and R. MacPherson conjectured or posed as a question the existence of a so-called {\em bivariant Chern class transformation\/} and {\em J.--P. Brasselet\/} \cite{Brasselet1} solved it:
\begin{thm} 
(J.-P. Brasselet) For the category of embeddable complex analytic varieties with {\em celluar morphisms\/}, there exists a Grothendieck transformation 
$$\ga : \bF \to \bH$$
such that for a morphism $f:X \to \{pt\}$ from a nonsingular variety $X$ to a point $\{pt\}$ and the bivariant constructible function 
$\jeden_f := 1_X$ the following normalization condition holds:
$$\ga(\jeden_f) = c^{*}(TX) \cap [X] \in  \bH_{*}(X)=H^{BM}_{*}(X;\bZ) \:.$$
\end{thm}
Since then, the {\em uniqueness\/} of the Brasselet bivariant Chern class and the problem of whether `` cellularness"  of morphisms
(which is not so easy to check) can be dropped or not have been unresolved. In \cite{Sabbah} {\em C. Sabbah\/} constructed a bivariant Chern class transformation
``micro-local analytically" in some cases. In \cite {Zhou1}, \cite{Zhou2} {\em J. Zhou\/} showed that the bivariant Chern classes constructed by J.-P. Brasselet \cite{Brasselet1} and the ones constructed by C. Sabbah \cite {Sabbah} in some cases are identical in the case when the target variety is a {\em nonsingular curve\/}. And in \cite [Theorem (3.7)]{Yokura-Documenta} we showed the following more general {\em uniqueness theorem\/} of bivariant Chern classes for morphisms whose target varieties are {\em nonsingular of any dimension\/}:
\begin{thm}
If there exists a bivariant Chern class transformation $\ga : \bF \to \bH$, then it is unique when restricted to morphisms whose target varieties are nonsingular; explicitly, for a morphism $f : X \to Y$ with $Y$ nonsingular and for any bivariant constructible function 
$\alpha \in \bF(X \xrightarrow {f } Y)$ the bivariant Chern class $\ga (\alpha)$ is expressed by
$$\ga (\alpha) = f^*s(TY) \cap c_*(\alpha)$$
where $s(TY):= c^{*}(TY)^{-1}$ is the Segre class of the tangent bundle.
\end{thm}
The twisted class $f^*s(TY) \cap c_*(\alpha)$ is called the {\em Ginzburg--Chern class\/} of $\alp$ (\cite{Ginzburg-1, Ginzburg-2} and \cite {Yokura-TAMS2, Yokura-Dedicata}).
Here, the above equality needs a bit of explanation. The left-hand-side $\ga(\alp)$ belongs to the bivariant homology group $\bH(X \xrightarrow {f} Y)$ and the right-hand-side $f^*s(TY) \cap c_*(\alp)$ belongs to the homology group $H^{BM}_*(X)$, and this equality is up to the isomorphism 
$$\begin{CD}
\bH(X \xrightarrow {f} Y) @>\bullet [Y]> \cong > \bH(X \to pt) @> \Cal A > \cong > H^{BM}_*(X)\:,
\end{CD}$$
where the first isomorphism is the bivariant product with the fundamental class $[Y]$ and the second isomorphism $\Cal A$ is the Alexander duality map.  Since we usually identify $\bH(X \to pt)$ as $H^{BM}_*(X)$ via this Alexander duality, we ignore this Alexander duality isomorphism, unless we have to mention it. Hence we have
$$\ga(\alp) \bullet [Y] = f^*s(TY) \cap c_*(\alp).$$
We remark that this formula follows from the {\em simple but crucial} observation that 
$$\ga_f(\alp) \bullet \ga_{Y\to pt}(\jeden_Y) = \ga_{X \to pt}(\alp)$$ 
and the fact that $\ga_{Y \to pt}$ is nothing but the 
Chern--Schwartz--MacPherson class transformation $c_*$.
And in \cite{BSY1} the above theorem is furthermore generalized to the case when the target variety can be singular but is ``like a manifold".

\begin{defn} (cf. \cite{Borho-MacPherson}) Let $A$ be a Noetherian ring. A complex variety $X$ is called an {\it $A$-homology manifold (of dimension $2n$)} or is said to be {\it $A$-smooth} if for all $x \in X$
$$
H_i(X, X \setminus {x}; A) = \begin{cases}
A \quad i = 2n \\
0 \quad \text{otherwise.}
\end{cases}
$$
\end{defn}
In this case $X$ has to be locally pure $n$-dimensional, where we consider $n$ as a locally constant function on $X$.
Just look at the regular part of $X$, because a pure $n$-dimensional complex manifold is a homology manifold of dimension $2n$.
Moreover the local orientation system $or_X$ with stalk
$or_{X,x}= H_{2n}(X, X \setminus {x}; A)\simeq A_X$ is then already trivial 
(on each connected component of $X$) so that $X$ becomes an {\em oriented $A$-homology manifold\/}.

\begin{ex} If $A = \bZ$, a $\bZ$-homology manifold is called simply a {\it homology manifold} (cf. \cite{Milnor-Stasheff}). There are singular complex varieties which are homology manifolds. Such examples are (products of) suitable singular hypersurfaces with isolated singualrities (see \cite{Milnor-hypersurface}). If $A = \bQ$, a $\bQ$-manifold is called a {\it rational homology manifold}. As remarked in \cite[\S1.4 Rational homology manifolds]{Borho-MacPherson}, examples of rational homology manifolds include surfaces with Kleinian singularities, the moduli space for curves of a given genus, and more generally {\it Satake's $V$-manifolds} or {\it orbifolds}. In particular, the quotient of a nonsingular variety by a finite group is a rational homology manifold. 
\end{ex}

\begin{thm}\label{homology-mfd} 
Let $Y$ be a complex analytic variety which is an {\em oriented $A$-homology manifold\/} for some commutative Noetherian ring $A$. If there exists a bivariant Chern class transformation
$\ga :\bF\otimes A  \to \bH\otimes A $, then for any morphism $f: X \to Y$ the bivariant Chern class 
$\ga_f: \bF(X\xrightarrow {f} Y)\otimes A \to \bH( X \xrightarrow {f} Y)\otimes A$ is uniquely determined and it is described by
$$\ga_f(\alp) = f^*c^*(Y)^{-1}\cap c_*(\alp) \:.$$
Here $c^*(Y)$ is the {\em unique cohomology class\/} such that $c_*(1_Y) = c^*(Y) \cap [Y].$ (Note that $c^*(Y)$ is invertible.)
\end{thm} 
When $Y$ is nonsingular, we see that the cohomolgy class $c^*(Y)$ is nothing but the total Chern class $c^*(TY)$ of the tangent bundle 
$TY$, hence the inverse $c^*(Y)^{-1}$ is the total Segre class $s(TY)$. Therefore the twisted class $f^*c^*(Y)^{-1}\cap c_*(\alp)$ shall also be called the {\em Ginzburg--Chern class\/} of $\alp$ and still denoted by $\ga^{\op {Gin}}(\alp)$.
Note that we also have in this more general context the isomorphism
$$\begin{CD}
\bH(X \xrightarrow {f} Y)\otimes A @>\bullet [Y]> \cong > \bH(X \to pt)\otimes A @> \Cal A > \cong > H^{BM}_*(X)\otimes A\:,
\end{CD}$$
since for an oriented $A$-homology manifold $Y$ the fundamental class $[Y]\in H^{BM}_*(X)\otimes A \simeq \bH(X \to pt)\otimes A$ 
is a {\em strong orientation\/} in the sense of bivariant theories (compare \cite{BSY1}).\\ 
$ $\\
{\bf Existence and uniqueness of bivariant characteristic classes.\/}
Note that the proof of Theorem \ref{homology-mfd} also applies in the real (semi-)algebraic or subanalytic context
to a bivariant {\em Stiefel-Whitney class transformation\/}
$\ga :\bF^{mod 2}  \to \bH^{mod 2}$ (with the obvious modification of the notations from $c^*,c_*$ to $w^*,w_*$).
In a similar manner, we can show the following theorem, which is an extended version of \cite[Theorem (3.7)]{Yokura-Documenta}:
\begin{thm} 
The Grothendieck transformation from the bivariant algebraic K-theory $\bK_{\op {alg}}$ of perfect complexes
$$\tau : \bK_{\op {alg}} \to \bH_{\bQ}$$
constructed in \cite[Part II]{Fulton-MacPherson} is unique on morphisms whose target varieties are rational homology manifolds. Explicitly, for a bivariant element $\alp \in \bK_{\op {alg}}(X \xrightarrow {f} Y)$ with $Y$ being a rational homology manifold
$$\tau(\alp) = f^*td^*(Y)^{-1} \cap td_*(\alp \bullet [\Cal O_Y]).$$
Here $[\Cal O_Y] \in \bK_{\op {alg*}}(Y) \simeq \bold G_{0}(Y)$ is the class of the structure sheaf and
the associated covariant transformation $\tau_*: \bold \bK_{\op {alg}*} (\quad)\simeq \bold G_{0}(\quad)  \to H^{BM}_*(\quad; \bQ)$ is Baum--Fulton--MacPherson's Todd class transformation $td_{*}$ constructed in \cite{Baum-Fulton-MacPherson}. Moreover 
$td^*(Y) \in H^*(Y; \bQ)$ is the Poincar\'e dual of the Todd class $td_*(Y):=td_*([\Cal O_Y])$, which is invertible.
\end{thm}

Conversely we ask ourselves whether the above Ginzburg--Chern class  becomes a Gro- thendieck transformation for morphisms whose target varieties are oriented $A$ - homology manifolds.
\begin{thm}\label{Gin}
 For a morphism of complex analytic varieties $f: X \to Y$ with $Y$ an oriented $A$-homology  manifold, we define $ \overline{\bF} (X \xrightarrow {f} Y)$ to be the set of all constructible functions $\alp \in F(X)$ satisfying the following two conditions ($\sharp$) and ($\flat$) : for any fiber square
$$\CD X' @> g' >> X \\@V f' VV @VV f V\\Y' @> g >> Y, \endCD$$
with $Y'$ an oriented  $A$-homology  manifold the following equalities hold:\\
($\sharp$)  for any constructible function $\be' \in F(Y')$:
$$\ga^{\op {Gin}} (g^{\star}\alp \bullet \be') = \ga^{\op {Gin}}(g^{\star}\alp) \bullet \ga^{\op {Gin}}(\be'),$$
($\flat$) 
$$\ga^{\op {Gin}}(g^{\star}\alp) = g^{\star} \ga^{\op {Gin}}(\alp).$$
\noindent Then $\overline {\bF}$ becomes a bivariant theory with the same operations as in $s\bF$ and furthermore the transformation 
$$\ga^{\op {Gin}}: \overline {\bF} \to \bH$$
is well-defined and becomes the unique Grothendieck transformation satisfying that $\ga^{\op {Gin}}$ for morphisms to a point is the Chern--Schwartz--MacPherson class transformation $c_*: F \to H_*$. And also $\overline {\bF}(X \to pt) = F(X)$. 
\end{thm}
The proof of the theorem is the same as in \cite{Yokura-CEJM}, in which the case when the target variety $Y$ is nonsingular is treated. Note that to prove $\overline {\bF}(X \to pt) = F(X)$ we need the {\em cross product formula\/} or {\em multiplicativity\/}
of the Chern--Schwartz--MacPherson class transformation $c_*$ due to {\em Kwieci\'nski\/} \cite{Kwiecinski1} (cf. \cite{KY}),
i.e. the commutativity of the following diagram:
$$\begin{CD} F(X)\times F(Y) @> \times >> F(X\times Y)\\
@V c_*\times c_* VV @VV c_*  V \\
H^{BM}_{*}(X;\bZ)\times H^{BM}_{*}(Y;\bZ) @> \times >> H^{BM}_{*}(X\times Y;\bZ) \:.
\end{CD} $$
The cross product formula for Stiefel-Whitney classes in the {\em real algebraic context\/} can be shown similarly
by using ``resolution of singularities", or the corresponding product formula for ``characteristic cycles" of constructible functions
so that a variant of this theorem also works in the real algebraic context.

And for a {\em much more general version \/}
of Theorem \ref{Gin}, see \cite{Schuermann-partial GT}.\\

The above theorem led us to another {\em uniqueness theorem\/}, which in a sense gives a positive solution to the general uniqueness problem concerning Grothendieck transformations posed in \cite[\S 10 Open Problems]{Fulton-MacPherson}.
\begin{thm}\label{bivariant} 
We define
$$\widetilde {\bF}(X \xrightarrow {f} Y)$$
to be the set consisting of all $\alp \in s\bF(X \xrightarrow {f} Y)$ satisfying the following condition: there exists a bivariant class $B_{\alp} \in \bH(X \xrightarrow {f}  Y)$ such that for any base change $g: Y' \to Y$ (without any requirement) of an independent square
$$\CD X'@> {g'} >> X \\@V {f'} VV @VV f V\\Y' @> g>> Y, \endCD$$
and for any $\be' \in F(Y')$ the following equality holds:
$$c_*(g^*\alp \bullet \be') = {g}^*B_{\alp} \bullet c_*(\be').$$
Then $\widetilde {\bF}$ is a bivariant theory. Furthermore $\widetilde {\bF}(X \to pt) = F(X)$.
\end{thm}
The above bivariant class $B_{\alp}$ should ideally be the unique bivariant Chern class of $\alp$. However, so far we still do not know if it is the case or not. So, provisionally we call $B_{\alp}$ {\it a pseudo-bivariant Chern class of $\alp$}.

\begin{ex}[VRR for smooth morphisms] \label{smooth}
Let $f: X\to Y$ be a {\em smooth} morphism of possible singular varieties. Then we have
$$\jeden_{f}:=1_{X} \in \widetilde {\bF}(X \xrightarrow {f} Y)$$
with $c^{*}(T_{f})\bullet [f]$ being a pseudo-bivariant Chern class of $\jeden_{f}$.
Here $T_{f}$ is the vector bundle of tangent spaces of fibers of $f$, and $[f]\in \bH(X \xrightarrow {f}  Y)$
is the {\em canonical orientation\/} of the smooth morphism $f$. Then as in 
Theorem \ref{bivariant} we have for $\be' \in F(Y')$:
$$
\aligned
c_*(g^*\jeden_{f} \bullet \be') & = c_*(f'^* \be') \\
& = c^{*}(T_{f'})\cap f'^{!}c_{*}(\be') \\
& = c^{*}(T_{f'})\bullet [f']\bullet c_{*}(\be') \\
& =g^{*}c^{*}(T_{f})\bullet g^{*}[f]\bullet c_{*}(\be') \\
& = g^{*}(c^{*}(T_{f})\bullet [f])\bullet c_{*}(\be').  
\endaligned
$$
Here $f'^{!}=[f'] \bullet : H^{BM}_{*}(Y')\simeq \bH_*(Y')\to \bH_*(X')\simeq H^{BM}_{*}(X')$
is the {\em smooth pullback\/} in Borel-Moore homology, and the equality
\begin{equation} \label{VRR2}
c_*(f'^* \be') = c^{*}(T_{f'})\cap f'^{!}c_{*}(\be')
\end{equation}
is the so-called {\em Verdier-Riemann-Roch theorem\/} for the smooth morphism $f'$
and the Chern class transformation $c_*$
(compare \cite{Fulton-MacPherson, Schuermann-VRR, Yokura-VRR-Chern}).
\end{ex}

In order to remedy this unpleasant possible non-uniqueness of the bivariant class $B_{\alp}$ above, we set
$$\aligned& \bPH(X \xrightarrow {f}  Y) := \\
& \Bigl  \{ B \in \bH(X \xrightarrow {f}  Y) | \text {$B$ is a pseudo-bivariant Chern class of some 
$\alp \in \widetilde {\bF}(X \xrightarrow {f}  Y)$ } \Bigr \}\endaligned$$
to be the set of all pseudo-bivariant Chern classes for the morphism $ f: X \to Y$.
It is clear that $\bPH$ is a {\em bivariant subtheory\/} of $\bH$, i.e, it is a subgroup stable under the three bivariant operations.Then we define
$$\widetilde {\bH}(X \xrightarrow {f}  Y) := \bPH(X \xrightarrow {f}  Y) / \sim$$
where the relation $\sim$ is defined by
$$B \sim B' \Longleftrightarrow g^*B \bullet c_*(\be') = g^*B' \bullet c_*(\be')$$
for all independent squares with $g:Y' \to Y$ and all $\be' \in F(Y')$. Certainly the relation $\sim$ is an equivalence relation. In other words, with this identification we want to make possibly many pseudo-bivariant Chern classes into one unique bivariant Chern class. Indeed we have
\begin{thm} $\widetilde {\bH}(X \xrightarrow {f}  Y)$ is an Abelian group and $\widetilde {\bH}$ is a bivariant theory with the canonical operations induced from those of $\bH$. Furthermore we have
$$\widetilde {\bH}(X \to {pt}) = \op {Image}(c_*:F(X) \to H^{BM}_*(X)).$$
\end{thm}
And we have the following theorem
\begin{thm}\label{uniqueness of GT}
There exists a unique Grothendieck transformation
$$\widetilde {\ga} : \widetilde {\bF} \to \widetilde {\bH}$$
whose associated covariant transformation is $c_*:F \to \op {Im}(c_*)$, where 
$$\op {Im}(c_*)(X) := \op {Image}\Bigl (c_*: F(X) \to H^{BM}_*(X) \Bigr )\:.$$
\end{thm}

\begin{rem} 
As mentioned above, a key for the above argument is the fact that $c_*(\alp) = \ga(\alp) \bullet c_*(\jeden_Y)$. So, very sloppy speaking, the bivariant class $\ga(\alp)$ is a kind of ``$c_*(\alp)$ divided by $c_*(\jeden_Y)$", whatever it is meant to be. In our previous paper \cite{Yokura-Documenta} we posed the problem of whether or not there is a reasonable bivariant homology theory so that such a ``quotient"
$$\frac {c_*(\alp)}{c_*(\jeden_Y)}$$
is well-defined. The above theory $\widetilde {\bH}$ is in a sense a positive answer to this problem.
\end{rem}

The above construction works for the following more general situation such that\\
\noindent (1) there exists a natural transformation $\tau_{*}: F_*(X) \to H_*(X)$ between two covariant functors $F_*$ and $H_*$
(covariant with respect to proper maps) such that $F_*(pt)$ and $H_*(pt)$ are commutative rings with unit and such that $\tau_*$ maps the unit to the unit,\\
\noindent (2) there are two bivariant theories $\bF$ and $\bH$ such that the associated covariant theories are
$$\bF(X \to pt) = F_*(X) \quad \text {and} \quad  \bH(X \to pt) = H_*(X),$$
\noindent (3) $\tau_{*}$ commutes with the bivariant exterior products, i.e., the following diagram commutes
$$\CD F_*(X) \times F_*(Y) @>  \times >> F_*(X \times Y) \\
@V \tau_{*} \times \tau_{*} VV @VV \tau_{*} V\\
H_*(X) \times H_*(Y) @> \times  >> H_*(X \times Y). \endCD$$
Here we assume that for $X = Y = \{pt\}$ a point this exterior product agrees with the given ring structure.\\

Certainly this construction works for the previous {\em motivic Chern class transformation\/}
$$mC_*: K_0(\Cal V/\quad) \to \bold G_0(\quad)\otimes \bZ[y] $$
and the {\em motivic Hirzebruch class transformation\/}
$${T_y}_*: K_0(\Cal V/\quad) \to H_*(\quad)\otimes \bQ[y]. $$
Indeed, the bivariant theory for  $K_0(\Cal V/\quad)$ is the simple bivariant theory
$$s\bold K_0(X \to Y) := K_0(\Cal V/X)\:,$$
the bivariant theory for $\bold G_0(\quad)\otimes \bZ[y]$ is the Fulton--MacPherson's bivariant algebraic K-theory $\bK_{\op {alg}}$
 tensored with $\bZ[y]$, 
and the bivariant theory for $H_*(\quad)\otimes \bQ[y]$ is of course the Fulton--MacPherson's bivariant homology theory $\bH$ tensored with $\bQ[y]$. It also applies in the real algebraic context to the {\em Stiefel-Whitney class transformation\/}
$$w_{*}: F^{mod 2}(\quad) \to H_*(\quad; \bZ_2) $$
by using the simple bivariant theory $s \bF^{mod 2}$ of $\bZ_{2}$-valued real algebraically constructible functions. 

\begin {rem}
Let $f: X\to Y$ be a {\em smooth} morphism of possible singular varieties. 
Then also the example \ref{smooth} works in this context, with 
$$\jeden_{f}:=\jeden_{X}=[id_{X}] \in s\bold K_0(X \xrightarrow {f} Y) \quad \text{or} \quad
\jeden_{f}:=1_{X}\in s \bF^{mod 2}(X \xrightarrow {f} Y) \:,$$
and $c \ell^{*}(T_{f})\bullet [f]$ being a pseudo-bivariant class of $\jeden_{f}$ for $c \ell^{*}(T_{f})=\lambda_{y}(T_{f}^*),
\widetilde td_{(y)}(T_{f})$ or $w^*(T_{f})$. Here the corresponding {\em Verdier-Riemann-Roch theorem\/} for the smooth morphism $f'$
follows for the motivic characteristic classes $mC_*$ and $T_{y*}$ from \cite[Corollary 2.1 and  Corollary 3.1]{BSY2}.
For the Stiefel-Whitney class transformation $w_*$ it can be shown as for Chern classes by using ``resolution of singularities"
or ``characteristic cycles of constructible functions".
\end{rem}

This {\em Verdier-Riemann-Roch theorem for smooth morphisms\/} is also very important for the definition of
{\em G-equivariant characteristic class transformations\/} in the equivariant algebraic context with $G$ a reductive linear algebraic
group. Here we refer to \cite{EG1, EG2, BZ} for the {\em equivariant Todd class transformation\/} $td_{*}^G$, and to \cite{Oh}
for the {\em equivariant Chern class transformation\/} $c_*^G$. In fact, in future work we will construct in this
equivariant algebraic context equivariant versions $mC_*^G$ and $T_{y*}^G$ of our motivic characteristic classes,
together with the equivariant version of Theorem \ref{T_y_*}, relating  $T_{-1*}^G$ with $c_*^G$ and
 $T_{0*}^G$ with $td_{*}^G$.\\
$ $\\
{\bf Bivariant L-classes.\/}
At the moment we have no bivariant version with values in bivariant homology
of the L-class transformation
$$L_{*}: \Omega(X)  \to H_{*}(X,\bQ) \:,$$
since we do not know a suitable bivariant theory, whose associated covariant theory reduces to the cobordism group
$\Omega(\quad)$ of selfdual constructible sheaf complexes. Note that in this case we {\em cannot\/} define a
{\em simple bivariant theory\/} $s\Omega$. Of course the Grothendieck group of constructible sheaf complexes
$K_c(\quad)$ satisfies the properties (SB1-4) with respect to the induced proper push down $f_{*}$, pullback $f^{*}$
and tensor product $\otimes$ so that one gets a simple bivariant theory $s \bold K_{c}$. But the problem is that
$f^{*}$ and  $\otimes$ do {\em not\/} commute with duality in general so that this approach doesn't apply
to $\Omega(\quad)$.\\

A similar problem appears in the context of real semialgebraic and subanalytic geometry for the group $F_{Eu}^{mod 2}(\quad)$ of 
$\bZ_{2}$-valued constructible functions satisfying the {\em mod 2 local Euler condition\/} (for a constant map),
which also can be interpretated as a ``duality" condition (compare \cite[p.135 and Remark 5.4.4, p.367]{Schuermann-book}).
This group (or condition) is also not stable under general pullback or product so that one {\em cannot\/} define a simple bivariant theory $s \bF_{Eu}^{mod 2}$ in this context (compareable to $s \bF^{mod 2}$ in the real algebraic context).
Nevertheless one can define a {\em Stiefel-Whitney class transformation\/}
$$w_*: F_{Eu}^{mod 2}(\quad) \to H^{BM}_*(\quad;\bZ_{2})$$
with the help of ``characteristic cycles of constructible functions" (compare \cite{FuMC}),
which is {\em multiplicative for exterior products\/} and satisfies the {\em Verdier Riemann-Roch theorem
for smooth morphisms\/}.\\

Similarly one can define in the complex algebraic or analytic context an {\em exterior product and smooth pullback\/}
for the cobordism group $\Omega(\quad)$ of selfdual constructible sheaf complexes (compare \cite{BSY2}),
and the L-class transformation $L_{*}$ is also {\em multiplicative\/} by an argument similarly as
in the recent paper \cite[p.26, Proposition 5.16]{Woolf}.
Also the corresponding Verdier Riemann-Roch theorem for smooth morphisms seems reasonable,
but at the moment we have no proof or reference for this.
Of course this VRR theorem holds on the image of the transformation $\omega: K_0(\Cal V/\quad) \to \Omega (\quad)$
from Theorem \ref{T_y_*} (compare \cite{BSY2}).\\

Then in both these cases, L-class and Stiefel-Whitney class transformations, we can apply the results
of \cite{Yokura-TopAppl} to get at least bivariant versions of these theories
for the corresponding {\em operational bivariant theories\/}.

\section{Characteristic classes of proalgebraic varieties}\label{section-pro}
A {\em pro-algebraic variety\/} is defined to be a projective system of complex algebraic varieties and a {\em proalgebraic variety\/} is defined to be the projective limit of a pro-algebraic variety. Proalgebraic varieties are the main objects in \cite{Gromov1}. A pro-category was first introduced by A. Grothendieck \cite{Grothendieck1} and it was used to develope the Etale Homotopy Theory \cite{Artin-Mazur} and Shape Theory (e.g., see \cite{Borusk}, \cite{Mardesic-Segal}, etc.) and so on. In [Grom 1] M. Gromov investigated the {\em surjunctivity\/}, i.e. being either surjective or non-injective, in the category of proalgebraic varieties. The original or classical surjunctivity theoremis the so-called {\em Ax' Theorem\/} \cite{Ax}, saying that every regular selfmapping of a complex algebraic variety is surjunctive; thus if it is injective then it has to be surjective.\\

A very simple example of a proalgebraic variety is the Cartesian product $X^{\bN}$ of countable infinitely many copies of a complex algebraic variety $X$, which is one of the main objects treated in \cite{Gromov1}. Then, what would be the {\em ``Chern--Schwartz--MacPherson class"\/} of $X^{\bN}$ ? In particular, what would be the {\em  ``Euler--Poincar\'e characteristic"\/} of $X^{\bN}$ ?
This simple question led us to a study of characteristic classes of proalgebraic varieties and it naturally led us to the so-called {\em motivic measures\/} (see \cite{Yokura-Chern motivic, Yokura-characteristic motivic}). The motivic measures/integrations have been actively studied by many people (e.g., see \cite{Craw}, \cite{Denef-Loeser1}, \cite{Denef-Loeser2}, \cite{Kontsevich}, \cite{Looijenga}, \cite{Veys-stringy} etc.).\\

In a general set-up one can deal with the so-called {\em bifunctors\/}. The bifunctors which we consider are bifunctors $\Cal F: \Cal C \to \Cal A$ from a category $\Cal C$ to the category $\Cal A$ of abelian groups, i.e., $\Cal F$ is a pair $(\Cal F_*, \Cal F^*)$ of a {\em covariant functor\/} $\Cal F_*$and a {\em contravariant functor\/} $\Cal F^*$ such that $\Cal F_*(X) = \Cal F^*(X)$ for any object $X$. Unless some confusion occurs, we just denote $\Cal F(X)$ for $\Cal F_*(X) = \Cal F^*(X)$. A typical example is the constructible function functor $F(X)$. Furthermore we assume that for a final object $pt \in Obj(\Cal C)$, $\Cal F(pt)$ is a commutative ring $\Cal R$ with a unit.  The morphism from an object $X$ to a final object $pt$ shall be denoted by $\pi_X : X \to pt$. Then the covariance of the bifunctor $\Cal F$ induces the homomorphism 
$\pi_{X*}:=\Cal F(\pi_X) : \Cal F(X) \to \Cal F(pt) = \Cal R$, which shall be denoted by
$$\chi_{\Cal F} : \Cal F(X) \to \Cal R$$
and called the {\em $\Cal F$-characteristic\/}, just mimicking the Euler--Poincar\'e characteristic 
(with compact support) $\chi: F(X) \to \bZ$ in the case when $\Cal F = F$.\\

Let $X_{\infty} = \varprojlim_{\la \in \La} \Bigl \{ X_{\la}, \pi_{\la \mu}: X_{\mu} \to X_{\la} \Bigr \}$ be a proalgebraic variety. Then we define
$$\Cal F^{\op {ind}}(X_{\infty}) := \varinjlim _{\la \in \La} \Bigl \{\Cal F(X_{\la}), {\pi_{\la \mu}}^*: \Cal F(X_{\la}) \to \Cal F(X_{\mu}) (\la < \mu) \Bigr \},$$
which {\it may not belong to} the category $\Cal A$. Another finer one can be defined as follows. Let $P = \bigl \{p_{\la \mu} \bigr \}$ be a {\em projective system\/} of elements of $\Cal R$ by the directed set 
$\La$, i.e., a set such that $p_{\la \la} = 1 \quad \text {(the unit) }$ and 
$p_{\la \mu} \cdot p_{\mu \nu} = p_{\la \nu} \quad (\la < \mu < \nu).$ For each $\la \in \La$ the {\em subobject\/} 
$\Cal F^{\op {st}}_P(X_{\la})$ of {\em $\chi_{\Cal F}$-stable elements\/} in $\Cal F(X_{\la})$ is defined to be 
$$\Cal F^{\op {st}}_P(X_{\la}) : = \Bigl \{ \alp_{\la} \in \Cal F(X_{\la}) | \;\chi_{\Cal F}\bigl ({\pi_{\la \mu}}^*\alp_{\la} \bigr ) = p_{\la \mu}\cdot \chi_{\Cal F}(\alp_{\la}) \; \text {for any $\mu$ such that} \; \la <\mu \Bigr \}.$$
The {\em inductive limit\/} 
$$\varinjlim_{\La}\Bigl \{\Cal F^{\op {st}}_P(X_{\la}) , \quad {\pi_{\la \mu}}^*: 
\Cal F^{\op {st}}_P(X_{\la}) \to \Cal F^{\op {st}}_P(X_{\mu}) \quad (\la < \mu) \Bigr \}$$ 
considered for a proalgebraic variety 
$X_{\infty} = \varprojlim _{\la \in \La} X_{\la}$ is denoted by 
$$\Cal F^{\op {st.ind}}_P(X_{\infty}).$$
Of course this definition is not intrinsic to the proalgebraic variety $X_{\infty}$,
but depends on the given projective system $\Bigl \{ X_{\la}, \pi_{\la \mu}: X_{\mu} \to X_{\la} \Bigr \}$. 
But for simplicity we use this notation. Our key observation, which is an application of standard facts on indcutive systems and limits, is the following:
\begin{thm}\label{pro-calF}
(i) For a proalgebraic variety $X_{\infty} = \varprojlim_{\la \in \La} \Bigl \{ X_{\la}, \pi_{\la \mu}: X_{\mu} \to X_{\la} \Bigr \}$  and a projective system $P = \bigl \{p_{\la \mu} \bigr \}$ of elements of $\Cal R$, we have the homomorphism
$$\chi^{\op {ind}}_{\Cal F}: \Cal F^{\op {st.ind}}_P(X_{\infty}) \to \varinjlim _{\la \in \La} \Bigl \{ \times p_{\la \mu} : \Cal R \to \Cal R \Bigr \},$$
which is called the {\em proalgebraic $\Cal F$-characteristic homomorphism\/}.\\
(ii) Assume $\La = \bN$. For a proalgebraic variety 
$X_{\infty} = \varprojlim_{n \in \bN} \Bigl \{ X_n, \pi_{nm}: X_m \to X_n \Bigr \}$ and a projective system $P = \{p_{nm} \}$ of elements of $\Cal R$, the proalgebraic $\Cal F$-characteristic homomorphism $\chi^{\op {ind}}_{\Cal F}: 
\Cal F^{\op {st.ind}}_P(X_{\infty}) \to \varinjlim _n \Bigl \{ \times p_{nm} : \Cal R \to \Cal R \Bigr \}$ is realized as the homomorphism
$$\widetilde {\chi^{\op {ind}}_{\Cal F}}: \Cal F^{\op {st.ind}}_P(X_{\infty}) \to \Cal R_P$$
defined by
$$\widetilde {\chi^{\op {ind}}_{\Cal F}}\Bigl ([\alp_n] \Bigr ) := \frac {\chi_{\Cal F}(\alp_n)}{p_{01}\cdot p_{12} \cdot p_{23} \cdots p_{(n-1)n}}.$$
Here $p_{01}:= 1$ and $\Cal R_P$ is the ring $\Cal R_S$ of fractions of $\Cal R$ with respect to the multiplicatively closed set $S$ consisting of all the finite products of powers of elements in $P$.\\
(iii) In particular, in the case when the above projective system $P = \{p^s \}$ consists of powers of an element $p$, we get the homomorphism
$$\widetilde {\chi^{\op {ind}}_{\Cal F}}: \Cal F^{\op {st.ind}}_P(X_{\infty}) \to \Cal R \Bigl [\frac {1}{p} \Bigr ]$$
defined by
$$\widetilde {\chi^{\op {ind}}_{\Cal F}} \Bigl ([\alp_n] \Bigr ) := \frac {\chi_{\Cal F}(\alp_n)}{p^{n-1}}.$$
Here $\Cal R \Bigl [\frac {1}{p} \Bigr ]$ is the 
localization by the multiplicatively closed set $S:=\{p^s|s\in \bN_{0}\}$.
\end{thm}

Note that $\Cal R_S$ or $\Cal R \Bigl [\frac {1}{p} \Bigr ]$ is the zero ring in the case when
$0\in S$ for the corresponding muliplicatively closed set $S$.
A typical example for the above theorem is the following.
\begin{ex}\label{typical example} Let $X_{\infty} = \varprojlim_{n \in \bN} \Bigl \{ X_n, \pi_{nm}: X_m \to X_n \Bigr \}$ be a proalgebraic variety such that for each $n$ the structure morphism $\pi_{n,n+1}: X_{n+1} \to X_n$ satisfies the condition that the Euler--Poincar\'e characteristics of the fibers of $\pi_{n,n+1}$ are non-zero (which implies the surjectivity of the morphism $\pi_{n,n+1}$) and constant; for example,  $\pi_{n,n+1}: X_{n+1} \to X_n$ is a locally trivial fiber bundle with fiber variety being $F_n$ and $\chi (F_n) \not = 0$
Let us denote the constant Euler--Poincar\'e characteristic of the fibers of the morphism $\pi_{n,n+1}: X_{n+1} \to X_n$ by $e_n$ and we set $e_0 :=1$. Then we get the canonical proalgebraic Euler--Poincar\'e characteristic homomorphism
$$\chi^{\op {ind}}: F^{\op {ind}}(X_{\infty}) \to \bQ$$
described by
$$\chi^{\op {ind}}\left ([\alp_n] \right ) = \frac {\chi(\alp_n)}{e_0 \cdot e_1 \cdot e_2 \cdots e_{n-1}}. $$
In particular, if the Euler--Poincar\'e characteristics $e_n$ are all the same, say $e_n = e$ for any $n$, then the canonical proalgebraic Euler--Poincar\'e characteristic homomorphism
$\chi^{\op {ind}}: F^{\op {ind}}(X_{\infty}) \to \bQ$
is described by
$\displaystyle \chi^{\op {ind}}\left ([\alp_n] \right ) = \frac {\chi(\alp_n)}{e^{n-1}}$, and furthermore the target ring $\bQ$ can be replaced by the ring $\bZ \Bigl [\frac {1}{e} \Bigr ]$.
\end{ex}

Note that this example applies especially to the Cartesian product $X^{\bN}$ of countable infinitely many copies of a complex algebraic variety $X$ with $\chi(X)\neq 0$.
In fact this example of Cartesian products is a special case of the following more general example:
\begin{ex}\label{typical example2}
We make the following additional assumptions for our bifunctor:\\

(1) The contravariant functor $\Cal F^*$ takes values in the category of {\em commutative rings with
unit\/}. The corresponding unit in $\Cal F (X)$ is denoted by $\jeden_X$, and
$\Cal F (X)$ becomes an $\Cal R:=\Cal F (pt)$-algebra by the pullback for $\pi_{X}: X\to pt$.\\

(2) $\Cal F^*$ and $\Cal F_*$ are related for a morphism $f:X\to Y$ by the {\em projection formula\/}
$$f_*(\alp \cdot f^*\be)=f_*(\alp)\cdot \be \quad \text{$\quad$ for all
$\alp \in \Cal F( X)$ and $\be \in \Cal F( Y)$}$$
so that $f_*: \Cal F( X)\to \Cal F( Y)$ is $\Cal F (Y)$- and $\Cal R$-linear.
(This is just a special case of our simple bivariant theories, where all morphisms
are ``proper" and only the ``trivial fiber squares" are ``independent".)

Consider a proalgebraic variety $X_{\infty} = \varprojlim_{n \in \bN} \Bigl \{ X_n, \pi_{nm}: X_m \to X_n \Bigr \}$ such that for each $n$ the structure morphism $\pi_{n,n+1}: X_{n+1} \to X_n$ satisfies the condition 
$$\pi_{n,n+1*}(\jeden_{X_{n+1}})=e_{n}\cdot \jeden_{X_{n}}\in \Cal F( X_{n})
\quad \text{for some $e_{n}\in \Cal R$, with $e_{0}:=\jeden_{pt}$.}$$ 
Then we get the canonical proalgebraic $\Cal F$-characteristic homomorphisms
$$\chi^{\op {ind}}_{\Cal F,X_1}: \Cal F^{\op {ind}}(X_{\infty}) \to \Cal F (X_{1})_E
\quad \text{and} \quad
\chi^{\op {ind}}_{\Cal F}: \Cal F^{\op {ind}}(X_{\infty}) \to \Cal R_E$$
described by
$$\chi^{\op {ind}}_{\Cal F,X_1}\left ([\alp_n] \right ) = \frac {{\pi_{1,n}}_*(\alp_n)}{e_0 \cdot e_1 \cdot e_2 \cdots e_{n-1}}
\quad \text{and} \quad
\chi^{\op {ind}}_{\Cal F}\left ([\alp_n] \right ) = \frac {\chi(\alp_n)}{e_0 \cdot e_1 \cdot e_2 \cdots e_{n-1}}. $$
Here $\Cal R_E$ (or $\Cal F (X_{1})_E$) is the ring of fractions of $\Cal R$ with respect to the multiplicatively closed set consisting of all the finite products of powers of the elements $e_{i}$ (or their pullbacks to $X_1$).
\end{ex}
Consider a bifunctor as in example \ref{typical example2}, with
$f:X \to Y$ being a morphism such $f_*(\jeden_{X})=e_{f}\cdot \jeden_{Y}$
for some $e_f\in \Cal R$. Then one gets any $\alp \in \Cal F( Y)$:
$$f_*f^*\alp = f_*(\jeden_{X}\cdot f^*\alp)= e_f \cdot \alp$$
so that for any morphism $g:Y\to Z$ (e.g. $g=\pi_{Y}: Y\to pt$):
$$
\aligned
(g\circ f)_* \bigl (f^* \alp \bigr ) & = g_* \bigl (f_*f^* \alp \bigr ) \\
 & = g_* (e_f \cdot \alp) \\
& = e_f \cdot g_* (\alp). 
\endaligned
$$
Hence if we set in the context of the example
\begin{equation*}
p_{nm} = \begin{cases}
1 & \quad n = m \\
e_n \cdot e_{n+1} \cdots e_{m-1} & \quad n<m,
\end{cases}
\end{equation*}
then $P:= \{p_{nm} \}$ is a projective system and  $\Cal F^{\op {st.ind}}_P (X_{\infty}) = 
\Cal F^{\op {ind}} (X_{\infty})$ for both notions of Euler characteristics working over the base space $X_1$ or over $pt$.
Thus the above description of $\chi_{\Cal F,X_1}^{\op {ind}}$ and
$\chi_{\Cal F}^{\op {ind}}$ follows from Theorem \ref{pro-calF}.\\

A ``motivic" version of the Euler--Poincar\'e characteristic $\chi: F(X) \to \bZ$ is the homomorphism $\Ga_X: F(X) \to K_0(\Cal V/X)$ ``tautologically" defined by 
$$\Ga_X (\sum_W a_W\jeden_W) := \sum_W a_W[W\subset X] \:,$$ 
or better is the composite $\Ga:=\pi_{X*}\circ \Ga_X: F(X) \to K_0(\Cal V)$.
Note that $\Ga_X$ commutes with pullback $f^*$ (but not with push down $f_*$).
Then we get the following theorem, which is a generalization of the (naive) motivic measure:
\begin{thm}\label{pro-Gro}
(i) For a proalgebraic variety $X_{\infty} = \varprojlim_{\la \in \La} \Bigl \{ X_{\la}, \pi_{\la \mu}: X_{\mu} \to X_{\la} \Bigr \}$ and a projective system $G = \bigl \{\ga_{\la \mu} \bigr \}$ of Grothendieck classes, we get the {\em proalgebraic Grothendieck class homomorphism\/}
$$\Ga^{\op {ind}}: F^{\op {st.ind}}_G(X_{\infty}) \to \varinjlim _{\la \in \La} \Bigl \{ \times \ga_{\la \mu}: K_0(\Cal V) \to 
K_0(\Cal V) \Bigr \}.$$
(ii) Assume $\La = \bN$. For a proalgebraic variety 
$X_{\infty} = \varprojlim_{n \in \bN} \Bigl \{ X_n, \pi_{nm}: X_m \to X_n \Bigr \}$ and a projective system $G = \{\ga_{n,m} \}$ of  Grothendieck classes, we have the following canonical proalgebraic Grothendieck class  homomorphism
$$\widetilde {\Ga^{\op {ind}}}: F^{\op {st.ind}}_G(X_{\infty}) \to K_0(\Cal V)_G$$
which is defined by
$$\widetilde {\Ga^{\op {ind}}}\Bigl ([\alp_n] \Bigr ) := \frac {\Ga(\alp_n)}{\ga_{0 1}\cdot \ga_{1 2} \cdot \ga_{2 3}\cdots 
\ga_{(n-1) n}}.$$
Here we set $\ga_{0 1}:= \jeden$ and $K_0(\Cal V)_G$ is the ring of fractions of $K_0(\Cal V)$ with respect to the multiplicatively closed set consisting of finite products of powers of elements of $G$.\\
(iii) Let $X_{\infty} = \varprojlim_{n \in \bN} \Bigl \{ X_n, \pi_{nm}: X_m \to X_n \Bigr \}$ be a proalgebraic variety such that each structure morphism $\pi_{n,n+1}: X_{n+1} \to X_n$ satisfies the condition:
$$\pi_{n,n+1*}([id_{X_{n+1}}]) = \ga_n \cdot [id_{X_n}]\in K_0(\Cal V/X_n)
\quad \text{for some $\ga_n \in K_0(\Cal V)$};$$ 
for example $\pi_{n,n+1}: X_{n+1} \to X_n$ is a {\em Zariski locally trivial fiber bundle with fiber variety being $F_n$\/} (in which case one can take $\ga_n := [F_n] \in K_0(\Cal V)$). Then the canonical proalgebraic Grothendieck class homomorphisms
$$\Ga^{\op {ind}}_{X_1}: F^{\op {ind}}(X_{\infty}) \to K_0(\Cal V/X_1)_G
\quad \text{and} \quad
\Ga^{\op {ind}}: F^{\op {ind}}(X_{\infty}) \to K_0(\Cal V)_G
$$
are described by
$$\Ga^{\op {ind}}_{X_1}\left ([\alp_n] \right ) = \frac {\pi_{1,n*}(\Ga_{X_n} (\alp_n))}
{\ga_0 \cdot \ga_1 \cdot \ga_2 \cdots \ga_{n-1}} 
\quad \text{and} \quad
\Ga^{\op {ind}}\left ([\alp_n] \right ) = \frac {\Ga (\alp_n)}{\ga_0 \cdot \ga_1 \cdot \ga_2 \cdots \ga_{n-1}} \:.
$$
Here $\ga_0 :=\jeden$ and $K_0(\Cal V)_G$ (or $K_0(\Cal V/X_1)_G$)
is  the ring of fractions of $K_0(\Cal V)$ with respect to the multiplicatively closed set consisting of finite products of powers of $\ga_m$ ($m =1, 2, 3 \cdots)$ (or their pullbacks to $X_1$).\\
(iv) In particular, if $\ga_n= \ga$ for all $n$, then the canonical proalgebraic  Grothendieck class homomorphisms
$$\Ga^{\op {ind}}_{X_1}: F^{\op {ind}}(X_{\infty}) \to K_0(\Cal V/X_1)_G
\quad \text{and} \quad
\Ga^{\op {ind}}: F^{\op {ind}}(X_{\infty}) \to K_0(\Cal V)_G
$$
are described by
$$\Ga^{\op {ind}}_{X_1}\left ([\alp_n] \right ) =\frac {\pi_{1,n*}(\Ga_{X_n} (\alp_n))}
{\ga^{n-1}}
\quad \text{and} \quad
\Ga^{\op {ind}}\left ([\alp_n] \right ) = \frac {\Ga (\alp_n)}{\ga^{n-1}}. $$
In this special case the quotient ring $K_0(\Cal V)_G$ (or $K_0(\Cal V/X_1)_G$)
shall be simply denoted by $K_0(\Cal V)_{\ga}$ (or $K_0(\Cal V/X_1)_{\ga}$).
\end{thm}

\begin{ex} The arc space $\Cal L(X)$ of an algebraic variety $X$ is
 defined to be the projective limit of the projective system consisting of the truncated arc varieties $\Cal L_n(X)$ of jets of order $n$ together with the canonical  projections $\pi_{n,n+1}: \Cal L_{n+1}(X) \to \Cal L_n(X)$. 
Note that $\Cal L_0(X)=X$ so that this time we use $\La = \bN_0$.
Thus the arc space is a nontrivial example of a proalgebraic variety. If $X$ is {\em nonsingular\/} and of complex dimension $d$, then the projection $\pi_{n,n+1}: \Cal L_{n+1}(X) \to \Cal L_n(X)$ is a Zariski locally trivial fiber bundle with fiber being $\bC^d$. Thus in this case, in (iv) of Theorem \ref{pro-Gro} the Grothendieck class $\ga$ is $\bL^d$, with $\bL:=[\bC]$.
\end{ex}

An element of $F^{\op {ind}}(X_{\infty}) = \varinjlim _{\la \in \La} F(X_{\la})$ is called and {\em indconstructible function \/} and up to now we have not discussed the role of functions, even though it is called ``function". In fact, the indconstructible function can be considered in a natural way as a function on the proalgebraic variety simply as follows: for $[\alp_{\la}] \in F^{\op {ind}}(X_{\infty}) = \varinjlim _{\la \in \La} F(X_{\la})$ the value of $[\alp_{\la}]$ at a point $(x_{\mu}) \in X_{\infty} = \varprojlim _{\la \in \La} X_{\la}$ is defined by
$$[\alp_{\la}] \Bigl ((x_{\mu}) \Bigr ) := \alp_{\la}(x_{\la})$$
which is well-defined. So, if we let $Fun(X_{\infty}, \bZ)$ be the abelian group of $\bZ$-valued functions on $X_{\infty}$, then the homomorphism
$$\Psi : \varinjlim _{\la \in \La} F(X_{\la}) \to Fun(X_{\infty}, \bZ) \quad \text {defined by} \quad
\Psi \left ([\alp_{\la}] \right ) \left ((x_{\mu}) \right ) := \alp_{\la}(x_{\la})$$
shall be called the ``functionization" homomorphism. 

One can describe this in a fancier way as follows.  Let $\pi_{\la}: X_{\infty} \to X_{\la}$ denote the canonical projection. Consider the following commutative diagram (which follows from $\pi_{\la} = \pi_{\la \mu} \circ \pi_{\mu} (\la < \mu)$):

$$\xymatrix{
F(X_\lambda) \ar[dd]_{\pi^*_{\la \mu}} \ar[dr] ^{\pi^*_\la} \\
& Fun(X_\infty,\bZ) \\
F(X_\mu) \ar[ur]_{\pi^*_\mu}
}$$

Then the ``functionization" homomorphism $\Psi : \varinjlim _{\la \in \La} F(X_{\la}) \to Fun(X_{\infty}, \bZ)$ is the unique homomorphism such that the following diagram commutes:

$$\xymatrix{
& F(X_{\la}) \ar [dl]_{ \rho^{\la}} \ar [dr]^{\pi_{\la}^*} \\
{F^{\op {ind}}(X_{\infty}) } \ar [rr] _{\Psi}& &  {Fun(X_{\infty}, \bZ)}.}
$$
To avoid some possible confusion, the image $\Psi \bigl ([\alp_{\la}] \bigr ) = \pi_{\la}^*\alp_{\la}$ shall be denoted by $[\alp_{\la}]_{\infty}$.  For a constructible set $W_{\la} \in X_{\la}$, by the definition we have
$$[\jeden_{W_{\la}}]_{\infty} = \jeden_{\pi_{\la}^{-1}(W_{\la})}.$$
${\pi_{\la}}^{-1}(W_{\la})$ is called a {\em proconstructible or a cylinder set\/}, mimicking [Cr]. And the characteristic function supported on a proconstructible set is called a {\em procharacteristic function\/} and a finite linear combination of procharacteristic functions is called a {\em proconstructible function\/}. Let $F^{\op {pro}}(X_{\infty})$ denote the abelian group of all proconstructible functions on the proalgebraic variety $X_{\infty} = \varprojlim_{\la \in \La} \Bigl \{ X_{\la}, \pi_{\la \mu}: X_{\mu} \to X_{\la} \Bigr \}$. Thus we have the following 

\begin{pro} For a proalgebraic variety $X_{\infty} = \varprojlim_{\la \in \La} \Bigl \{ X_{\la}, \pi_{\la \mu}: X_{\mu} \to X_{\la} \Bigr \}$ 
$$F^{\op {pro}}(X_{\infty}) = \op {Image}\left (\Psi : F^{\op {ind}}(X_{\infty}) \to Fun (X_{\infty}, \bZ) \right ) = \bigcup_{\mu} \pi_{\mu}^* \bigl  (F(X_{\mu}) \bigr ).$$
If the structure morphisms $\pi_{\la \mu}: X_{\mu} \to X_{\la}$ ($\la < \mu$) are all surjective, then we have
$$F^{\op {ind}}(X_{\infty}) \cong F^{\op {pro}}(X_{\infty}).$$
\end{pro}

In the case of the arc space $\Cal L(X)$ of a nonsingular variety $X$, since each structure morphism $\pi_{n,n+1}: \Cal L_{n+1}(X) \to \Cal L_n (X)$ is always surjective, we get the following 

\begin{cor} \label{naive-motivic}
Assume $X$ is a {\em nonsingular\/} variety of dimension $d$. Then we have
for the arc space $\Cal L(X)$ the canonical isomorphism
$$F^{\op {ind}}\bigl (\Cal L(X) \bigr )\cong F^{\op {pro}}\bigl (\Cal L(X) \bigr ),$$
together with the following canonical Grothendieck class homomorphisms
$$\Ga_X^{\op {ind}}: F^{\op {pro}}(\Cal L(X)) \to K_0(\Cal V/X)_{[\bL^d]}
\quad \text{and} \quad
\Ga^{\op {ind}}: F^{\op {pro}}(\Cal L(X)) \to K_0(\Cal V)_{[\bL^d]}$$
described by
$$\Ga_X^{\op {ind}}\left ([\alp_n]_{\infty} \right ) = \frac {\pi_{0,n*}(\Ga_{\Cal L_n(X)} (\alp_n))}{[\bL]^{nd}}
\quad \text{and} \quad
\Ga^{\op {ind}}\left ([\alp_n]_{\infty} \right ) = \frac {\Ga (\alp_n)}{[\bL]^{nd}}. $$
In particular, we get that 
$\Ga_X^{\op {ind}} \left (\jeden_{\Cal L(X)} \right ) = [id_X]$ and
$\Ga^{\op {ind}} \left (\jeden_{\Cal L(X)} \right ) = [X].$
\end{cor}

So $\Ga_X^{\op {ind}}$ and $\Ga^{\op {ind}}$ define finitely additive measures $\mu_X$ and $\mu$
on the algebra of cylinder sets in the arc space $\Cal L(X)$ of a {\em nonsingular\/} variety $X$,
which are called {\em naive motivic measures\/}. So we can rewrite $\Ga_X^{\op {ind}}(\alp)$ and $\Ga^{\op {ind}}(\alp)$
for $\alp \in F^{\op {pro}}\bigl (\Cal L(X) \bigr )$ as motivic integrals
$$\Ga_X^{\op {ind}}(\alp)=\int_{\Cal L(X)} \: \alp \:d\mu_{X}
\quad \text{and} \quad
\Ga^{\op {ind}}(\alp)=\int_{\Cal L(X)} \: \alp\: d\mu \:.$$
Therefore we see that our proalgebraic Grothendieck class homomorphisms of Theorem \ref{pro-Gro} are a generalization of these naive motivic measures.
Here for ``naive" we point out that for
the applications of a good motivic integration theory (e.g., as described in the next section) 
one needs to take values in a suitable {\em completion\/} of $K_0(\Cal V/X)_{[\bL^d]}$ or
$K_0(\Cal V)_{[\bL^d]}$ so that more general sets than just cylinder sets become
``measurable". Also the use of the ``relative measure" $\Ga_X^{\op {ind}}$ over the
base space $X$ due to Looijenga \cite{Looijenga} is more recent, and will become
important in the next section. \\

When we extend the {\em Chern--Schwartz--MacPherson class transformation\/} [Mac1] 
to a category of proalgebraic varieties, we appeal to the {\em Bivariant Theory\/}.
To fit it in with the notion of {\em bifunctors\/} used before, we assume for simplicity that {\em all\/} morphisms in the underlying category
are ``proper", e.g. in the topological context we work only with {\em compact\/} spaces.
More generally, applying {\em bivariant characteristic classes, namely Grothendieck transformations\/} (as in Theorem \ref{uniqueness of GT}), given in the previous section, we can get a general theory of {\em characteristic classes of proalgebraic varieties\/} as follows:

For a morphism $f:X \to Y$ and a bivariant class $b \in \bB(X \xrightarrow {f} Y)$,
the pair $(f;b)$ is called a {\it bivariant-class-equipped morphism} and we just express $(f;b): X \to Y$. Let $\bB$ be a bivariant theory having units. 
If a system $\bigl \{b_{\la \mu} \bigr \}$ of bivariant classes satisfies that
$$b_{\la \la} = 1_{X_{\la}} \quad \text {and} \quad b_{\mu \nu} \bullet b_{\la \mu} = b_{\la \nu} \quad (\la < \mu < \nu), $$
then we call the system {\it a projective system of bivariant classes}.
If $\bigl \{\pi_{\la \mu}: X_{\mu} \to X_{\la} \bigr \}$ and $\bigl \{b_{\la \mu} \bigr \}$ are projective systems, then the system $\bigl \{(\pi_{\la \mu}; b_{\la \mu}): X_{\mu} \to X_{\la} \bigr \}$ shall be called {\it a projective system of bivariant-class-equipped morphisms}.

For a bivariant theroy $\bB$ having units on the category $\Cal C$ and for a projective system $\bigl \{(\pi_{\la \mu}; b_{\la \mu}): X_{\mu} \to X_{\la} \bigr \}$ of bivariant-class-equipped morphisms, the inductive limit
$$\varinjlim_{\La} \Bigl \{\bB_*(X_{\la}),  b_{\la \mu} \bullet: \bB_*(X_{\la}) \to \bB_*(X_{\mu})\Bigr \}$$
shall be denoted by
$$\bB_*^{\op {ind}}\Bigl (X_{\infty}; \{b_{\la \mu} \} \Bigr )$$
emphasizing the projective system $\{b_{\la \mu} \}$ of bivariant classes, because the above inductive limit surely depends on the choice of it. So we make the covariant functor $\bB_*$ into a bifunctor using the functorial ``Gysin homomorphisms"
$b_{\la \mu} \bullet: \bB_*(X_{\la}) \to \bB_*(X_{\mu})$ induced by the projective system $\bigl \{b_{\la \mu} \bigr \}$.
For example, in the above Example \ref{typical example} we have that 
$$F^{\op {ind}}(X_{\infty}) = \bF_*^{\op {ind}}\Bigl (X_{\infty}; \bigl \{\jeden_{\pi_{\la \mu}} \bigr \} \Bigr ).$$

\begin{defn} Let $\{f_{\la}: X_{\la} \to Y_{\la}\}_{\la \in \La}$ be a pro-morphism of pro-algebraic varieties $\Bigl \{ X_{\la}, \pi_{\la \mu}: X_{\mu} \to X_{\la} \Bigr \}$ and $\Bigl \{ Y_{\la}, \rho_{\la \mu}: Y_{\mu} \to Y_{\la} \Bigr \}$. If the following commutative diagram for $\la <\mu$

$$\CD
X_{\mu}@> {f_{\mu}} >> Y_{\mu}\\
@V {\pi_{\la \mu}}VV @VV {\rho_{\la \mu}}V\\
X_{\la}@>> {f_{\la}} > Y_{\la} \endCD
$$
is a fiber square, then we call the pro-morphism $\{f_{\la}:  X_{\la} \to Y_{\la} \}_{\la \in \La}$ a {\it fiber-square pro-morphism}, abusing words.
\end{defn}

With these definitions we have the following theorem:
\begin{thm} 
(i) Let $\ga: \bB \to \bB'$ be a Grothendieck transformation between two bivariant theories $\bB, \bB': \Cal C \to \Cal A$ and let 
$\bigl \{(\pi_{\la \mu}; b_{\la \mu}): X_{\mu} \to X_{\la} \bigr \}$ be a projective system of bivariant-class-equipped morphisms. Then we get the following pro-version of the natural transformation $\ga_*:\bB_* \to \bB'_*$:
$$\ga_*^{\op {ind}}: \bB_*^{\op {ind}}\Bigl (X_{\infty}; \{b_{\la \mu} \} \Bigr ) \to 
{\bB'_*}^{\op {ind}}\Bigl (X_{\infty}; \{\ga (b_{\la \mu}) \} \Bigr ).$$
(ii)  Let $\{f_{\la}:  Y_{\la} \to X_{\la} \}$ be a fiber-square pro-morphism between two projective systems 
$\bigl \{(\rho_{\la \mu}; d_{\la \mu}): Y_{\mu} \to Y_{\la} \bigr \}$ and 
$\bigl \{(\pi_{\la \mu}; b_{\la \mu}): X_{\mu} \to X_{\la} \bigr \}$ of bivariant-class-equipped morphisms such that 
$d_{\la \mu}= f_{\la}^{\bigstar}b_{\la \mu}$. Then we have the following commutative diagram:
$$\CD \bB_*^{\op {ind}}(Y_{\infty}; \{d_{\la \mu} \}) @> {\ga_*^{\op {ind}}} >> 
{\bB'}_*^{\op {ind}}(Y_{\infty}; \{\ga (d_{\la \mu}) \})\\
@V {{f_{\infty}}_*}VV @VV {{f_{\infty}}_*}V\\\bB_*^{\op {ind}}(X_{\infty}; \{b_{\la \mu} \}) @>> {\ga_*^{\op {ind}}} > {\bB'}_*^{\op {ind}}(X_{\infty}; \{\ga (b_{\la \mu}) \}). \endCD$$
(iii) Let $\bB_*(pt) = \bB'_*(pt)$ be a commutative ring $\Cal R$ with a unit and we assume that the homomorphism 
$\ga: \bB_*(pt) \to \bB'_*(pt)$ is the identity. Let $P = \{ p_{\la \mu} \}$ be a projective system of elements 
$p_{\la \mu} \in \Cal R$. Then we get the commutative diagram
$$\xymatrix{\bB_{*, P}^{\op {st.ind}}
\Bigl (X_{\infty}; \{b_{\la \mu} \} \Bigr )  \ar[dr]_ {\chi^{\op {ind}}_{\bB_*}} \ar[rr]^ {\ga_*^{\op {ind}}} && 
\bB'^{\op {st.ind}}_{*, P}\Bigl (X_{\infty}; \{\ga (b_{\la \mu}) \} \Bigr )  \ar[dl]^{\chi^{\op {ind}}_{\bB'_*}}\\
& \varinjlim _{\la \in \La} \Bigl \{ \times p_{\la \mu} : \Cal R \to \Cal R \Bigr \}.}$$
\end{thm}

If we apply this theorem to the Brasselet's bivariant Chern class \cite{Brasselet1} or to the one of \cite{BSY1}, we get a proalgebraic version $c_*^{\op {ind}}$ of the Chern--Schwartz--MacPherson class transformation $c_*: F \to H_*$.
But of course we also can apply it to the bivariant versions of our motivic characteristic class transformations
$mC_*$ and $T_{y*}$. 

As a very simple example, consider a proalgebraic variety
$X_{\infty} = 
\varprojlim_{\la \in \La} \Bigl \{ X_{\la}, \pi_{\la \mu}: X_{\mu} \to X_{\la} \Bigr \}$, whose structure maps $\pi_{\la \mu}$ are smooth (and therefore ``Euler morphisms") and proper.
Then we can apply the proalgebraic Chern--Schwartz--MacPherson class transformation 
$c_*^{\op {ind}}$ to 
$$F^{\op {ind}}(X_{\infty}) = \bF_*^{\op {ind}}\Bigl (X_{\infty}; \bigl \{\jeden_{\pi_{\la \mu}} \bigr \} \Bigr ).$$
Note that in this case $\gamma(\jeden_{\pi_{\la \mu}})=c^{*}(T_{\pi_{\la \mu}})\bullet
[\pi_{\la \mu}]$ by the {\em Verdier Riemann-Roch theorem for a smooth morphism\/}, so that ${\bH_*}^{\op {ind}}\Bigl (X_{\infty}; \{\ga (\jeden_{\la \mu}) \} \Bigr )$
is just the inductive limit of the following system of ``twisted" smooth pullbacks in
homology:
$$\pi_{\la \mu}^{!!}:=
c^{*}(T_{\pi_{\la \mu}})\cap \pi_{\la \mu}^!: H_*(X_{\la};\bZ)\to H_*(X_{\mu};\bZ)\:.$$
Suitable modifications of such ``inductive limits of twisted smooth pullback morphisms"
are closely related to the construction of {\em equivariant characteristic classes\/}
(compare for example with \cite[\S 3.3, p.12-13]{Oh}).

\section{Stringy and arc characteristic classes of singular spaces}
In this last section we explain another and more recent extension of characteristic classes to singular spaces.
These are {\em not\/} functorial theories as before, but have a better ``birational invariance", in particular
for {\em K-equivalent\/} manifolds, i.e. $M_i$ ($i=1,2$) are irreducible (or pure dimensional) complex algebraic manifolds
dominated by a third such manifold $M$, with $\pi_i: M\to M_i$ proper birational ($i=1,2$) such that the pullbacks of their
canonical bundles (or divisors) $\pi_1^*K_{M_1}\simeq \pi_2^*K_{M_2}$ are isomorphic (or linearly equivalent).
For example $M_1$ and $M_2$ are both {\em Calabi-Yau manifolds\/} in the sense that their canonical bundle is trivial.
In fact the origin of these classes and invariants goes back to two different generalizations
of Hirzebruch's $\chi_y$-genus (which was related to our motivic characteristic classes $mC_*$ and $T_{y*}$).\\

The first one is the {\em E-polynomial or Hodge characteristic\/} $E(X)(u,v)\in \bZ[u,v]$ defined in terms of Deligne's {\em mixed Hodge structure\/} \cite{De-hodgeI, De-hodgeII} for the cohomology with compact support $H^*_c(X,\bQ)$ of a complex algebraic variety. We have that $E(X)(1,1)=\chi(X)$ for any variety $X$ and $E(X)(-y,1)=\chi_y(X)$ for $X$ smooth and compact.
In the 90's {\em V. Batyrev\/} \cite{Batyrev-stringy} extended this E-polynomial to a {\em stringy E-function $E_{str}$ and stringy Euler 
numbers $\chi^{str}$\/} of ``log-terminal pairs" $(X,D)$ relating them in some cases known as the ``McKay correspondence" to orbifold invariants of suitable
quotient varieties. He also used in \cite{Batyrev-birat} methods from p-adic integration theory to prove that different ``crepant resolutions" of a given singular space, and also {\em birationally equivalent Calabi-Yau manifolds\/}, have equal Betti numbers.
Later on  {\em  M. Kontsevich\/} \cite{Kontsevich} invented ``motivic integration" (with some analogy to  p-adic integration)
for extending these results from Betti numbers to Hodge numbers.\\

The other generalization of the $\chi_y$-genus is the (complex) {\em elliptic genus $ell_k$\/} studied by I. Krichever
\cite{Krich} and  G. H\"{o}hn \cite{Hoehn}. As observed by {\em Totaro\/} \cite{Totaro}
(and compare with \cite{BF}), this is the most general genus on the complex cobordism ring $\Omega^{U}_*\otimes \bQ$, which can be invariant under a suitable notion of ``flips". Later on this was extended by 
{\em L. Borisov and A. Libgober\/}
\cite{BL-elliptic} and {\em C.-L. Wang\/} \cite{Wang} for showing the invariance of this  elliptic genus $ell_k$ for $K-equivalent$
complex algebraic manifolds, a notion coming from ``minimal model theory". Both works use the very deep ``weak factorization theorem" (\cite{AKMW} and \cite{Wlodarczyk}) for the comparison of different resolution spaces. They also introduced in this way
the {\em elliptic homology class\/} $\Cal E ll_*(X)$ of a {\em $\bQ$-Gorenstein log-terminal\/} singular complex algebraic variety
$X$ \cite{BL-McKay, Wang}. Here {\em $\bQ$-Gorenstein\/} for a normal irreducible (or pure dimensional) variety $X$ just means
that some multiple $r\cdot K_X$ ($r\in \bN$) of the canonical Weil divisor $K_X$ is already a Cartier divisor, with $r=1$ corresponding
to a {\em Gorenstein variety\/} (e.g. $X$ is smooth). Here $K_X$ is just the closure of a canonical divisor on the regular part.
In fact, Borisov-Libgober proved in \cite{BL-McKay} for this elliptic homology class a very general
version of the ``McKay correspondence".\\

More recently simpler {\em stringy Chern classes\/} $c^{str}_*(X)$ were introduced by {\em Aluffi\/} \cite{Aluffi-modification}, based 
on the ``weak factorization theorem", and independently by {\em de Fernex, Lupercio, Nevins and Uribe\/} \cite{FLNU}, based on ``motivic integration" and MacPherson's
functorial Chern class transformation $c_*$. In fact Aluffi pointed out that there are two possible notions of such classes,
depending on two different choices of a system of ``relative canonical divisors"
$K_{\pi}$ for suitable resolution of singularities $\pi: M\to X$
(i.e. $\pi$ is proper and $M$ smooth), which he calls the ``$\Omega$-flavor"  and ``$\omega$-flavor".\\ 
The ``$\omega$-flavor" is related to ``stringy invariants and characteristic classes" (like $E_{str}, \Cal E ll_*$
and $c^{str}_*$). Here one assumes $X$ is irreducible and $\bQ$-Gorenstein so that the relative canonical divisor
$K_{\pi}:= K_M- \pi^*K_X$ is at least a $\bQ$- Cartier divisor (class). 
Moreover it is supported on the exceptional locus $E$
of the resolution, which is supposed to be (contained in) a normal crossing divisor with smooth irreducible components $E_i$.
Then $K_{\pi}\simeq \sum_i \:a_i\cdot E_i$ for some fixed $a_i\in \bQ$ (depending on the resolution).
And for the definition of all these ``stringy invariants"  one needs the condition $a_i>-1$ for all $i$,
which exactly means that $X$ has only {\em log-terminal singularities\/}. If this condition holds for one such resolution,
then it is true for any resolutions of this type.
A resolution $\pi$ is called {\em crepant\/}, if $K_M\simeq \pi^*K_X$, e.g. all $a_i=0$ for $E$ a normal crossing divisor
as before.\\
The ``$\Omega$-flavor" is related to what we call ``arc invariants and arc characteristic classes",
because these generalize corresponding ``arc invariants" of {\em Denef and Loeser\/} \cite[\S 6]{Denef-Loeser1} and
\cite[\S 4.4.1]{Denef-Loeser2}, which they introduced already before by their work on ``motivic integration".
In this case $X$ is only assumed to be pure $d$-dimensional and $K_{\pi}$ is defined for all resolutions $\pi$
such that the canonical map $\pi^*\Omega^d_X\to \Omega^d_M$ of K\"{a}hler differentials has an image $\Cal I\otimes \Omega^d_M$
with $\Cal I$ a principal ideal in $\Cal O_M$ (this can always be achieved by Hironaka \cite{Hironaka}). Then $K_{\pi}$ is defined by
$\Cal I=\Cal O_M(-K_{\pi})$. The effective Cartier divisor $K_{\pi}$ is again supported on the exceptional locus $E$
of the resolution, which can also be supposed to be (contained in) a normal crossing divisor with smooth irreducible components $E_i$.
Then one can introduce the $a_i\in \bN_{0}$ as before.\\

For $X$ already smooth, both notations of a relative canonical divisor $K_{\pi}$ agree with the divisor of the Jacobian of $\pi$
defined by the section $s$ of $K_M\otimes \pi^*K_X^*$ corresponding to the canonical map $\pi^*\Omega^d_X\to \Omega^d_M$.
Note that in both cases the corresponding resolutions $\pi: M\to X$ as above form a directed set, i.e. two of them
can be dominated by a third one of this type (and taking suitable limits over this directed set corresponds to the view point of
Aluffi \cite{Aluffi-modification}).
If $\pi':M'\to M$ is a proper birational map with $\pi$ and $\pi\circ \pi'$
as above, then the relative canonical divisors have (in both cases) the following crucial transitivity property:
\begin{equation} \label{reldivisor}
K_{\pi\circ \pi'}\simeq K_{\pi'} + \pi'^*K_{\pi}\:.
\end{equation}
Then  all these new invariants $I(X)$ for a singular space $X$ as above are defined as
$$I(X):=\pi_*(I(M)\cdot J(\{E_i,a_i\})) \in \bB_*(X)$$
for such a special resolution $\pi: M\to X$, with $E$ a normal crossing divisor with smooth irreducible components $E_i$, 
where $I(M)\in \bB_*(M)$ is the corresponding invariant of the
smooth space  $M$, together with some ``correction term" $J(\{E_i,a_i\})\in \bB^*(M)$ depending on the exceptional divisor $E$ and the multiplicities
$a_i$ defined by the relative canonical divisor $K_{\pi}$. 
Here $\bB_*$ and $\bB^*$ are suitable covariant and contravariant theories taking values in the category of Abelian groups 
and commutative rings with unit, related by the projection formula as in Example \ref{typical example2}. Typical examples are
\begin{enumerate}
\item $\bB_*(X)=\bB^*(X)=\Lambda$ is a commutative ring with unit (with all pullbacks and push downs the identity transformation $id_R$)
so that $I(M)\in R$ corresponds to a suitable generalized ``Euler characteristic type invariant".
\item  $\bB_*$ and $\bB^*$ correspond to a suitable (co)homology theory like $(\bB_*(X),\bB^*(X))=
(H_*^{BM}(X)\otimes \Lambda, H^*(X)\otimes \Lambda)$ or $(\bB_*(X),\bB^*(X))=(\bG_0(X)\otimes \Lambda,\bold K^0(X)\otimes \Lambda)$
so that $I(M)\in B_*(M)$ is a suitable characteristic class of $M$.
\item $\bB(X):=\bB_*(X)=\bB^*(X)$ is a bifunctor as in Example \ref{typical example2}, e.g. like constructible
functions $\bB(X)=F(X)\otimes \Lambda$ or relative Grothendieck rings of varieties  $K_0( \Cal V /X)\otimes \Lambda$
coming up from ``motivic integrals".
\end{enumerate}
If $I(X)\in \bB_*(X)$ is such an invariant not depending on the choice of the resolution $\pi$,
then the same is true for $\gamma_*(I(X))\in \bB'(X)$ for any natural transformation of
covariant theories $\gamma_*: \bB_*\to \bB'_*$. For example $I(X)\in H_*(X)\otimes \Lambda$
is a characteristic homology class with $X$ compact, and $deg:=\gamma_{*}:  H_*(X)\otimes \Lambda
\to H_*(\{pt\})\otimes \Lambda=\Lambda$ is just its degree (or push down to a point).
Or we apply suitable ``completions" of our motivic characteristic class transformations
$mC_*$ and $T_{y*}$ to invariants $I(X)$ coming from motivic integration!\\

For showing that the final result $I(X)$ does not depend on the choice of the resolution, either ``motivic integration with its
transformation rule" related to the ``Jacobian factor" $J(\{E_i,a_i\})$ is used:
\begin{equation} \label{transformation}
\int_{\Cal L(M)} \: \bL^{-\alpha} \:d\tilde{\mu}_{M} = \pi'_* \int_{\Cal L(M')} \: \bL^{-(\pi'^*\alpha+K_{\pi'})} \:d\tilde{\mu}_{M'}
\end{equation}
for $\pi': M'\to M$ a proper birational map of manifolds and $\bL:=[\bC] \in K_0( \Cal V)$. This suggests to think of
$I(X)$ as the push down of an ``integral with respect to the invariant $I(M)$":
$$I(X)= \pi_* \int_{M} \: \bL^{-K_{\pi}} \:dI(M)\:. $$
Or the "weak factorization theorem" is used,
in which case only the invariance under suitable ``blowing ups" has to be checked.\\

Moreover $J(\{E_i,a_i\})=1$ in case all $a_i=0$, so that $I(X)=\pi_*(I(M))$ in case of a {\em crepant resolution\/}.
In particular $\pi_*(I(M))$ does {\em not\/} depend on the choice of this crepant resolution.
Suppose two maybe {\em singular\/} spaces $X_i$ ($i=1,2$) are K-equivalent in the sense that they are
dominated by a manifold $M$, with $\pi_i: M\to X_i$   a resolution of singularities such that
the relative canonical divisors $K_{\pi_i}$ are defined ($i=1,2$) and equal. After taking another resolution of $M$,
we can even assume that the exceptional locus of both maps is contained in a normal crossing divisor $E$ with smooth
irreducible components $E_i$ (here we use the transitivity property of the relative canonical divisors).
But then the correction factor $J(\{E_i,a_i\})$ for both maps is the same, so that
$$I(X_1)=\pi_{1*}(I(M)\cdot J(\{E_i,a_i\})) 
\quad \text{and} \quad 
I(X_2)=\pi_{2*}(I(M)\cdot J(\{E_i,a_i\}))\:,$$
i.e. both invariants $I(X_1)$ and $I(X_2)$ are ``dominated" by the same element coming from $M$.
In particular 
$$I(X_1)=I(X_2)$$ 
in case of ``Euler characteristic type invariants", and 
$$deg(I(X_1))= deg(I(X_2))$$ 
in case of ``characteristic homology classes" for compact spaces $X_i$.
If we are working in the ``$\omega$-flavor" of stringy homology classes $I(X_i)\in H_*^{BM}(X_i)\otimes \Lambda$
for $\bQ$-Gorenstein varieties $X_i$, we can use the first Chern class
$c^1(K_{X_i}):=c^1(r\cdot K_{X_i})/r \in H^2(X_i; \bQ)$ (for $\Lambda$ a $\bQ$-algebra) to modify $I(X_i)$ into
$$ I'(X_i):=f(c^1(K_{X_i})) \cdot I(X_i) \in H_*^{BM}(X_i)\otimes \Lambda \:.$$
By the {\em projection formula\/} also these new invariants $I'(X_1)$ and $I'(X_2)$ are ``dominated" by the same element coming from $M$,
where $f\in \Lambda[[z]]$ can be any power series. If $X_i$ are both Gorenstein, we can do the same thing for corresponding
invariants $I(X_i) \in \bG_0(X)\otimes \Lambda$ by using polynomials in the (inverse) classes 
$[K_{X_i}^{\pm 1}]\in \bold K^0(X_i)$ of the canonical Cartier
divisors (instead of their first Chern classes).\\
 
Note that the approach by resolution of singularities is different from our
approach to functorial ``motivic characteristic classes" based on ``additivity" (i.e. decomposing a singular space into smooth
pieces), but nevertheless they nicely fit together as we now explain.\\

\subsection{Elliptic classes} Let us start with the definition of the {\em (complex) elliptic class\/} $\Cal E LL(E)$ of a complex vector bundle $E\to X$.
Consider the formal power series
$$\Lambda_{t}(E):= \sum_{n\geq 0}\: t^{n}\Lambda^{n}E 
\quad \text{and} \quad
S_{t}(E):= \sum_{n\geq 0}\: t^{n}S^{n}E \:,$$
with $\Lambda^{n}E$ and $S^{n}E$ the corresponding exterior and symmetric power of $E$
(so $\Lambda^{n}E=0$ for $n>\op {rank} \;E$, with $\Lambda_{t}$ the total Lambda class coming up in our definition
of the motivic Chern class transformation $mC_*$ in Corollary \ref{mC*}). Then one has
$$\Lambda_{t}(E\oplus F) = \Lambda_{t}(E)\Lambda_{t}(F), \quad
S_{t}(E\oplus F) = S_{t}(E)S_{t}(F), \quad \text{and} \quad
\Lambda_{t}(E)S_{-t}(E)=1 \:.$$
So these operations extend to the Grothendieck group of complex vector bundles (and similarly in the algebraic context):
$$\Lambda_{t},S_t:  (\bold K(X),\oplus)\to (1+\bold K(X)[[t]],\otimes) \subset (\bold K(X)[[t]],\otimes)\:.$$

Then we define the {\em complex elliptic class\/} 
$$\Cal E LL(E)=\Cal E LL(y,q)(E)\in \bold K(X)[[q]][y^{\pm 1}]$$ 
of a complex vector bundle $E\to X$ as $\Cal E LL(E):= \Lambda_{y}(E^*)\otimes \Cal W(E)$, with
\begin{equation} \label{Ell-c}
\Cal W(E):= \bigotimes_{n\geq 1}\:
\Bigl(\Lambda_{yq^n}(E^*)\otimes \Lambda_{y^{-1}q^n}(E)\otimes S_{q^n}(E^*)\otimes S_{q^n}(E)\Bigr)\:.
\end{equation}
More generally the {\em elliptic class of order $k$\/} 
$$\Cal E LL_k(E)=\Cal E LL_k(y,q)(E)\in \bold K(X)[[q]][y^{\pm 1}] \quad \text{with $k\in \bZ$}$$ 
of a complex vector bundle $E\to X$  is defined as the twisted class
\begin{equation} \label{Ell}
\Cal E LL_k(E):= det(E)^{\otimes -k}\otimes \Cal E LL(E)\:,
\end{equation}
with $det(E):=\Lambda^{\op {rank} \;E}(E)$ being the determinant line bundle of $E$. 
So $\Cal E LL(E)$ (or $\Cal E LL_k(E)$) is a one (or two) parameter deformation
of the total Lambda class $\Lambda_y(E^*)$, with 
$$\Cal E LL_0(E) = \Cal E LL(E)
\quad \text{and} \quad
\Cal E LL(E)|_{q=0}= \Lambda_y(E^*) \:.$$

For $M$ a complex projective  algebraic manifold
(or a compact almost complex manifold) one can introduce as in \S \ref{HRR-GRR} the $\chi=T$-characteristic 
$$\chi(\Cal E LL_k(E))\in \bQ[[k,q]][y^{\pm 1}]$$
of $\Cal E LL_k(E)$ as
\begin{equation*} \begin{split}
\chi(\Cal E LL_k(E))
& :=\int_{M} \:ch^*(\Cal E LL_k(E))\cdot td^*(TM)\cap [M]\\
& =\int_{M} \:e^{-k\cdot c^1(E)}\cdot ch^*(\Cal E LL(E))\cdot td^*(TM)\cap [M] \:.
\end{split} \end{equation*}
Note that in the last term one can introduce $k$ as a formal parameter.
$ch^*(\Cal E LL_k(E))$ and $ch^*(\Cal E LL(E))$)  are {\em multiplicative\/} (but not normalized) characteristic classes so that we get the induced Krichever--H\"{o}hn {\em elliptic genus\/}
$$ell_k: \Omega_*^U\otimes \bQ \to \bQ[[k,q]][y^{\pm 1}] \:,$$
with
\begin{equation}\label{ell}
ell_k(M):= \int_{M} \:e^{-k\cdot c^1(TM)}\cdot ch^*(\Cal E LL(TM))\cdot td^*(TM)\cap [M] \:.
\end{equation}
The corresponding {\em complex elliptic genus\/} $ell:= ell_0: \Omega_*^U\otimes \bQ \to \bQ[[q]][y^{\pm 1}]$ given by
\begin{equation*} \begin{split} 
ell_0(M)
&= \int_{M} \:ch^*(\Cal W(E))\cdot ch^*(\Lambda_y T^*M)\cdot td^*(TM)\cap [M]\\
&=: \chi_{y}(M,\Cal W(E))\\
&= \chi_{y}\Bigl(M, \bigotimes_{n\geq 1}\:
\bigl(\Lambda_{yq^n}(TM^*)\otimes \Lambda_{y^{-1}q^n}(TM)\otimes S_{q^n}(TM^*)\otimes S_{q^n}(TM) \bigr)\Bigr)
\end{split}
\end{equation*}
was formally interpretated by Witten as the $S^1$-equivariant
$\chi_y$-genus $\chi_{y}(S^1,LM)$ of the free loop space $LM=\{f: S^1\to M|f\; \text{smooth}\}$ of $M$ (compare \cite[Appendix III]{HBJ} and \cite{BF}).  
$$\chi_{k,y}(M):= ell_k(M)|_{q=0} \in \bQ[y][[k]]$$
is called the {\em twisted $\chi_y$-genus\/} of $M$:
\begin{equation} \label{twistedchi-y}
\chi_{k,y}(M)= \int_{M} \:e^{-k\cdot c^1(TM)}\cdot ch^*(\Lambda_y(T^*M))\cdot td^*(TM)\cap [M] \:.
\end{equation}
Another specialization is the {\em real elliptic genus\/} $ell|_{y=1}$, which factorizes over the oriented cobordism ring
$$ell|_{y=1}: \Omega_*^{SO}\otimes \bQ \to \bQ[[q]]\:.$$
This one parameter genus interpolates between the signature genus (for $q\to 0$) and the
$\hat{A}$-genus (for $q\to \infty$),
and was formally interpretated by Witten as the $S^1$-equivariant signature $\sigma(S^1,LM)$ of the free loop space $LM$ of 
the oriented manifold $M$ (compare \cite[\S 6]{HBJ} and \cite{BF}).

\begin{rem} \label{norm-ell}
We should point out that there are many different normalizations of the elliptic genus and classes in the literature.
First of all many authors (like \cite{BL-elliptic, BL-McKay, Totaro, Wang}) use $-y$ instead of $y$ so that their
elliptic genus is related to the $\chi_{-y}$-genus. But what is maybe more important, we do not work with
``normalized characteristic classes", i.e. the power series $f(z)\in \bQ[[k,q]][y^{\pm 1}][[z]]$ in the variable $z=c^1$
corresponding to the multiplicative characteristic class $ch^*(\Cal E LL_k(\quad))$ has a constant coefficient $a:=f(0)\neq 1$,
since $ch^*(\Cal E LL(E))|_{q=0}=ch^*(\Lambda_y(E^*))$ implies $a(k=0,q=0)=1+y \in \bQ[y^{\pm 1}]$.
So twisting $f(z)$ to a normalized power series $f(z)/a$ (as used in \cite{BF, Totaro, Wang})
would change the elliptic genus only to
$ell_k(M)/a^n$ for $M$ an (almost) complex manifold of complex dimension $n$,
and similarly a characteristic homology class $cl_i(\quad)\in H^{BM}_{2i}(\quad)\otimes \Lambda$ would just be multiplied by $a^{-i}$.
For example in theorem \ref{T_y_*} we could have started with the natural  transformation
(with respect to proper maps):
$$\tilde{T}_{y*}:= td_*\circ mC_*: K_0(\Cal V/\quad) \to H^{BM}_{2*}(\quad)\otimes \bQ[y] \:,$$
satisfying for $M$ nonsingular the normalization 
$$\tilde{T}_{y*}([M \xrightarrow  {\op {id}} M]) = ch^*(\Lambda_y T^*M)\cdot td^*(TM) \cap [M].$$
And ``twisting" by $1+y$ would then give our motivic characteristic class transformation $T_{y*}$ 
with 
\begin{equation}
T_{y,i}(\quad)= (1+y)^{-i}\cdot \tilde{T}_{y,i}(\quad) \in H^{BM}_{2i}(\quad)\otimes \bQ[y,(1+y)^{-1}] \:.
\end{equation}
But since we work in this section only with
pure dimensional spaces, this ``twisting"  does not matter for the question of getting invariants of pure dimensional
singular complex algebraic varieties. Similarly it will be enough to consider only the complex elliptic genus and classes
corresponding to $k=0$ (as in \cite{BL-elliptic, BL-McKay}),
since the case of general $k$ follows then from the projection formula (as already explained before).
So the elliptic classes $\Cal E ll^*(z,\tau)$ used in \cite{BL-elliptic, BL-McKay} correspond in our notation to
$$\Cal E ll^*(z,\tau)(TM):= y^{-dim(M)/2}\cdot td^*(TM)\cdot ch^*(\Cal E LL(TM))(-y,q)\:,$$ 
with $y=e^{2\pi iz}$ and $q=e^{2\pi i\tau}$.
\end{rem} 

With these notations, we can now explain the definition of Libgober and Borisov 
(\cite[Definition 3.2]{BL-McKay} with $G:=\{id\}$)
for their {\em elliptic class\/}
$\Cal E ll_*((X,D))$ of a ``Kawamata log terminal pair $(X,D)$", i.e. $X$ is a normal irreducible 
complex algebraic variety, with $D$ a $\bQ$-Weil divisor on $X$ such that $K_X+D$ is a $\bQ$-Cartier divisor 
satisfying the following condition: There is a resolution of  singularities $\pi: M\to X$ with the exceptional locus $E$
and the support of $K_{\pi}(D):=K_M-\pi^*( K_X+D)$ contained in a normal crossing divisor with smooth irreducible
components $E_i$ ($i\in I$) such that $K_{\pi}(D)\simeq \sum_i\: a_i\cdot E_i$, with all $a_i\in \bQ$ satisfying the inequality
$a_i>-1$. Note that the last condition is then independent of the choice of such a resolution
(compare \cite[Definition 2.34, Corollary 2.31]{Kollar-Mori}), with the case $D=0$
corresponding to the case ``$X$ is $\bQ$-Gorenstein with only log-terminal singularities".
Moreover, the ``relative canonical divisor $K_{\pi}(D)$ of $D$" also satisfies the transitivity property
\begin{equation} \label{reldivisor-D}
K_{\pi\circ \pi'}(D)\simeq K_{\pi'} + \pi'^*K_{\pi}(D)
\end{equation}
for $\pi':M'\to M$ a proper birational map with $\pi$ and $\pi\circ \pi'$ as before. Then
\begin{equation} \label{ell-BL}
\Cal E ll_*((X,D)):=\pi_*\bigl (\Cal E ll^*(TM)\cap [M]) \cap \prod_i \: J(E_i,a_i) \bigr) \:,
\end{equation}
with 
$$J(E_i,a_i)(z,\tau):= \frac{\theta(\frac{e_i}{2\pi i}-(a_i+1)z,\tau)\theta(-z,\tau)}
{\theta(\frac{e_i}{2\pi i}-z,\tau)\theta(-(a_i+1)z,\tau)} \in H^*(M;\bQ)[[y,q]]\:.$$
Here $\theta(z,\tau)$ is the Jacobi theta function in $y=e^{2\pi iz}$ and $q=e^{2\pi i\tau}$,
with $e_i=c^1(E_i)\in H^2(M,\bZ)$ the first Chern class of the smooth divisor $E_i$.\\

The proof of the independence of the resolution $\pi$ uses the ``weak factorization theorem"
for reducing it to the comparison with a suitable blowing up along a smooth center. 
Using some modularity properties of the $\theta$-function, this is finally reduced to the
vanishing of a suitable residue (of an elliptic function with exactly one pole,
compare \cite[p.11]{BL-McKay} and \cite[\S 4]{Wang}). If $X$ is compact, then
\begin{equation}
ell((X,D)):=deg\bigl( \Cal E ll_*((X,D)) \bigr)
\end{equation}
is just the {\em singular eliptic genus\/} of the Kawamata log terminal pair $(X,D)$
as defined in \cite[Definition 3.1]{BL-elliptic} (up to a normalization factor).\\

Later on we only need the following limit formula (with $y=e^{2\pi iz}$):
\begin{equation} \label{limit-theta}
\begin{split} 
\lim_{\tau\to i\infty} \: J(E_i,a_i)(z,\tau) &= \frac{(y-1)(1-y^{a_i+1}e^{-e_i})}{(y^{a_i+1}-1)(1-ye^{-e_i})} \\
&= 1+ \frac{(y-y^{a_i+1})(1-e^{-e_i})}{(y^{a_i+1}-1)(1-ye^{-e_i})} \:.
\end{split}
\end{equation}
Note that the multiplicative characteristic class 
$$\tilde{T}_y^*(E):= ch^*(\Cal E LL(E))|_{q=0}\cdot td^*(E)=
ch^*(\Lambda_y(E^*))\cdot td^*(E)$$
exactly corresponds to the power series $f(z)=\frac{z(1+ye^{-z})}{1-e^{-z}}$
in the variable $z=c^1$ (compare with section \ref{gHRR}).
If we denote for $J\subset I$ the closed embedding $i_J: E_J:= \bigcap_{i\in J}\;E_i \to M$
of the submanifold $E_J$ (with $E_\emptyset:=M$), then one has by the ``adjunction formula"
$i_{J*}i^*_J= \prod_{i\in J}\:e_i\;\cap\;$, with $TE_J= i_J^*(TM-\sum_{i\in J}\; \Cal O(E_i))$ 
(compare \cite[p.36]{HBJ}):
$$i_{J*}(\tilde{T}_y^*(TE_J)\cap [E_J])=(\tilde{T}_y^*(TM)\cap [M])\cap 
\prod_{i\in J} \:\frac{1-e^{-e_i}}{1+ye^{-e_i}} \:.$$
So altogether we get the following ``limit formula" (with $y=e^{2\pi iz}$):
\begin{equation} \label{limit-ell}
\lim_{\tau\to i\infty}\: y^{dim(X)/2}\cdot\Cal E ll_*((X,D))=
\pi_*\bigl ( \sum_{J\subset I}\: i_{J*}(\tilde{T}_{-y*}(E_J))\cdot 
\prod_{i\in J} \:\frac{y-y^{a_i+1}}{y^{a_i+1}-1} \bigr)\:.
\end{equation}
Recall that we use the notation $cl_*(E_J)=cl^*(TE_J)\cap [E_J]$ for the characteristic
homology class of a manifold (corresponding to a characteristic class $cl^*$ of vector bundles).

\subsection{Motivic integration} 
Motivic integration was invented by Kontsevich \cite{Kontsevich} for showing that birational equivalent Calabi-Yau manifolds
have equal Hodge numbers. In all details with many different applications it was developed by Denef-Loeser
(e.g. \cite{Denef-Loeser1, Denef-Loeser2, Denef-Loeser3}), with some improvements by Looijenga \cite{Looijenga},
who in particular introduced the calculus of relative Grothendieck rings $K_0( \Cal V /X)$ of algebraic varieties.
For a nice introduction to ``stringy invariants of singular spaces" we recommend \cite{Veys-zeta,Veys-stringy}.
Even though motivic integration can be directly studied on singular spaces, we restrict ourselves to the simpler case of
smooth spaces, which will be enough for our applications. Also in this way it can easily be compared to results
coming from the use of the ``weak factorization theorem". For a quick introduction to ``motivic integration on
smooth spaces" compare with \cite{Craw} (where by Corollary \ref{naive-motivic}
all arguments of \cite{Craw} extend to the framework of ``relative motivic measures).\\

Let $M$ be a pure $d$-dimensional complex algebraic manifold and $E=\sum_{i=1}^{k}\:a_{i}E_{i}$
be an effective normal crossing divisor (e.g. $a_{i}\in \bN_{0}$) on $M$, with
smooth irreducible components $E_{i}$. Then one can introduce on the arc space $\Cal L(M)=\{\gamma_u|u\in M\}$ the order function
along $E$:
$$ord(E):=\sum_i\; a_i\cdot ord(E_i): \Cal L(M)\to \bN_0\cup \infty \:,$$
with $ord(E_i)(\gamma_u):=ord_0\: f_i\circ \gamma_u(t)$ the zero order of $f_i\circ \gamma_u(t)\in \bC[[t]]$,
if $f_i$ is a local defining equation of $E_i$ near the point $u\in M$.
In particular
$$ord(D_i)(\gamma_u)=0 \Leftrightarrow u\notin D_i 
\quad \text{and} \quad 
ord(D_i)(\gamma_u)=\infty  \Leftrightarrow \gamma_u\subset D_i \:.$$ 
Then $\{ord(E)=n\}\subset \Cal L(M)$ is for all $n\in \bN_0$ a {\em proconstructible or cylinder set\/}
in the sense of \S \ref{section-pro}. Then one would like to introduce the following motivic integral:
\begin{equation} \label{motint1}
\int_{\Cal L(M)} \: \bL^{-ord(E)}d\mu_{M} :=
\sum_{p\in \bN_0}\: \mu_{M}(\{ord(E)=p\}) \cdot \bL^{-p}
\end{equation}
with values in the localized ring $K_0(\Cal V/M)_{[\bL^d]}$ as in Corollary \ref{naive-motivic}.
Recall that we normalized the (naive) motivic measure $\mu_{M}$ in such a way that we get for $E=0$:
$$\int_{\Cal L(M)} \: 1 \:d\mu_{M}=[M]\in K_0(\Cal V/M)_{[\bL^d]}\:.$$
But the problem with the definition (\ref{motint1}) is that this is not a finite series, and
that $\{ord(E)=\infty\}$ is {\em not\/} a cylinder set in $\Cal L(M)$. Both problems are solved by taking a suitable completion
of $K_0(\Cal V/M)_{[\bL^d]}$. More precisely for $X$ a complex algebraic variety let $\widehat{\bM}(\Cal V/X)$ be the completion  
of $K_0(\Cal V/X)[\bL^{-1}]$ with respect to the following dimension filtration
(for $k\to -\infty$):
$$F_{k}( K_0(\Cal V/X)[\bL^{-1}])\quad \text{is generated by} \quad
[X'\to X]\bL^{-n} \quad \text{with $dim(X')-n\leq k$.}$$

\begin{rem} \label{completion}
Here we consider $K_0(\Cal V/X)$ as an algebra over 
$K_0(\Cal V):=K_0(\Cal V/\{pt\})$ by the pullback $const^*$ for
$const: X\to \{pt\}=Spec(\bC)$ the constant structure map. If $S\subset K_0(\Cal V)$ is a multiplicatively closed subset,
then we can localize the commutative ring  $K_0(\Cal V/X)$ with respect to the induced multiplicatively closed subset
$const^*(S)\subset K_0(\Cal V/X)$, or we can localize $K_0(\Cal V/X)$ as an $K_0(\Cal V)$-module with respect to $S$.
Both localizations can be identified, since $const^*$ is injective (compose with any map $\{pt\}\to X$),
and are denoted by $K_0(\Cal V/X)_{S}$. In case $S=\{\bL^n|n\in \bN_0\}$, with $\bL:=[\bC]\in K_0(\Cal V)$, we also use
the notation $K_0(\Cal V/X)[\bL^{-1}]$ above.

Also note that the filtration and completion as above are compatible with push down $f_*$ and exterior product $\times$
so that in particular $\widehat{\bM}(\Cal V/X)$ is a $\widehat{\bM}(\Cal V):=\widehat{\bM}(\Cal V/\{pt\})$-module,
with an induced $\widehat{\bM}(\Cal V)$-linear
push down $f_*: \widehat{\bM}(\Cal V/X)\to \widehat{\bM}(\Cal V/Y)$ for $f: X\to Y$ an algebraic morphism.
\end{rem}

Let us come back to our motivic integral (\ref{motint1}) on the manifold $M$.
The composed {\em relative motivic measure\/}
$$\tilde{\mu}_{M}: F^{pro}(\Cal L(M))\to \widehat{\bM}(\Cal V/M)$$
can now be extended from cylinder sets to a more general class of ``measureable subsets" of the arc space
$\Cal L(M)$ in such a way that $\{ord(E)=\infty\}$ becomes measureable with measure $0$, and the series
(\ref{motint1}) above converges in $\widehat{\bM}(\Cal V/M)$. So now one can define
\begin{equation} \label{motint2}
\int_{\Cal L(M)} \: \bL^{-ord(E)}d\tilde{\mu}_{M} :=
\sum_{p\in \bN_0}\: \tilde{\mu}_{M}(\{ord(E)=p\}) \cdot \bL^{-p} \: \in \widehat{\bM}(\Cal V/M)\:.
\end{equation}
Moreover it can easily be computed:
\begin{equation} \label{eq:motint} \begin{split}
\int_{\Cal L(M)} \: \bL^{-ord(E)}d\tilde{\mu}_{M} &=
\sum_{I\subset \{1,\dots,k\}} \: [E_{I}^{o}\to M] \cdot \prod_{i\in I}\:
\frac{\bL-1}{\bL^{a_{i}+1}-1} \\
&= \sum_{I\subset \{1,\dots,k\}} \: [E_{I}\to M] \cdot \prod_{i\in I}\:
(\frac{\bL-1}{\bL^{a_{i}+1}-1} -1) \:.
\end{split} \end{equation}
Here we use the notation:
$$E_{I}:=\bigcap _{i\in I}\: E_{i} \quad \text{(with $E_{\emptyset}:=M$), and}
\quad
E_{I}^{o}:=E_{I}\backslash \ \bigcup_{i\in \{1,\dots,k\}\backslash I} \:E_{i} \:,$$
and the factor $(\bL^{a_{i}+1}-1)^{-1}=\bL^{-(a_{i}+1)}\cdot (1-\bL^{-(a_{i}+1)})^{-1}$  has to be developed
as the corresponding geometric series in $\widehat{\bM}(\Cal V)$.
Moreover one gets the last equality in (\ref{eq:motint}) by multiplying out
the following products:
\begin{equation} \label{eq:motint3} \begin{split}
& \prod_{i=1}^{k}\:\Bigl( b_{i}\cdot[E_{i}\to M]+ [M\backslash E_{i}\to M]\Bigr) =\\
& \prod_{i=1}^{k}\: \Bigl( (b_{i}-1)\cdot [E_{i}\to M]+ [id_{M}] \Bigr) \in 
\widehat{\bM}(\Cal V/M)\:,
\end{split} \end{equation}
with $b_{i}:=(\bL-1)(\bL^{a_{i}+1}-1)^{-1} \in \widehat{\bM}(\Cal V)$.
Recall that multiplication in $\widehat{\bM}(\Cal V/M)$ is induced from taking the
fiber product over $M$.\\

The other piece of information that we need is the {\em transformation rule\/}
\begin{equation} \label{transformation2}
\int_{\Cal L(M)} \: \bL^{-ord(E)} \:d\tilde{\mu}_{M} = \pi'_* \int_{\Cal L(M')} \: \bL^{-ord(\pi'^*E+K_{\pi'})} \:d\tilde{\mu}_{M'}
\end{equation}
for $\pi': M'\to M$ a proper birational map of pure dimensional complex algebraic manifolds such that
$\pi'^*E+K_{\pi'}$ is a normal crossing divisor with smooth irreducible components.\\

Assume now that we have a proper birational map $\pi: M\to X$, with $X$ pure dimensional but maybe singular,
together with a Cartier divisor $D$ on $M$ such that $D$ and the exceptional locus of $\pi$ are contained
in (the support of) $E$. Finally we assume 
$$K_{\pi}(D):= K_{\pi} -D \simeq \sum_i\: a_i\cdot E_i \:,$$ 
with all $a_i\in \bZ$ satisfying the inequality $a_i>-1$ (i.e. $a_i\in \bN_0$). 
Here we use of course the relative canonical divisor $K_{\pi}$ in the ``$\Omega$-flavor".
Then we define the following {\em motivic arc invariant\/} 
$$\Cal E^{arc}((X,D))\in \widehat{\bM}(\Cal V/X)$$ 
of the pair $(X,D)$:
\begin{equation} \label{inv-arc}
\Cal E^{arc}((X,D)):= \pi_*\bigl( \int_{\Cal L(M)} \: \bL^{-ord(K_{\pi}(D))}d\tilde{\mu}_{M} \bigr) \:,
\end{equation}
which more explicitly can be calculated as in (\ref{eq:motint}).
This invariant is ``independent" of the choice of $\pi$ in the following sense. Let $\pi': M'\to M$
be a proper birational map of pure dimensional complex algebraic manifolds such that
$\pi'^*D$ and the exceptional locus of $\pi\circ\pi': M'\to X$ is contained in a normal crossing divisor with smooth irreducible components. Then 
$$K_{\pi\circ \pi'}(\pi'^*D)= K_{\pi\circ \pi'}-\pi'^*D = \pi'^*K_{\pi}(D)+K_{\pi'}$$ 
is also an effective Cartier divisor with 
$$\Cal E^{arc}((X,D))=\Cal E^{arc}((X,\pi'^*D))$$ 
by the transformation rule. So this is an invariant of the pair $(X,D)$, if we consider $D$ as a Cartier divisor 
(in the sense of Aluffi \cite{Aluffi-modification}) on the 
directed set of all such resolutions $\pi: M\to X$. 
In particular $\Cal E^{arc}(X):=\Cal E^{arc}((X,0))$ is an invariant of the singular space $X$.
In fact in the language of \cite[sec.6]{Denef-Loeser1} and \cite[sec.4.4]{Denef-Loeser2} it is just the ``motivic volume of the arc space $\Cal L(X)$" of the singular space $X$:
$$\Cal E^{arc}(X)= \int_{\Cal L(X)} \: 1 \:d\tilde{\mu}_{X} \:.$$
And this fits with our general description in the introduction of this section, if we set 
$$I(M):=[id_M]\in \widehat{\bM}(\Cal V/M)\:, \quad \text{with} \quad
J(\{E_i,a_i\}):= \int_{\Cal L(M)} \: \bL^{-ord(K_{\pi})}d\tilde{\mu}_{M}\:.$$

For the corresponding ``stringy invariant" in the ``$\omega$-flavor", one has first to extend these motivic integrals
to $\bQ$-Cartier divisors supported on a normal crossing divisor with smooth irreducible components $E_i$,
i.e. we start with a strict normal crossing divisor $E=\sum_{i=1}^{k}\:a_{i}E_{i}$ on the smooth manifold $M$,
with $a_{i}\in \bQ$ such that $r\cdot E$ is a Cartier divisor for some $r\in \bN$, i.e. $r\cdot a_i\in \bZ$
for all $i$. Add a formal variable $\bL^{1/r}$ to  $\widehat{\bM}(\Cal V)$ (and $const^*\bL^{1/r}$ to
$\widehat{\bM}(\Cal V/X)$),
with $(\bL^{1/r})^r= \bL$. Then one can introduce and evaluate the integral
\begin{equation} \label{motint2-str}
\int_{\Cal L(M)} \: \bL^{-ord(E)}d\tilde{\mu}_{M} :=
\sum_{p\in \bZ}\: \tilde{\mu}_{M}(\{ord(rE)=p\}) \cdot (\bL^{1/r})^{-p}  \:,
\end{equation}
with value in $\widehat{\bM}(\Cal V/M)[\bL^{1/r}]$,
if $a_i>-1$ for all $i$. Moreover the corresponding formula (\ref{eq:motint}) 
with $\bL^{a_i+1}:=(\bL^{1/r})^{r\cdot (a_i+1)}$, and transformation rule (\ref{transformation2})
are also true in this more general context (compare with \cite[Appendix]{Veys-zeta} for more details).\\
 
With these improvements, one can introduce for a ``Kawamata log terminal pair $(X,D)$"
the corresponding {\em motivic stringy invariant\/} (for a suitable $r\in \bN$):
$$\Cal E^{str}((X,D))\in \widehat{\bM}(\Cal V/X)[\bL^{1/r}]\:.$$
Let $D$ be a $\bQ$-Weil divisor on the normal and irreducible complex variety $X$ such that $K_X+D$ is a $\bQ$-Cartier divisor (with $r\cdot (K_X+D)$ a Cartier divisor)
satisfying the following condition: There is a resolution of  singularities $\pi: M\to X$ with the exceptional locus $E$
and the support of $K_{\pi}(D):=K_M-\pi^*( K_X+D)$ contained in a normal crossing divisor with smooth irreducible
components $E_i$ ($i\in I$) such that $K_{\pi}(D)\simeq \sum_i\: a_i\cdot E_i$, with all $a_i\in \bQ$ satisfying the inequality
$a_i>-1$. Then we set
\begin{equation} \label{inv-str}
\Cal E^{str}((X,D)):= \pi_*\bigl( \int_{\Cal L(M)} \: \bL^{-ord(K_{\pi}(D))}d\tilde{\mu}_{M} \bigr) \:,
\end{equation}
which more explicitly can be calculated as in (\ref{eq:motint}).
Once more this is an invariant of the pair $(X,D)$, not depending on the resolution $\pi$ by the transformation
rule! In the language of \cite{Denef-Loeser1, Denef-Loeser2, Denef-Loeser3} it is for $D=0$ just 
the ``motivic Gorenstein volume of the arc space $\Cal L(X)$" of the singular space $X$, i.e.
the following ``motivic integral" on the singular space $X$:
$$\Cal E^{str}((X))= \int_{\Cal L(X)} \: \bL^{-ord(K_{X})} \:d\tilde{\mu}_{X} \:.$$
Note that by our conventions $\Cal E^{str}((X,D))=\Cal E^{arc}((X,D))$ in case $D$ a Cartier divisor (with strict normal crossing)
on a smooth manifold $X=M$.

\subsection{Stringy/arc E-function and Euler characteristic}
By application of suitable transformations, one can build from the motivic invariants 
$\Cal E^{str}((X,D))$ and $\Cal E^{arc}((X,D))$ 
other invariants. For example by pushing down by a constant map:
$$const_*:  \widehat{\bM}(\Cal V/X)[\bL^{1/r}]\to \widehat{\bM}(\Cal V)[\bL^{1/r}] \:,$$
one can transform these ``relative invariants over $X$" to "absolute invariants"
(with $r=1$ in the case of ``arc invariants").
And then one can apply for example the ``E-function characteristic" 
$$E: \widehat{\bM}(\Cal V)[\bL^{1/r}]\to \bZ[u,v][[(uv)^{-1}]][(uv)^{1/r}] \:,$$
which is defined with the help of Deligne's mixed Hodge theory. Then 
$$E_{str}((X,D)):= E\bigl(\Cal E^{str}((X,D))\bigr)$$
becomes Batyrev's {\em stringy E-function\/} of the Kawamata log terminal pair $(X,D)$
(as in \cite{Batyrev-stringy}). Similarly
$$E_{arc}(X):= E\bigl(\Cal E^{arc}(X)\bigr)$$
is the ``Hodge-arc invariant" of $X$ in the sense of 
\cite[\S 6]{Denef-Loeser1} and
\cite[\S 4.4.1]{Denef-Loeser2} (up to a normalization factor $(uv)^{dim(X)}$ coming from
a different normalization of the motivic measure).\\ 

Here $E: K_{0}(\Cal V)\to \bZ[u,v] $ is induced from
\begin{equation} \label{eq:addgr}
X\mapsto E(X):= \sum_{i,p,q\geq 0}\:(-1)^{i}\cdot 
dim_{C}\left(gr_{F}^{p}gr_{p+q}^{W}H^{i}_{c}(X^{an},\bC)\right) u^{p}v^{q} \:,
\end{equation} 
with $F$ the decreasing Hodge filtration and $W$ the increasing weight filtration
of Deligne's canonical and functorial {\em mixed Hodge structure\/} on
$H^{i}_{c}(X^{an},\bQ)$ \cite{De-hodgeI, De-hodgeII}. Here $X^{an}$ means the complex algebraic variety $X$ with its classical (and not the Zariski) topology. 
This E-polynomial satisfies the defining ``additivity" relation of $K_{0}(\Cal V)$,
because the corresponding long exact cohomology sequence is strictly compatible with
the filtrations $F$ and $W$ (i.e. the sequence remains exact after application of
$gr_{F}^{p}gr_{p+q}^{W}$).

In particular, $E(-1,-1)(X)=\chi(X)$ is the topological Euler characteristic of $X$.
Finally classical Hodge theory implies, for $X$ {\em smooth and compact\/},
the ``purity result" $gr_{p+q}^{W}H^{i}(X^{an},\bC)=0$ for $p+q\neq i$,
together with
\begin{align*}
h^{p,q}(X):= &\sum_{i\geq 0}\: (-1)^{i}(-1)^{p+q}\cdot
dim_{C}\left(gr_{F}^{p}gr_{p+q}^{W}H^{i}_{c}(X^{an},\bC)\right)\\
= & dim_{C}\left(gr_{F}^{p}H^{p+q}(X^{an},\bC)\right)
= dim_{C}H^{q}(X^{an},\Lambda^{p}T^{*}X^{an})\\
= &  dim_{C}H^{q}(X,\Lambda^{p}T^{*}X)\:.
\end{align*}

\begin{rem}
One can get the transformation $E: K_{0}(\Cal V)\to \bZ[u,v]$
also as an application of Theorem \ref{blowup}
(but in a less explicit way), since the invariant
$$d_X:=E(X)= \sum_{p,q\geq 0}\:(-1)^{p+q}\cdot 
dim_{C}H^{q}(X,\Lambda^{p}T^{*}X) u^{p}v^{q}$$
for $X$ compact and smooth 
satisfies the corresponding properties (iii-1) and (iii-2).
\end{rem}

In particular, $\chi_{y}(X)= E(-y,1)(X)$ for $X$ smooth and compact
by ({\bf g-HRR}), so that this E-function is another generalization of the $\chi_{y}$-genus.
But the classes $[X]$  for $X$ smooth and compact generate $K_{0}(\Cal V)$ so
that we get the following Hodge theoretic description for any $X$
(with $T_{y*}$ our Hirzebruch class transformation of Theorem \ref{T_y_*}):
\begin{equation} \label{eq:TyH}
T_{y*}([X]) = \sum_{i,p\geq 0}\:(-1)^{i} 
dim_{C}\left(gr_{F}^{p}H^{i}_{c}(X^{an},\bC)\right) (-y)^{p} =  E(-y,1)(X)\:.
\end{equation}
Moreover $\chi_{y}(X):= E(-y,1)(X)$ is for $X\neq \emptyset$ of dimension $d$ a polynomial of degree $d$,
with $E(\bL)=E(\bC)=uv\in \bZ[u,v]$ so that one gets an induced map
$$E: \widehat{\bM}(\Cal V)[\bL^{1/r}]\to \bZ[u,v][[(uv)^{-1}]][(uv)^{1/r}] \:.$$
By (\ref{eq:motint}) we get the following explicit description of $E_{str}((X,D))$,
with $\pi: M\to X$ a resolution of singularities such that $K_{\pi}(D)\simeq \sum_i\: a_i\cdot E_i$ 
is a strict normal crossing divisor with $a_i>-1$ for all $i$ as before
(and similarly for $E_{arc}((X))$):
\begin{equation} \label{E-str} \begin{split}
 E_{str}((X,D)) &=
\sum_{I\subset \{1,\dots,k\}} \: E(E_{I}^{o}) \cdot \prod_{i\in I}\:
\frac{uv-1}{(uv)^{a_{i}+1}-1} \\
&= \sum_{I\subset \{1,\dots,k\}} \: E(E_{I}) \cdot \prod_{i\in I}\:
(\frac{uv-1}{(uv)^{a_{i}+1}-1} -1) \:.
\end{split} \end{equation}
Putting $(u,v)=(-y,1)$ gives a similar formula for (or defines) the ``stringy $\chi_y$-characteris\-tic" $\chi_y^{str}((X,D))$
(or the ``arc $\chi_y$-characteristic" $\chi_y^{arc}((X))$),
and also the limit $u, v\to 1$ exists with
 \begin{equation} \label{chi-str} \begin{split}
\chi^{str}((X,D)) &:=\lim_{u,v\to 1}\: E_{str}((X;D))\\ 
&= \sum_{I\subset \{1,\dots,k\}} \: \chi(E_{I}^{o}) \cdot \prod_{i\in I}\:\frac{1}{a_{i}+1} \\
&= \sum_{I\subset \{1,\dots,k\}} \: \chi(E_{I}) \cdot (-1)^{|I|}\cdot \prod_{i\in I}\:
\frac{a_i}{a_{i}+1} \:.
\end{split} \end{equation}
This $\chi^{str}((X,D))$ is just Batyrev's {\em stringy Euler 
number\/} of the log-terminal pair $(X,D)$
(as defined in \cite{Batyrev-stringy}). Similarly  $\chi^{arc}(X)$ is just the
{\em arc Euler characteristic\/} of $X$ in the sense of
\cite[\S 6]{Denef-Loeser1} and \cite[\S 4.4.1]{Denef-Loeser2}.
Finally note that (\ref{E-str}) and the ``limit formula" (\ref{limit-ell}) for the elliptic class
$\Cal E ll((X,D))$ of the pair $(X,D)$ imply for $X$ compact (with $y=e^{2\pi iz})$):
\begin{equation} 
\lim_{\tau\to i\infty}\: y^{dim(X)/2}\cdot ell((X,D))=
\chi_{-y}^{str}((X,D))= E_{str}((X,D))(y,1)\:.
\end{equation}

\subsection{Stringy and arc characteristic classes}\label{stringy-arc classes}
Recall our motivic characteristic class transformations $mC_*$ form Corollary \ref{mC*},
$T_{y*}$ from Theorem \ref{T_y_*} and $\tilde{T}_{y*}$ from Remark \ref{norm-ell}.
Here $T_{y,i}(\quad) = (1+y)^{-i}\cdot \tilde{T}_{y,i}(\quad)$ for all $i$,  so that both classes carry the same information.
These classes all satisfy $cl_*([\bC])=-y$,  so that they induce similar transformations on
$K_{0}(\Cal V/X)[\bL^{-1}]$:
\begin{equation*}\begin{split}
mC_{*} &: K_{0}(\Cal V/X)[\bL^{-1}]\to \bold G_{0}(X)\otimes \bZ[y,y^{-1}]\:,\\
T_{y*}, \tilde{T}_{y*}
&: K_{0}(\Cal V/X)[\bL^{-1}]\to H^{BM}_{*}(X)\otimes \bQ[y,y^{-1}]\:.\\
\end{split} \end{equation*}
And these extend by \cite[Corollary 2.1.1, Corollary 3.1.1]{BSY2} to the completions
\begin{equation} \label{eq:compl} \begin{split}
mC^{\wedge}_{*} &: \widehat{M}(\Cal V/X)[\bL^{1/r}]
\to \bold G_{0}(X)\otimes \bZ[y][[y^{-1}]][(-y)^{1/r}]\:,\\
T^{\wedge}_{y*}, \tilde{T}^{\wedge}_{y*} 
&: \widehat{M}(\Cal V/X)[\bL^{1/r}]
\to H^{BM}_{*}(X)\otimes \bQ[y][[y^{-1}]][(-y)^{1/r}]\:.\\
\end{split} \end{equation}
So we can introduce for $cl_{*}= mC_{*},T_{y*}, \tilde{T}_{y*}$
the corresponding
{\em stringy characteristic homology class\/} $cl_{*}^{str}((X,D))$ of the Kawamata log terminal pair $(X,D)$  by
\begin{equation} \label{eq:stringcl} \begin{split}
cl_{*}^{str}((X,D)):&=
cl_{*}^{\wedge}\bigl( \Cal E^{str}((X,D))\bigr) \:.
\end{split} \end{equation}
Moreover these transformations $cl^{\wedge}_{*}$ commute with proper push down and exterior products so that
$$cl^{\wedge}_{*}\bigl( f_*(\alpha \cdot const^*\beta)\bigr)=
\bigl(f_*(cl^{\wedge}_{*}(\alpha))\bigr)\cdot const^*\beta$$
for $f: X\to Y$ proper, with $\alpha \in \widehat{M}(\Cal V/X)$ and
$\beta\in \widehat{M}(\Cal V)$.
By (\ref{eq:motint}) we get the following explicit description of $cl_{*}^{str}((X,D))$,
with $\pi: M\to X$ a resolution of singularities such that 
$K_{\pi}(D)\simeq \sum_i\: a_i\cdot E_i$ 
is a strict normal crossing divisor with $a_i>-1$ for all $i$ as before:
\begin{equation} \label{eq:motint2} \begin{split}
 cl_{*}^{str}((X,D))=
\sum_{I\subset \{1,\dots,k\}} \: cl_{*}([E_{I}^{o}\to X]) \cdot \prod_{i\in I}\:
\frac{(-y)-1}{(-y)^{a_{i}+1}-1} \\
= \sum_{I\subset \{1,\dots,k\}} \: cl_{*}([E_{I}\to X]) \cdot \prod_{i\in I}\:
\frac{(-y)-(-y)^{a_{i}+1}}{(-y)^{a_{i}+1}-1}  \:.
\end{split} \end{equation}
But $E_{I}$ is a closed smooth submanifold of $M$ so that $cl_{*}([E_{I}\to X])$ is just the 
proper pushforward to $X$ of the corresponding characteristic (homology) class
$$cl_{*}(E_{I})=cl^{*}(TE_{I})\cap [E_{I}] 
\quad \text{for} \quad cl_{*}= mC_{*},T_{y*}, \tilde{T}_{y*}\:.$$

The {\em stringy Hirzebruch classes\/} $T^{str}_{y*}((X,D))$ and $\tilde{T}^{str}_{y*}((X,D))$
interpolate by (\ref{limit-ell}) and (\ref{eq:motint2}) in the following sense between the 
{\em elliptic class\/} $\Cal E ll_*((X,D))$ of Borisov-Libgober defined in (\ref{ell-BL}):
\begin{equation} \label{eq:ell=hirze}
\lim_{\tau\to i\infty} \: y^{dim(X)/2}\cdot \Cal E ll((X,D))(z,\tau) = 
\tilde{T}^{str}_{-y*}((X,D)) \quad \text{for $y=e^{2\pi iz}$} \:,
\end{equation}
and for compact $X$ the {\em stringy $E$-function\/} $E_{str}((X,D))$
of Batyrev as in (\ref{E-str}):
\begin{equation} \label{eq:Efct} \begin{split}
\chi_{-y}^{str}((X,D)):&= deg\bigl( T^{str}_{-y*}((X,D)) \bigr)\\
&= deg\bigl( \tilde{T}^{str}_{-y*}((X,D)) \bigr) = E_{str}((X,D))(y,1) \:.
\end{split} \end{equation}

So these stringy Hirzebruch classes are ``in between" the elliptic class
and the stringy $E$-function, and as suitable limits they are ``weaker"
than these more general invariants. But they have the following good properties of both of them:
\begin{itemize}
\item The stringy Hirzebruch classes come from a functorial ``additive"
characteristic homology class.
\item The stringy $E$-function comes from the ``additive" {\em E-polynomial\/} 
defined by Hodge theory, which does not have a homology class version
(compare with \cite[\S 5]{BSY2}).
\item The elliptic class is a homology class, which does not come
from an ``additive" characteristic class (of vector bundles), since the corresponding
{\em elliptic genus\/} is more general than the {\em Hirzebruch $\chi_{y}$-genus\/},
which is the most general ``additive" genus of such a class.
\end{itemize}
 
Finally the stringy Hirzebruch class $T^{str}_{y*}((X,D))$ specializes for $y=-1$
in the following way to the {\em stringy Chern class\/} $c_{*}^{str}((X,D))$ of $(X,D)$ as introduced in \cite{Aluffi-modification, FLNU}:
\begin{equation} \label{eq:hirze=ch}
\lim_{y\to -1}\:T^{str}_{y*}((X,D)) = c_{*}^{str}((X,D)) \in H^{BM}_{*}(X)\otimes\bQ \:.
\end{equation}
In fact
 \begin{equation} \label{c-str} \begin{split}
\lim_{y\to -1}\: T^{str}_{y*}((X,D))
&= \sum_{I\subset \{1,\dots,k\}} \: T_{-1*}([E_{I}^{o}\to X]) \cdot
 \prod_{i\in I}\:\frac{1}{a_{i}+1} \\
&= \sum_{I\subset \{1,\dots,k\}} \: T_{-1*}([E_{I}\to X]) \cdot 
(-1)^{|I|}\cdot \prod_{i\in I}\:\frac{a_i}{a_{i}+1} \:.
\end{split} \end{equation}
So by Theorem \ref{unification theorem} (for $y=-1$) we get:
\begin{equation} \label{c-str2} \begin{split}
\lim_{y\to -1}\: T^{str}_{y*}((X,D))
&= c_{*}\bigl(\sum_{I\subset \{1,\dots,k\}} \: \pi_*(1_{E_{I}^{o}}) \cdot
 \prod_{i\in I}\:\frac{1}{a_{i}+1} \bigr) \\
&= \sum_{I\subset \{1,\dots,k\}} \:  
(-1)^{|I|}\cdot \prod_{i\in I}\:\frac{a_i}{a_{i}+1} \cdot \pi_*(c_*(E_{I}))\:.
\end{split} \end{equation}

And the right hand side is just $c_{*}^{str}((X,D))$
by \cite[\S \S 3.4,5.5,6.5]{Aluffi-modification} and \cite[Corollary 2.5, \S 4]{FLNU}.
In a similar way one gets for $cl_{*}= mC_{*},T_{y*}, \tilde{T}_{y*}$
the {\em arc characteristic classes\/}
\begin{equation} \label{eq:arccl} \begin{split}
cl_{*}^{arc}((X,D)):&=
cl_{*}^{\wedge}\bigl( \Cal E^{arc}((X,D))\bigr) \:,
\end{split} \end{equation}
with
\begin{equation} \label{eq:hirze=ch2}
\lim_{y\to -1}\:T^{arc}_{y*}((X,D)) = c_{*}^{arc}((X,D)) \in H^{BM}_{*}(X)\otimes\bQ 
\end{equation}
the Chern class $\int_{X} \: \jeden(-D)\; dc_X$ of the pair $(X,-D)$ as
introduced and studied in \cite[\S \S 3.3,5.5]{Aluffi-modification},
with ``$\bL^{-ord(K_{\pi}(D))}$ corresponding to $\jeden(-D)$ for $\bL\to -y \to 1$".\\

Of course it is also natural to look at the other specializations $y\to 0$ and
$y\to 1$ of the stringy and arc characteristic classes $cl^{str/arc}_{*}((X,D))$ 
for $cl_{*}= mC_{*},T_{y*}, \tilde{T}_{y*}$. But the limit $y\to 1$
doesn't exist in general so that one can {\em not\/} introduce
``stringy or arc L-classes and signature" in this generality.  
But if we specialize in (\ref{eq:motint2}) for $D=0$ to $y=0$, then we get by ``additivity":
$$\lim_{y\to 0}\: mC^{str}_{*}(X) = \pi_{*}([\Cal O_{M}])= \lim_{y\to 0}\: mC^{arc}_{*}(X)$$
and
$$\lim_{y\to 0}\: T^{str}_{y*}(X) = \pi_{*}(Td^{*}(TM)\cap [M]) 
= \lim_{y\to 0}\: T^{arc}_{y*}(X) \:.$$
In particular the middle terms are independent of a resolution $\pi: M\to X$,
whose exceptional locus is contained in a strictly normal crossing divisor.
And by the ``weak factorization theorem" one can even conclude 
(compare \cite[Corollary 3.2]{BSY2}):
\begin{pro} \label{tdbirat}
Let $\pi: M\to X$ be a resolution of singularities of the pure dimensional
complex algebraic variety $X$. Then the classes 
$$ \pi_{*}([\Cal O_{M}])\in \bold G_0(X)
\quad \text{and} \quad
\pi_{*}(Td^{*}(TM)\cap [M]) \in H^{BM}_{*}(X)\otimes\bQ $$
are independent of $\pi$.
\end{pro}


\begin{thebibliography}{99999}

\bibitem[AKMW]{AKMW}
D. Abramovich, K. Karu, K. Matsuki and J. W\l odarczyk,
{\it Torification and factorization of birational maps},
J. Amer. Math. Soc., {\bf 15} (2002), 531--572.

\bibitem [Alu1]{Aluffi-lecture note}
P. Aluffi, {\it Characteristic classes of singular varieties},
Topics in Cohomological Studies of Algebraic Varieties (Ed. P. Pragacz), Trends in Mathematics, Birkh\"auser (2005), 1--32.

\bibitem [Alu2]{Aluffi-IMRN}
P. Aluffi, {\it Chern classes of birational varieties},
Int. Math. Res. Notices {\bf 63}(2004), 3367--3377.

\bibitem [Alu3]{Aluffi-limits}
P. Aluffi, {\it Limits of Chow groups, and a new construction of Chern-Schwartz-MacPherson classes},
math.AG/0507029.

\bibitem[Alu4]{Aluffi-modification}
P. Aluffi, {\it Modification systems and integration in their Chow groups},
arXiv:math.AG/0407150, to appear in Selecta Math.

\bibitem[AM]{Artin-Mazur} 
M. Artin and B. Mazur, {\it Etale Homotopy},
Springer Lecture Notes in Math. No. {\bf 100}, Springer-Verlag, Berlin, 1969.

\bibitem[AS]{AS}
M. F.  Atiyah, I. M. Singer, {\it The index of elliptic operators on compact manifolds},
Bull. Amer. Math. Soc. {\bf 69} (1963), 422-433. 

\bibitem [Ax]{Ax}
J. Ax, {\it Injective endomorphisms of varieties and schemes},
Pacific J. Math. {\bf  31} (1969), 1--17.

\bibitem[Ban]{Ban}
M. Banagl, {\it Extending intersection homology type invariants to non-Witt spaces},
Mem. Amer. Math. Soc. 160, 2002.

\bibitem[Bat1]{Batyrev-stringy}
V. Batyrev, {\it Non-Archimedean integrals and stringy Euler numbers of log-terminal pairs},
J. Eur. Math. Soc. {\bf 1} (1999, 5-33.  

\bibitem[Bat2]{Batyrev-birat}
V. Batyrev, {\it Birational Calabi-Yau n-folds have equal Betti numbers},
In {\it New trends in algebraic geometry}, London Math. Soc. Lecture Note Ser. {\bf 264} (1999),
1--11.

\bibitem[BFM1]{Baum-Fulton-MacPherson}
P. Baum, W. Fulton and R. MacPherson, {\it Riemann--Roch for singular varieties},
Publ. Math. I.H.E.S. {\bf 45} (1975), 101--145.

\bibitem[BFM2]{Baum-Fulton-MacPherson2}
P. Baum, W. Fulton and R. MacPherson, {\it Riemann-Roch and topological K-theory for singular varieties},
Acta Math. {\bf 143} (1979), 155-192. 

\bibitem [Bit]{Bittner}
F. Bittner, {\it The universal Euler characteristic for varieties of characteristic zero},
Compositio Math.  {\bf 140} (2004), 1011--1032.

\bibitem[BoSe]{Borel-Serre}
A. Borel and J.-P. Serre, {\it Le th\'eor\`eme de Riemann-Roch (d'apres Grothendieck)},
Bull.Soc.Math.France, {\bf 86} (1958), 97--136.

\bibitem[BM]{Borho-MacPherson}
W. Borho and R. MacPherson,
{\it Partial Resolutions of Nilpotent Varieties},
Ast\'erisque {\bf 101--102} (1983), 23--74.

\bibitem[BL1]{BL-elliptic}
L. Borisov, A. Libgober, {\it Elliptic genera for singular varieties},
Duke Math. J. {\bf 116} (2003), 319-351. 

\bibitem[BL2]{BL-McKay}
L. Borisov, A. Libgober, {\it  McKay correspondence for elliptic genera},
arXiv:math.AG/0206241.

\bibitem[Bor]{Borusk}
K. Borusk, {\it Theory of Shape},
Lecture Notes No. {\bf 28}, Matematisk Inst. Aarhus Univ. (1971), 1--145.

\bibitem[BF]{BF}
H.W. Braden, K.E. Feldman, {\it Functional equations and the generalized elliptic genus\/},
math-ph/0501011.

\bibitem[Br1]{Brasselet1} 
J.-P. Brasselet, {\it Existence des classes de Chern en th\'eorie bivariante}, Ast\'erisque,{\bf 101--102} (1981), 7--22.

\bibitem[Br2]{Brasselet-lecture note}
J.-P. Brasselet, {\it From Chern classes to Milnor classes -- a history of characteristic classes for singular varieties},
Singularities--Sapporo 1998, Adv. Stud. Pure Math. {\bf 29} (2000), 31--52

\bibitem[BLSS]{BLSS}
J.-P. Brasselet, D. Lehman, J. Seade and T. Suwa,
{\it Milnor classes of local complete intersection},
Tran. Amer. Math. Soc. {\bf 354} (2001), 1351--1371.

\bibitem[BSY1]{BSY1} 
J.-P. Brasselet, J. Sch\"urmann and S. Yokura, {\it Bivariant Chern classes and Grothendieck transformations},  math. AG/0404132.

\bibitem[BSY2]{BSY2} 
J.-P. Brasselet, J. Sch\"urmann and S. Yokura, {\it Hirzebruch classes and motivic Chern classes for singular spaces}, math. AG/0503492

\bibitem[BrSc]{Brasselet-Schwartz}
J.-P. Brasselet and M.-H. Schwartz, {\it Sur les classes de Chern d'une ensemble analytique complexe},
Ast\'erisque {\bf 82--83}(1981), 93--148.

\bibitem[BZ]{BZ}
J.-L. Brylinski, B. Zhang, {\it Equivariant Todd Classes for Toric Varieties},
math. AG/0311318.

\bibitem[CS1]{Cappell-Shaneson1}
S. Cappell and J. L. Shaneson, {\it Stratifiable maps and topological invariants},
J. Amer. Math. Soc. {\bf 4} (1991),  521--551.

\bibitem[CS2]{Cappell-Shaneson2}
S. Cappell and J. L. Shaneson, {\it Genera of algebraic varieties and counting lattice points},
Bull. Amer. Math. Soc. {\bf 30} (1994),  62--69.

\bibitem[Ch1]{Chern1}
S. S. Chern, {\it Characteristic classes of Hermitian manifolds},
Ann. Math. {\bf 47} (1946),  85--121.

\bibitem[Ch2]{Chern2}
S. S. Chern, {\it On the multiplication in the characteristic ring of a sphere bundle},
Ann. Math. {\bf 49} (1948),  362--372.

\bibitem[Ch3]{Chern3}
S. S. Chern, {\it A simple intrinsic proof of the Gauss--Bonnet formula for closed Riemann manifolds},
Ann. of Math. {\bf  45} (1944), 747--752.

\bibitem[CG]{Chriss-Ginzburg}
N. Chriss and V. Ginzburg,
{\it Representation theory and complex geometry},
Birkh\"auser, 1997.

\bibitem [Cr]{Craw}
A. Craw, {\it An introduction to motivic integration},
in {\it Strings and Geometry}, Clay Math. Proc., {\bf 3} (2004), Amer. Math. Soc., 203--225.

\bibitem [DL1]{Denef-Loeser1}
J. Denef and F. Loeser, {\it Germs of arcs on singular algebraic varieties and motivic integration},
Invent. Math. {\bf 135} (1999), 201--232

\bibitem [DL2]{Denef-Loeser2}
J. Denef and F. Loeser, {\it Geometry on arc spaces of algebraic varieties},
European Congress of Mathematicians (Barcelona, 2000), {\bf 1} (2001) Birkh\"auser, 327--348 
\bibitem [DL3]{Denef-Loeser3}
J. Denef and F. Loeser, {\it Motivic integration, quotient singularitis and the McKay correspondence},
Comp. Math. {\bf 131} (2002), 267--290.

%\bibitem [Dim]{Dimca}
%A. Dimca, {\it Sheaves in Topology},
%Springer-Verlag, 2004.

\bibitem[EG1]{EG1}
D. Edidin, W. Graham, {\it Riemann-Roch for equivariant Chow groups},
Duke Math. J. {\bf 102} (2000), 567--594.

\bibitem[EG2]{EG2}
D. Edidin, W. Graham, {\it Riemann-Roch for quotients and Todd classes of simplicial toric varieties},
Commun. Algebra {\bf 31}, No.8 (2003), 3735--3752.

\bibitem[De1]{De-hodgeI}
P. Deligne, {\it Th\'{e}orie des Hodge II},
Publ. Math. IHES {\bf 40} (1971), 5-58. 

\bibitem[De2]{De-hodgeII}
P. Deligne, {\it  Th\'{e}orie des Hodge III},
Publ. Math. IHES {\bf 44} (1974), 5-78. 

\bibitem[DuBo]{DuBois}
Ph. DuBois,
{\it Complexe de De Rham filtr\'e d'une vari\'et\'e singuli\`ere},
Bull. Soc. Math. France, {\bf 109} (1981), 41--81.

\bibitem[Er]{Ernstrom}
L. Ernstr\"om, {\it Topological Radon transforms and the local Euler obstruction}, Duke Math. J. {\bf 76} (1994), 1--21

\bibitem[FLNU]{FLNU}
T. de Fernex, E. Lupercio, T. Nevins, B. Uribe, {\it Stringy Chern classes of singular varieties},
arXiv:math.AG/0407314.

\bibitem [Fu1] {Fulton-intersection} 
W. Fulton, {\it  Intersection Theory},
Springer Verlag, 1981.

\bibitem [Fu2] {Fulton-toric} 
W. Fulton, {\it Introduction to Toric Varieties},
Ann. of Math. Studies, No. {\bf 131}, Princeton Univ. Press, 1993.
\bibitem[FM]{Fulton-Johnson}
W. Fulton and K. Johnson, {\it Canonical classes on singular varieties},
Manus. Math. {\bf 32} (1980), 381--389.

\bibitem[FM]{Fulton-MacPherson}
W. Fulton and R. MacPherson, {\it Categorical frameworks for the study of singular spaces},
Memoirs of Amer. Math. Soc. {\bf 243}, 1981.

\bibitem[FuMC]{FuMC}
J. Fu and C. McCrory, {\it Stiefel-Whitney classes and the conormal cycle of a singular variety},
Trans. Amer. Math. Soc. {\bf 349} (1997), 809--835  

\bibitem[Gi1]{Ginzburg-1}
V. Ginzburg, {\it G--Modules, Springer's Representations and Bivariant Chern Classes}, Adv. in Math., {\bf  61} (1986), 1--48.

\bibitem[Gi2]{Ginzburg-2}
V. Ginzburg, {\it Geometric methods in the representation theory of Hecke algebras and quantumgroups}, In: A. Broer and A. Daigneault (Eds.): {\it Representation theories and algebraic geometry (Montreal, PQ, 1997)}, Kluwer Acad. Publ., Dordrecht, 1998, 127--183.

\bibitem[GM1]{Goresky-MacPherson1}
M. Goresky and R. MacPherson, {\it Intersection homology theory},
Topology {\bf 149} (1980), 155--162.

\bibitem[GM2]{Goresky-MacPherson2}
M. Goresky and R. MacPherson, {\it Intersection homology, II},
Inventiones Math. {\bf 149} (1980), 77--129.

%\bibitem[Got]{Gottschalk}
%W. Gottschalk, {\it Some general dynamical notions}, 
%Springer LNM.  No. 318 (1973), 120--125.

\bibitem[Grom]{Gromov1}
M. Gromov, {\it Endomorphisms of symbolic algebraic varieties},
J. Eur. Math. Soc. {\bf 1} (1999), 109--197.

%\bibitem[Grom2]{Gromov2}
%M. Gromov, {\it Topological invariants of dynamical systems and spaces of holomorphic maps: I}, 
%Mathematical Physics, Analysis and Geometry {\bf 2} (1999), 323--415.

\bibitem[Gros]{Gros}
M. Gros, {\it Classes de Chern et classes de cycles en cohomologie de Hodge-Witt logarithmique},
Bull. Soc. Math. France Mem. 21, 1985.

\bibitem[Grot1]{Grothendieck1}
A. Grothendieck, {\it Technique de descente et th\'eor\`emes d'existence en g\'eom\'etrie alg\'ebrique, II},
S\'eminaire Bourbaki, 12 \`eme ann\'ee, expos\'e 190-195 (1959-60).

\bibitem[Grot2]{Grothendieck2}
A. Grothendieck, {\it R\'ecoltes er Semailles-- R\'eflexions et T\'emoignages sur un pass\'e de math\'ematicien}, Reprint, 1985.

\bibitem[GNA]{GNA}
F. Guill\'en, V. Navarro Aznar, {\it Un crit\`ere d'extension des foncteurs d\'efinis
sur les sch\'emas lisses}, Publ. Math. I.H.E.S. {\bf 95}  (2002), 1--91.

\bibitem[H\"{o}hn]{Hoehn}
G. H\"{o}hn, {\it Komplexe elliptische Geschlechter und $S^1$-\"{a}quivariante Kobordismustheorie},
Diplomarbeit, Bonn 1991.

\bibitem[Hi]{Hironaka}
H. Hironaka, {\it Resolution of singularities of an algebraic variety over a field of characteristic zero}, Ann. of Math. {\bf 79} (1964), 109--326.

\bibitem[Hir1]{Hirzebruch1}
F. Hirzebruch, {\it Some problems on differentiable and complex manifolds},
Ann. of Math. {\bf 60} (1954), 213--236.

\bibitem[Hir2]{Hirzebruch2}
F. Hirzebruch, {\it Topological Methods in Algebraic Geometry, 3rd ed. (1st German ed. 1956)},
Springer-Verlag, 1966.

\bibitem[Hir3]{Hirzebruch3}
F. Hirzebruch, {\it The Signature Theorem: Reminiscences and Recreation}, 
in ``Prospects in Mathematics", Ann. of Math. Studies, No.{\bf 70}, Princeton Univ. Press, 1971.

\bibitem[HBJ]{HBJ}
F. Hirzebruch, T. Berger and R. Jung, {\it Manifolds and Modular Forms},
Vieweg, 1992.

\bibitem[Hus]{Hus}
D. Husemoller, {\it Fibre Bundles, 3rd ed.},
Springer-Verlag, 1994.

\bibitem [KS]{Kashiwara-Schapira}
M. Kashiwara and P. Schapira, {\it Sheaves on Manifolds},
Springer-Verlag, Berlin, Heidelberg, 1990.
\bibitem[Ken]{Ken}
G. Kennedy, {\it MacPherson's Chern classes of singular varieties},
Com. Algebra. {\bf 9} (1990), 2821--2839. 

\bibitem[KMY]{KMY}
G. Kennedy, C. McCrory and S. Yokura, {\it Natural transformations from constructible functions to homology},
C. R. Acad. Sci. Paris, {\bf  t.319} (1994), 969--973.

\bibitem[KM]{Kollar-Mori}
J. Koll\'{a}r, S. Mori, {\it Birational geometry of algebraic varieties},
Cambridge Tracts in Math. 124, Cambridge Univ. Press, 1998. 

\bibitem[Kon]{Kontsevich} 
M. Kontsevich, {\it Lecture at Orsay}, 
1995.

\bibitem[Krich]{Krich} 
I. Krichever, {\it Generalized elliptic genera and Baker-Akhiezer functions},
Math. Notes {\bf 47} (1990), 132--142.

\bibitem[Kw1]{Kwiecinski1} 
M. Kwieci\'nski, {\it Formule du produit pour les classes caract\'eristiquesdeChern-Schwartz-MacPherson et homologie d'intersection}, 
C. R. Acad. Sci. Paris, {\bf 314} (1992),  625--628.

\bibitem[Kw2]{Kwiecinski2} 
M. Kwieci\'nski, {\it Sur le transform\'e de Nash et la construction du graph de MacPherson}, 
In Th\`ese, Universit\'e de Provence, 1994.

\bibitem[KY]{KY} 
M. Kwieci\'nski and S. Yokura, {\it Product formula of the twisted MacPherson class}, 
Proc. Japan Acad., {\bf 68} (1992), 167--171.
\bibitem[Levy]{Levy}
R. Levy, {\it Riemann-Roch theorem for complex spaces},
Acta Math. {\bf 158} (1987), 149-188.

\bibitem[Lo]{Looijenga}
E. Looijenga, {\it Motivic measures},
S\'eminaire Bourbaki, Ast\'erisque {\bf 276} (2002), 267--297.

\bibitem[Mac1]{MacPherson1}
R. MacPherson, {\it Chern classes for singular algebraic varieties},
Ann. of Math. {\bf 100}(1974), 423--432.

\bibitem[Mac2]{MacPherson2}
R. MacPherson, {\it Characteristic classes for singular  varieties},
Proceedings of the 9-th Brazilian Mathematical Colloquium (Pocos de Caldas, 1973) Vol.II, Instituto de Matem\'atica Pura e Aplicada, 
S\~ao Paulo, (1977), 321--327.

\bibitem[Man]{Manin}
Yu. I. Manin, {\it Lectures on the K-functor in algebraic geometry},
Russ. Math. Survey, {\bf 24} (1969), 1--89.

\bibitem[MaSe]{Mardesic-Segal}
S. Mardesi\'c and J. Segal, {\it Shape Theory},
North-Holland, 1982.

\bibitem[Mi1]{Milnor}
J. W. Milnor , {\it Topology from the differentiable Viewpoint},
Univ. Virginia Press 1965.

\bibitem[Mi2]{Milnor-hypersurface}
J. W. Milnor , {\it Singular points of complex hypersurfaces},
Ann. of Math. Studies, No. {\bf 61}, Princeton Univ. Press, 1968.

\bibitem[MiSt]{Milnor-Stasheff}
J. W. Milnor and J. D. Stasheff, {\it Characteristic classes},
Ann. of Math. Studies {\bf  76}, Princeton Univ. Press 1974.

\bibitem[Nov]{Nov}
S. P. Novikov, {\it Topological invariance of rational classes of Pontrjagin},
Dokl. Akad. Nauk SSSR {\bf 163} (1965), 298--300,
English translation, Soviet Math. Dokl. {\bf 6} (1965), 921--923. 


\bibitem[Oh]{Oh}
T. Ohmoto, {\it Equivariant Chern classes of singular algebraic varieties with group actions},
math.AG/0407348.

\bibitem[Pa]{Parusinski-lecture note}
A. Parusi\'nski, {\it Characteristic classes of singular varieties},
Singularity Theory and Its Applications, Sapporo, September 16--25, 2003.

\bibitem[Pontr]{Pontr}
L. Pontrjagin, {\it Characteristic cycles on differentiable manifolds},
Mat. Sbornik N. S. {\bf 21}, (1947), 233--284, A.M.S. Translation {\bf 32} (1950).

\bibitem[Sa]{Sabbah} 
C. Sabbah, {\it Espaces conormaux bivariants}, 
Th\`ese, l'Universit\'e Paris, 1986.

\bibitem[Sai]{Saito}
M. Saito, {\it Mixed Hodge Modules},
Publ. RIMS., Kyoto Univ. {\bf 26} (1990), 221--333.

\bibitem[Sch1]{Schuermann-VRR}
J. Sch\"urmann, {\it A generalized Verdier-type Riemann--Roch theorem for Chern-- Schwartz-- MacPherson classes},
math.AG/0202175.

\bibitem[Sch2]{Schuermann-partial GT} 
J. Sch\"urmann, {\it A general construction of partial Grothendieck transformations}, 
math. AG/0209299.

\bibitem[Sch3]{Schuermann-book}
J. Sch\"urmann, {\it Topology of singular spaces and constructible sheaves},
Monografie Matematyczne {\bf 63} (New Series), Birkh\"auser, Basel, 2003.

\bibitem [Sch4]{Schuermann-lecture note}
J. Sch\"urmann, {\it Lecture on characteristic classes of constructible functions},
Topics in Cohomological Studies of Algebraic Varieties (Ed. P. Pragacz), Trends in Mathematics, Birkh\"auser (2005), 175--201.


\bibitem[Sch5]{Schuermann-index}
J. Sch\"urmann, {\it Characteristic cycles and indices of $1$-forms on singular spaces},
in preparation.

\bibitem[Schw1]{Schwartz1}
M.-H. Schwartz, {\it  Classes caract{\'e}ristiques d{\'e}finies par unestratiification d'une vari{\'e}ti{\'e} analytique complex},
C. R. Acad. Sci. Paris {\bf t. 260} (1965), 3262-3264, 3535--3537.

\bibitem[Schw2]{Schwartz2}
M.-H. Schwartz, {\it Classes et caract\`eres de Chern des espaces lin\'eaires},  
Pub. Int. Univ. Lille, {\bf 2 Fasc. 3}(1980).

\bibitem[Serre]{Serre-GAGA}
J.-P. Serre, {\it G\'eom\'etrie alg\'ebrique et g\'eom\'etrie analytique},
Ann. Inst. Fourier {\bf 6} (1956), 1--42.

\bibitem[Sh]{Shaneson}
J. L. Shaneson, {\it Characteristic classes, lattice points and Euler--Maclaurin formulae},  
Proceedings of the ICM'94 (Z\"urich, Switerland) Birkh\"auser Verlag (1995), 612--624.

\bibitem[Si]{Si}
P. H. Siegel, {\it Witt spaces: A geometric cycle theory for KO-homology at odd primes},
Amer. J. of Math. {\bf 105}, (1983) 1067-1105.

\bibitem[Sti]{Stiefel}
E. Stiefel, {\it Richtungsfelder und Fernparallelismus in Mannigfaltigkeiten},
Comm. Math. Helv. {\bf 8} (1936), 3--51.

\bibitem[Stong]{Stong}
R. E. Stong, {\it Notes on Cobordism Theory},
Princeton Math. Notes, Princeton Univ. Press 1968.

\bibitem[Sull]{Sullivan}
D. Sullivan, {\it Combinatorial invariants of analytic spaces},
Springer Lecture Notes in Math., {\bf 192} (1970), 165--168.

\bibitem[Su]{Suwa-lecture note}
T. Suwa, {\it Characteristic classes of singular varieties},
Sugaku Expositions, Amer. Math. Soc. {\bf 16} (2003), 153--175.

\bibitem[Thom1]{Thom1}
R. Thom, {\it Espaces fibr\'{e}s en spheres et carr\'{e}s de Steenrod},
Ann. Sci. Ecole Norm. Sup. {\bf 69} (1952), 109--181. 

\bibitem[Thom2]{Thom}
R. Thom, {\it Les Classes Caracteristique de Pontrjagin des Vari\'et\'es Triangul\'es},
Symp. Intern. de Topologia Algebraica. Unesco (1958).

\bibitem[To]{Totaro}
B. Totaro, {\it Chern numbers for singular varieties and elliptic homology},
Ann. Math. {\bf 151} (2000), 757-791. 

\bibitem[Ve1]{Veys-zeta}
W. Veys, {\it Zeta functions and `Kontsevich invariants' on singular varieties},
Can. J. Math. {\bf 53} (2001), 834--865.

\bibitem[Ve2]{Veys-stringy}
W. Veys, {\it Arc spaces, motivic integration and stringy invariants},
Singularity Theory and Its Applications, Sapporo, September 16--25, 2003, math.AG/0401374.

\bibitem[Wang]{Wang}
C.-L. Wang, {\it K-equivalence in birational geometry and characterization of complex elliptic genera},
J. Alg. Geometry {\bf 12} (2003), 285--306.

\bibitem[Wh1]{Whitney}
H. Whitney, {\it Sphere spaces},
Proc. Nat. Acad. Sci. {\bf 21} (1935), 462--468.

\bibitem[Wh2]{Whitney2}
H. Whitney, {\it On the theory of sphere bundles},
Proc. Nat. Acad. Sci. {\bf 26} (1940), 148--153.

\bibitem[Wi]{Witten}
E. Witten, {\it Supersymmetry and Morse theory},
J. Diff. Geom. {\bf 17} (1982), 661--692.

\bibitem[W]{Wlodarczyk}
J. W\l odarczyk, {\it Toroidal varieties and the weak factorization theorem},
Invent. Math., {\bf 154} (2003), 223--331.

\bibitem[Wo]{Woolf}
J. Woolf, {\it Witt groups of sheaves on topological spaces},
math. AT/0510196.

\bibitem [Y1] {Yokura-RIMS}
S. Yokura, {\it A generalized Grothendieck--Riemann--Roch theorem for Hirzebruch's $\chi_y$-characteristic and $T_y$-characteristic},
Publ. Res. Inst. Math. Sci., {\bf 30} (1994), 603--610.

\bibitem [Y2] {Yokura-TAMS}
S. Yokura, {\it On Cappell--Shaneson's homology $L$-class of singular algebraic varieties},
Trans. Amer. Math. Soc., {\bf 347} (1995), 1005--1012.

\bibitem [Y3] {Yokura-Banach Center}
S. Yokura, {\it A singular Riemann--Roch theorem for Hirzebruch characteristics},
Banach Center Publ., {\bf 44} (1998), 257--268.

\bibitem [Y4] {Yokura-VRR-Chern}
S. Yokura, {\it On a Verdier-type Riemann--Roch for Chern--Schwartz--MacPherson class},
Topology and Its Applications., {\bf 94} (1999), 315--327.

\bibitem[Y5] {Yokura-Documenta} 
S. Yokura, {\it On the uniqueness problem of the bivariant Chern classes},
Documenta Mathematica, {\bf 7} (2002), 133--142.

\bibitem[Y6]{Yokura-TopAppl} 
S. Yokura:, {\it Bivariant theories of constructible functions and Grothendieck transformations}, 
Topology and Its Applications, {\bf 123} (2002), 283--296.

\bibitem[Y7]{Yokura-TAMS2} 
S. Yokura:, {\it On Ginzburg's bivariant Chern classes}, 
Trans. Amer. Math. Soc., {\bf 355} (2003), 2501--2521.

\bibitem[Y8]{Yokura-Dedicata} 
S. Yokura:, {\it On Ginzburg's bivariant Chern classes,II}, 
Geometriae Dedicata, {\bf 101} (2003), 185--201.

\bibitem[Y9]{Yokura-CEJM} 
S. Yokura, {\it Bivariant Chern classes for morphisms with nonsingular target varieties}, 
Central European J. Math., {\bf 3} (2005), 614--626.

\bibitem [Y10] {Yokura-Chern motivic}
S. Yokura, {\it Chern classes of proalgebraic varieties and motivic measures},
math.AG/0407237.

\bibitem [Y11] {Yokura-characteristic motivic} 
S. Yokura, {\it Characteristic classes of proalgebraic varieties and motivic measures},
preprint (2005).

\bibitem[You]{You} 
B. Youssin, {\it Witt Groups of Derived Categories},
K-Theory,  {\bf 11} (1997), 373-395.

\bibitem[Zh]{Zhang}
W. Zhang, {\it Lecture on Chern--Weil Theory and Witten Deformations},
Nankai Tracts in Mathematics Vol.4, World Scientific, 2001.

\bibitem[Z1]{Zhou1} 
J. Zhou, {\it Classes de Chern en th\'eorie bivariante}, 
Th\`ese, Universit\'e Aix-Marseille, 1995.

\bibitem[Z2]{Zhou2} 
J. Zhou, {\it Morphisme cellulaire et classes de Chern bivariantes},
Ann. Fac. Sci. Toulouse Math., {\bf 9} (2000), 161--192.

\end{thebibliography}
\end{document}